\documentclass[mnsc,nonblindrev]{informs3Mod}

\OneAndAHalfSpacedXI

\usepackage{subfig}
\usepackage{caption}
\usepackage{amsmath, amssymb,amstext}
\usepackage{bm}
\usepackage{longtable}
\usepackage{algorithm}
\usepackage{algorithmic}
\usepackage{enumerate}
\usepackage{enumitem}
\usepackage{appendix}
\usepackage{secdot}
\usepackage{calc}
\usepackage{bbm}
\usepackage{multirow}
\usepackage{color}
\usepackage{xr, xr-hyper}
\usepackage{hyperref}
\hypersetup{
	colorlinks=true,
	citecolor=blue,
	linkcolor=black,
	urlcolor=blue
}

\newcommand{\E}{\mathbb{E}}

\newcommand{\Pc}{\mathcal{P}}

\newcommand{\bx}{{\mathbf x}}
\newcommand{\bz}{{\mathbf z}}
\newcommand{\by}{{\mathbf y}}

\newcommand{\bw}{{\mathbf w}}
\newcommand{\bv}{{\mathbf v}}
\newcommand{\ba}{{\mathbf a}}
\newcommand{\bs}{{\mathbf s}}

\newcommand{\bxi}{{\bm \xi}}
\newcommand{\bzeta}{{\bm \zeta}}

\newcommand{\vertiii}[1]{{\left\vert\kern-0.25ex\left\vert\kern-0.25ex\left\vert #1
		\right\vert\kern-0.25ex\right\vert\kern-0.25ex\right\vert}}
\newcommand{\R}{\mathbb{R}}

\newcommand{\alpsoln}{\ensuremath{\theta}}
\newcommand{\Space}{\mathcal{S}}
\newcommand{\action}{\mathcal{A}}
\newcommand{\epsoracle}{\epsilon_{\mathcal{A}}}

\usepackage{natbib}
 \bibpunct[, ]{(}{)}{,}{a}{}{,}%
 \def\bibfont{\small}%

\TheoremsNumberedThrough     
\ECRepeatTheorems

\EquationsNumberedThrough    

\MANUSCRIPTNO{MS-17-01263.R2}
\begin{document}

\RUNAUTHOR{}

\RUNTITLE{Stochastic Level-Set Method}

\TITLE{A Data Efficient and Feasible Level Set Method for Stochastic Convex Optimization with Expectation Constraints}

\ARTICLEAUTHORS{%
\AUTHOR{Qihang Lin}
\AFF{Tippie College of Business, The University of Iowa, 21 East Market Street, Iowa City, IA 52242, USA, \EMAIL{qihang-lin@uiowa.edu}}
\AUTHOR{Selvaprabu Nadarajah}
\AFF{College of Business Administration,  University of Illinois at Chicago, 601 South Morgan Street, Chicago, Illinois, 60607, USA, \EMAIL{selvan@uic.edu}}
\AUTHOR{Negar Soheili}
\AFF{College of Business Administration,  University of Illinois at Chicago, 601 South Morgan Street, Chicago, Illinois, 60607, USA, \EMAIL{nazad@uic.edu}}
\AUTHOR{Tianbao Yang}
\AFF{CDepartment of Computer Science, The University of Iowa, 21 East Market Street, Iowa City, IA 52242, USA, \EMAIL{tianbao-yang@uiowa.edu }}
} 

\ABSTRACT{Stochastic convex optimization problems with expectation constraints (SOECs) are encountered in statistics and machine learning, business, and engineering. In data-rich environments, the SOEC objective and constraints contain expectations defined with respect to large datasets. Therefore, efficient algorithms for solving such SOECs need to limit the fraction of data points that they use, which we refer to as algorithmic data complexity. Recent stochastic first order methods exhibit low data complexity when handling SOECs but guarantee near-feasibility and near-optimality only at convergence. These methods may thus return highly infeasible solutions when heuristically terminated, as is often the case, due to theoretical convergence criteria being highly conservative. This issue limits the use of first order methods in several applications where the SOEC constraints encode implementation requirements. We design a stochastic feasible level set method (SFLS) for SOECs that has low data complexity and emphasizes feasibility before convergence. Specifically, our level-set method solves a root-finding problem by calling a novel first order oracle that computes a stochastic upper bound on the level-set function by extending mirror descent and online validation techniques. We establish that SFLS maintains a high-probability feasible solution at each root-finding iteration and exhibits favorable iteration complexity compared to state-of-the-art deterministic feasible level set and stochastic subgradient methods. Numerical experiments on three diverse applications validate the low data complexity of SFLS relative to the former approach and highlight how SFLS finds feasible solutions with small optimality gaps significantly faster than the latter method.\looseness=-1
}

\maketitle
\vspace{-0.4in}

\section{Introduction}\label{sec:intro}

Consider the stochastic optimization problem with expectation constraints (SOEC) 
\begin{eqnarray}
\label{eq:gco}
f^*:=
\min_{\bx\in\mathcal{X}}\left\{f_0(\bx)=\E\left[F_0(\bx,\xi_0)\right]\right\}\quad\text{s.t.}\quad f_i(\bx):=\E\left[F_i(\bx,\xi_i)\right]\leq r_i,~ i=1,2,\dots,m,
\end{eqnarray}
where $\mathcal{X}\subset\mathbb{R}^d$ is a nonempty closed convex set, $\xi_i$, $i=0,1,\ldots,m$, is a random vector whose probability distribution is supported on set $\Xi_i\subseteq \mathbb{R}^{q_i}$, and $F_i(\bx,\xi_i):\mathcal{X}\times\Xi_i\rightarrow \mathbb{R}$ is  continuous and convex in $\bx$ for each realization of $\xi_i$ for $i=0,1,2,\dots,m$. Given $\epsilon >0$, a solution $\bx_\epsilon\in\mathcal{X}$ is called \emph{$\epsilon$-feasible} if $\max_{i=1,\dots,m}\{f_i(\bx_\epsilon)-r_i\}\leq \epsilon$. A solution $\bx_\epsilon\in\mathcal{X}$ is referred to as \emph{$\epsilon$-optimal} if $f_0(\bx_\epsilon)-f^*\leq \epsilon$. Alternatively, optimality can be measured relative to an initial feasible solution  $\bx^0 \in\mathcal{X}$. In this case, we say $\bx_\epsilon\in\mathcal{X}$ is \emph{relative $\epsilon$-optimal} with respect to $\bx^0$ if $(f(\bx_\epsilon) - f^*)/(f(\bx^0) - f^*)  \leq \epsilon$.

Problem \eqref{eq:gco} is pervasive in stochastic optimization and appears as a central challenge in semi-supervised learning~\citep{chapelle2009semi}, shape-restricted regression~\citep{seijo2011nonparametric,sen2017testing,lim2014convergence,Cotter2016,fard2016fast},  Neyman-Pearson classification~\citep{tong2016survey,rigollet2011neyman,tong2013plug,zhao2015neyman}, approximate linear programming and related relaxations~\citep{de_farias_linear_2003,adelman2013dynamic,nadarajah2015relaxations}, portfolio selection~\citep{markowitz1952portfolio,abdelaziz2007multi}, risk management~\citep{rockafellar2000optimization}, supply chain design \citep{azaron2008multi}, and multi-objective stochastic programming~\citep{marler2004survey,abdelaziz2012solution,mahdavi2013stochastic,barba2017multi}. In this paper, we focus on overcoming the challenges of applying existing methods for solving SOECs in settings that are both data rich and where expectation constraints capture requirements that cannot be violated during real-world implementation.

In data-rich environments, each expectation appearing in \eqref{eq:gco} is defined by a data set containing a large number of data points (possibly infinite). The number of data points used when solving SOEC is an important computational bottleneck, which we refer to as the data complexity of an algorithm. Traditional approaches for solving SOECs can lead to large data complexity. For instance, consider the popular strategy of replacing each expectation in \eqref{eq:gco} by a sample average approximation (SAA; \citealp{shapiro2013sample,oliveira2017sample}) and solving the resulting model using a deterministic iterative method (see, e.g., \citealp{nesterov2004introductory,soheili2012smooth}, and references therein). If the number of samples used to construct SAAs is small, the solution from the deterministic approximation may be highly infeasible to the original SOEC, in addition to being suboptimal ~\citep{shapiro2013sample,oliveira2017sample}. Instead, if a large number of samples are used in each SAA, then the data complexity becomes large because the gradient or objective function evaluation at each iteration requires using a significant portion of each of the data sets.

In contrast, stochastic first order methods for tackling stochastic optimization problems have low per-iteration cost and data complexity and thus play a central role in machine learning packages such as TensorFlow and PyTorch ~\citep{RobMon51,nemirov09,Lan:2010:Optimal,ghadimi2012,DBLP:journals/siamjo/Lan13b,DBLP:conf/nips/ChenLP12,lan2012validation,Schmidt13,Pegasos,LanZ15,DBLP:journals/oms/LinCP14,duchi-2009-efficient,xiao2014proximal,xiao2010dual,DBLP:journals/jmlr/HazanK11a,DBLP:conf/nips/BachM13,Allen-Zhu:2017:STOC,goldfarb2017linear}.
These methods update solutions using stochastic gradients that can be computed using a small number of sampled data points. Stochastic first order methods typically ensure feasibility via projections onto a convex set at each iteration, where the convex set is assumed to be simple (e.g. a box or ball) for computational tractability. This assumption limits the applicability of first order methods for solving SOECs with general non-linear constraints. Recently, \citet{lan2016algorithms} and \citet{yu2017online} developed stochastic subgradient (SSG) methods devoid of projections for solving \eqref{eq:gco} with single ($m=1$) and multiple constraints ($m > 1$), respectively. The SSG methods in these papers guarantee an $\epsilon$-optimal and $\epsilon$-feasible solution only at convergence. 

In practice, SSG methods are terminated before their conservative theoretical conditions are met. Premature termination may lead to highly infeasible and sub- or super- optimal solutions. While some deviation from optimality is likely acceptable, a highly infeasible solution may not be implementable. Such situations arise in several data science applications in machine learning, as well as, across business (e.g., operations and finance) and engineering domains. We elaborate on the practical need for feasibility in a few cases below.
\begin{itemize}
	\item Fairness constraints: Enforcing fairness criteria when learning classifiers across multiple classes (e.g., male and female) has become important in machine learning \citep{goh2016satisfying}. This learning problem can be cast as an SOEC where fairness is modeled via expectation constraints. Constraint violations lead to classifiers that are biased towards one or more classes.
	\item Risk constraints: Planning problems in supply chain management and portfolio optimization often include bounds on the Conditional Value at Risk (CVaR), which can be cast as expectation constraints \citep{fabian2008handling,chen2010cvar}. Such constraints also arise when modeling distributionally robust versions of chance constraints \citep{wiesemann2014distributionally} and when limiting misclassification risk (i.e., misclassification rates) in multi-class Neyman Pearson classification \citep{weston1998multi,crammer2002learnability}. The aforementioned problems can be formulated as SOECs. Solutions violating risk constraints will likely fail stress tests that are performed before implementation. 
	\item Bounding property: Approximate linear programs (ALPs) are well-known models for approximating the value function of high-dimensional Markov decision processes \citep{schweitzer_generalized_1985,de_farias_linear_2003}, and in particular, are SOECs. A solution satisfying the ALP constraints provides an optimistic bound on the optimal policy value, which is useful to evaluate the suboptimality of heuristic policies. Infeasibility in an ALP setting thus voids this desirable bounding property.
\end{itemize}


Motivated by the importance of feasibility and the status quo of stochastic first order methods, we design an approach for solving SOECs that has low data complexity and provides high probability feasible solutions before convergence. As a first step, we cast SOEC as a root-finding problem involving a min-max level set function, which is challenging to solve because it is non-smooth and includes high-dimensional expectations in the SOEC objective and constraints. To solve this reformulation, we develop a stochastic feasible level-set method (SFLS) for root finding that requires evaluating a ``good" upper bound (we will make this notion of goodness precise in later sections) on the challenging level set function at each iteration. We show that employing the mirror descent method \citep{nemirov09} for computing such an upper bound requires approximating expectations in SOEC using SAAs at each iteration, which as already discussed above, leads to high data complexity. To overcome this issue, we introduce an SSG method to upper bound the level-set function by combining mirror-decent and online validation techniques, and in particular, extending the latter technique, originally proposed for minimization problems \citep{lan2012validation}, to handle saddle point formulations. This method only requires stochastic values and gradients of the objective and constraint functions, respectively, which can be constructed at low cost using a small number of samples of $\xi_i$ in \eqref{eq:gco}, that is, it has low data complexity. Calls to our SSG method return high-probability feasible solutions, which allows it to maintain an implementable solution at each root-finding iteration. 

We analyze the iteration complexity of SFLS to find a feasible solution path (i.e., sequence of feasible solutions) that becomes relative $\epsilon$-optimal with high probability. It is encouraging that the dependence of this complexity on $\epsilon$ is $1/\epsilon^2$, which is comparable to the method by \cite{yu2017online} (labeled YNW\footnote{We abbreviate this method by YNW using the first letters of the last names of the authors.}) that also finds an $\epsilon$-optimal solution but only guarantees $\epsilon$-feasibility at convergence. In other words, the intermediate solutions generated by YNW are not necessarily feasible. There is indeed a cost for ensuring feasibility in SFLS, which appears in the form of its iteration complexity depending on a condition measure. Such condition measures do not influence the complexity of YNW. \looseness=-1

For deterministic constrained convex optimization problems, the level-set method (DFLS) of \cite{lin2018level} also guarantees a feasible solution path with its iteration complexity depending on a condition measure. In principle, these DFLS based approaches can be applied to solve SOECs by viewing them as deterministic problems. This perspective is restrictive because it entails computing expectations in $f_i$ for $i=0,1,\dots,m$ exactly or replacing them by SAAs. In either case, the data complexity of DFLS will be high for reasons analogous to the ones already discussed above related to the use of SAAs. Therefore, a fully stochastic approach is required to achieve low data complexity when solving SOECs. \cite{pmlr-v80-lin18c} extend DFLS using variance-reduced sampling, which requires the functions have a finite-sum structure with each summand taking a specific form.\footnote{In particular, \cite{pmlr-v80-lin18c} require each summand has the form of $\phi(\bx^\top\mathbf{\xi})$.} Unfortunately, as a result, their method cannot be applied to SOECs with generic expectation while our method does not have such limitation and assumes little structure on the problems. We are not aware of prior efforts to develop a fully stochastic versions of level set methods for SOECs -- SFLS in this paper fills this gap. 

To assess the performance of SFLS, we provide implementation guidelines with supporting theory and numerically evaluate SFLS on three applications: (i) approximate linear programming for Markov decision processes, (ii) Neyman-Pearson multi-class classification with risk constraints, and (iii) learning a classifier with fairness constraints. Feasibility plays a key role in each of these applications for reasons mentioned earlier in the introduction. Approximate linear programs in the first application are known special cases of SOECs. For the latter two applications, we propose formulations that are SOECs. As algorithmic benchmarks, we consider YNW and DFLS. We find that SFLS delivers feasible solutions quicker than YNW and in several cases also leads to smaller optimality gaps. Moreover, when YNW computes infeasible solutions it is challenging to interpret its objective value since it can be superoptimal, an issue that does not arise with SFLS. Both SFLS and DFLS maintain feasible solution paths (with outer iterates) but SFLS produces feasible solutions with much smaller optimality gaps due to its lower data complexity. In other words, DFLS requires significantly more data passes to reduce the suboptimality of its solutions and will thus not be practical for solving SOECs based on large data sets. Our findings underscore two important algorithmic insights: (i) feasible SOEC solutions can be computed well before theoretical convergence criteria are satisfied but doing this hinges on methods being able to emphasize feasibility; and (ii) ensuring that these early feasible solutions have small optimality gaps requires approaches with low data complexity. Both these properties are true for SFLS, while only the first and second properties, respectively, hold for DFLS and YNW. \looseness=-1


This paper is organized as follows. In \S\ref{sec:SFLS}, we introduce SFLS, analyze its oracle complexity, and present a saddle-point reformulation of an SOEC. In \S\ref{StochasticOracle}, we discuss how the well-known stochastic mirror descent algorithm provides an idealized stochastic oracle for SFLS and highlight issues that complicate its use. In \S\ref{sec:compFeasStochOracle}, we propose and analyze a new stochastic oracle to overcome these issues. In \S\ref{StopComplexityOracle}, we analyze SFLS combined with this oracle and provide implementation guidelines. In \S\ref{sec:num}, we perform a computational study to understand the performance of SFLS across three applications relative to two benchmark methods. We conclude in \S\ref{sec:conclusion}.
\section{Stochastic Feasible Level-set Method}
\label{sec:SFLS}
Level-set methods tackle a constrained convex optimization problem by transforming it into a one-dimensional root-finding problem that is a function of a scalar level parameter $r$ \citep{Lemarechal1995,nesterov2004introductory}. We develop in this section a stochastic and feasible level set method that adds to this framework. We make the following standard assumption throughout the paper, which ensures that a strictly feasible and sub-optimal solution exists.
\begin{assumption}[Strict Feasibility] 
	There exists a strictly feasible solution $\tilde \bx \in \mathcal{X}$ such that
	\label{assum1}
	$\max_{i=1,\ldots,m}\{f_i(\tilde \bx)-r_i\}<0$ and $f_0(\tilde \bx) > f^*$. 
\end{assumption} 
The root-finding reformulation of \eqref{eq:gco} relies on the \emph{level-set function}
\begin{eqnarray}
\label{eq:gcols}
H(r):=\min_{\bx\in\mathcal{X}}\Pc(r,\bx)
\end{eqnarray}
where $r\in\mathbb{R}$ is a \emph{level parameter} and 
\[\Pc(r,\bx):=\max\left\{f_0(\bx)-r,f_1(\bx)-r_1,\dots,f_m(\bx)-r_m\right\}.\] Note that the expectation constraints of SOEC are now in the objective function of \eqref{eq:gcols}. For a given $(r, \bx) \in \mathbb{R} \times \mathcal{X}$, if $\Pc(r,\bx) \leq 0$ then $\bx$ is a feasible solution to \eqref{eq:gco}. Formulations \eqref{eq:gco} and \eqref{eq:gcols} are further linked by known properties of $H(r)$, which are summarized in the following lemma (based on lemmas 2.3.4 and 2.3.6 in \citealp{nesterov2004introductory} and Lemma 1 in \citealp{lin2018level}).
\begin{lemma}
	\label{lem:knownPropsOfL}
	It holds that 
	\begin{itemize}
		\item[(a)] $H(r)$ is non-increasing and convex in $r$;
		\item[(b)] $H(f^*)=0$;
		\item[(c)] $H(r) > 0$, if $r < f^*$ and $H(r) < 0$, if $r> f^*$.
	\end{itemize}
\end{lemma}
\noindent Part (a) of Lemma \ref{lem:knownPropsOfL} highlights that $H(r)$ is non-increasing and convex. Moreover, its part (b) implies that $r = f^*$ is the unique root of $H(r) = 0$. Therefore, one can use a root finding procedure to generate both a sequence of level parameters $r^{(1)},r^{(2)},\dots$ that converges to $f^*$ and an associated vector $\bx^{(k)} :=\argmin_{\bx\in\mathcal{X}}\Pc(r^{(k)},\bx)$ at each iteration $k$. Computationally, when a level parameter $r^{(k^*)}\approx f^*$ is found, the solution $\bx^{(k^*)} :=\argmin_{\bx\in\mathcal{X}}\Pc(r^{(k^*)},\bx)$ provides an ``approximate'' solution to \eqref{eq:gco}. From the perspective of feasibility, it is important whether we have $r^{(k^*)} < f^*$ or $r^{(k^*)} > f^*$. To elaborate, if $r^{(k^*)} < f^*$, then $H(r^{(k^*)}) > 0$ by Lemma \ref{lem:knownPropsOfL}(c) and the corresponding solution $\bx^{(k^*)}$ need not be feasible to~\eqref{eq:gco}. On the other hand, if $r^{(k^*)} >f^*$, we have $H(r^{(k^*)}) = \Pc(r^{(k^*)},\bx^{(k^*)}) < 0$ from Lemma \ref{lem:knownPropsOfL}(c) and the vector $\bx^{(k^*)}$ is indeed a feasible solution. A root finding scheme that ensures $r^{(k)} > f^*$ at each iteration $k$ will thus return a sequence of feasible solutions $\bx^{(1)}, \bx^{(2)}, \ldots, \bx^{(k^*)}$, that is a feasible solution path, where $k^*$ is such that $f^*<r^{(k^*)} < f^* + \epsilon$ for a given $\epsilon > 0$ and, in addition, we have $f_0(\bx^{(k^*)}) \leq r^{(k^*)}$ from $\Pc(r^{(k^*)},\bx^{(k^*)}) < 0$. These inequalities imply that $f_0(\bx^{(k^*)}) - f^* \leq \epsilon$. Thus, $\bx^{(k^*)}$ is an $\epsilon$-optimal and feasible solution to \eqref{eq:gco} and it follows that solving SOEC can be cast as a root-finding problem involving $H(r)$. \looseness = -1

Applying a root-finding algorithm to solve $H(r) = 0$ requires the exact computation of $H(r)$ at each iteration, which is difficult due to the nontrivial stochastic optimization in \eqref{eq:gcols}. Hence, we consider an inexact root-finding method, henceforth \emph{stochastic feasible level set method} (SFLS), extending what is done in~\cite{lin2018level} and~\cite{aravkin2019level} in a deterministic setting. Level set methods require an oracle to compute an approximation $U(r)$ of $H(r)$. This approximation is used to update $r$. A key element that we develop as part of SFLS is the notion of a stochastic oracle, which we introduce next.\looseness = -1

\begin{definition}[Stochastic Oracle]
	\label{def:feasOracle}
	Given $r>f^*$, $\epsilon>0$, and $\delta\in(0,1)$, a \textbf{stochastic oracle} $\mathcal{A}(r, \epsilon, \delta)$ returns a value $U(r)$ and a vector $\hat \bx\in\mathcal{X}$ that satisfy the inequalities $\Pc(r,\hat \bx)- H(r)\leq\epsilon$ and $|U(r)-H(r)| \leq \epsilon$ with a probability of at least $1-\delta$.\looseness = -1
\end{definition}

Lemma \ref{lem:feasOracle} clarifies the importance of the conditions underpinning the above definition to ensure a feasible solution to \eqref{eq:gco}.
\begin{lemma}\label{lem:feasOracle}
	Given $r>f^*$, $0< \epsilon\leq -\frac{\theta-1}{\theta+1}H(r)$, $\delta\in(0,1)$, and $\theta > 1$, the vector $\hat \bx\in\mathcal{X}$ returned by a stochastic oracle $\mathcal{A}(r, \epsilon, \delta)$ defines a feasible solution to \eqref{eq:gco} with probability of at least $1-\delta$. \looseness = -1
\end{lemma}
This lemma states that a stochastic oracle can recover a high probability feasible solution provided the optimality tolerance $\epsilon$ is less than $-\frac{\theta-1}{\theta+1}H(r)$.  


Algorithm~\ref{alg:lsdecreasing} formalizes the steps of SFLS to find an approximate root to $H(r) = 0$. Its inputs include a stochastic oracle $\mathcal{A}$; an initial level parameter value $r^{(0)}  > f^*$, which exists because we can set $r^{(0)} = f_0(\tilde \bx)$ by Assumption \ref{assum1}; optimality and error tolerances $\epsilon_{\text{opt}}$ and $\epsoracle$, respectively; a probability $\delta$; and a parameter $\theta$ that defines a step length as $1/2\theta$. SFLS begins from the level set defined by $r^{(0)}$. At each iteration $k$ it executes lines 3 though 9. In line 3, SFLS computes a probability $\delta^{(k)}$ that is used in the stochastic oracle call of line 4 to obtain an approximation $U(r^{(k)})$ and a high probability feasible solution $x^{(k)}$. The probability $\delta^{(k)}$ decreases with the iteration count $k$, that is, the probabilistic guarantee required of the stochastic oracle becomes more stringent to ensure the entire solution path is feasible with probability of at least $1-\delta$. Lines 5-7 model the termination condition, which involves checking whether the approximation $U(r^{(k)})$ is greater than or equal to $-\epsilon_{\text{opt}}$. If this condition holds, then the algorithm halts and returns the incumbent solution $x^{(k)}$. Otherwise, $r^{(k)}$ is updated to $r^{(k+1)}$ in line 8 using $U(r^{(k)})$ and $\theta$. Line 9 increments the iteration counter. While SFLS belongs to the family of level set approaches, it differs from known deterministic level set methods (see, e.g., \citealp{lin2018level} and \citealp{aravkin2019level}) in its update step, termination criterion, and stochastic oracle.\looseness=-1

\begin{algorithm}[h]
	\caption{Stochastic Feasible Level-Set Method (SFLS)} 
	\label{alg:lsdecreasing}
	\begin{algorithmic}[1]
		\STATE {\bfseries Inputs:} A stochastic oracle $\mathcal{A}$, a level parameter $r^{(0)}>f^*$, an optimality tolerance $\epsilon_{\text{opt}}>0$, an oracle error $\epsoracle>0$, a probability $\delta\in(0,1)$, and a step length parameter $\theta>1$. 
		\FOR{$k=0,1,\dots,$}
		\STATE $\delta^{(k)}=\dfrac{\delta}{2^k}$.
		\STATE $\left(U(r^{(k)}),\bx^{(k)}\right)=\mathcal{A}\left(r^{(k)},\epsoracle,\delta^{(k)}\right)$.
		\IF {$U(r^{(k)}) \geq -\epsilon_{\text{opt}}$} 
		\STATE Halt and return $\bx^{(k)}$.  
		\ENDIF
		\STATE $r^{(k+1)}\leftarrow r^{(k)}+U(r^{(k)})/(2\theta)$.
		\STATE $k\leftarrow k+1$.
		\ENDFOR
	\end{algorithmic}
\end{algorithm}

We define the notion of an \emph{input tuple} to ease the exposition of theoretical statements in the rest of the paper.
\begin{definition}[Input tuple] A tuple containing a subset of the elements $r, r^{(0)},\epsilon,\epsoracle,\delta,\theta$, and $\gamma_t$ is an {\bf input tuple} if its respective components satisfy $r>f^*$, $r^{(0)}>f^*$, $\epsilon > 0$, $\epsoracle > 0$, $\delta\in(0,1)$, $\theta > 1$, and $\gamma_t = 1/(M\sqrt{t+1})$, where $M > 0$ is a constant that is formally defined in \eqref{eq:defM}. \looseness=-1
\end{definition}
Theorem~\ref{generalcomplexity} provides the maximum number of calls to the stochastic oracle by Algorithm~\ref{alg:lsdecreasing} to obtain a feasible and relative $\epsilon$-optimal solution, which depends on a \emph{condition measure} $\beta$ of SOEC \eqref{eq:gco} defined as 
\begin{equation}
\label{eq:defbeta}
\beta := -\dfrac{H(r^{(0)})}{r^{(0)} - f^*}\in(0,1].
\end{equation}
It is easy to see that $\beta$ provides an assessment of the slope of $H(r)$ at $r = f^*$. Intuitively, for an SOEC instance with a large $\beta$ (i.e., well conditioned case), a root-finding method will be able to move towards the root of $H(r)$ faster compared to an instance with a small $\beta$ (i.e., ill-conditioned case). See Figure 2.1 of \citet{lin2018level} for a graphical illustration of this statement. 


\begin{theorem}
	\label{generalcomplexity}
	Given an input tuple $(r^{(0)},\epsilon,\delta,\theta)$, suppose $\epsilon_{\text{opt}}= -\frac{1}{\theta} H(r^{(0)})\epsilon$ and $\epsoracle= - \frac{\theta-1}{2\theta^2(\theta+1)}  H(r^{(0)}) \epsilon$.
	Algorithm~\ref{alg:lsdecreasing} generates a feasible solution at each iteration with a probability of at least $1-\delta$. Moreover, 
	it returns a relative $\epsilon$-optimal and feasible solution with this probability in at most $$\dfrac{2\theta^2}{\beta}\ln\left(\dfrac{\theta^2}{\beta\epsilon}\right)$$ calls to oracle $\mathcal{A}$. 
\end{theorem}
The bound on the number of oracle calls increases with $\theta$ because both the step-length $1/2\theta$ and the optimality tolerance $\epsilon_{\text{opt}}$ decrease with $\theta$. The maximum number of oracle calls is also a decreasing function of both the condition measure $\beta$ and tolerance $\epsilon$, that is, SFLS requires fewer iterations for problems that are better conditioned and when $\epsoracle$ and $\epsilon_{\text{opt}}$ are larger. Here, both $\epsoracle$ and $\epsilon_{\text{opt}}$  require knowledge of $H(r^{(0)})$, which is difficult to compute exactly. We want to point out that the dependence of $\epsoracle$ and $\epsilon_{\text{opt}}$ on $H(r^{(0)})$ are introduced here only to simplify the theorem and its proof, which helps readers to understand the main idea behind our technique. In \S\ref{StopComplexityOracle}, we will show that SFLS has a similar complexity even if $H(r^{(0)})$ in $\epsoracle$ and $\epsilon_{\text{opt}}$ is replaced by an upper bound $\bar U$ with $H(r^{(0)})\leq \bar U<0$ and $\bar U$ can be computed (by Algorithm \ref{alg:SMDwithOnlineValidationPDStop}) in a low cost independent of $\epsilon$.

\looseness=-1

SFLS relies on the availability of a valid stochastic oracle $\mathcal{A}$. Standard subgradient methods cannot be used as oracles to solve \eqref{eq:gcols} since computing a deterministic subgradient of $\Pc(r,\bx)$ requires exact evaluations of $f_i$ for $i=0,1,\dots,m$ (see \citealp{bertsekas1999nonlinear} or \citealp[p.737]{danskin2012theory}), which is challenging due to the high-dimensional expectations in the definition of these functions. 
Indeed, the expectation in each $f_i$ can be replaced by a direct SAA to obtain a sampled version $\hat\Pc(r,\bx)$ of $\Pc(r,\bx)$. This replacement is also problematic as subgradients of $\hat\Pc(r,\bx)$ provide biased subgradients of $\Pc(r,\bx)$ due to the maximization in the definition of the latter function. 

To avoid this issue, we reformulate \eqref{eq:gcols} into the equivalent min-max (i.e., saddle-point) form\looseness = -1
\begin{eqnarray}
\label{eq:fssaddleold}
H(r)=\min_{\bx\in\mathcal{X}}\max_{\by\in\mathcal{Y}}\left\{\sum_{i=0}^m y_i(f_i(\bx)-r_i)\right\},
\end{eqnarray}	 
where $r_0 := r$ and $\mathcal{Y}:=\left\{\by=(y_0,\dots,y_m)^\top\in\mathbb{R}^{m+1}|\sum_{i=0}^m y_i=1,y_i \geq 0 \right\}$. Given $\bx\in\mathcal{X}$, it is easy to check that $\by^* \in \arg\max_{\by\in\mathcal{Y}} \sum_{i=0}^m y_i(f_i(\bx)-r_i)$ can be chosen as a unit vector with 1 corresponding to an index $i^* \in \argmax_{i = 1,\ldots,m}\{f_i(\bx) - r_i\}$ and zeros for the remaining indices. Let $\Xi := \Xi_0\times\Xi_1\times \ldots\times\Xi_m$, $\bxi=(\xi_0,\xi_1,\dots,\xi_m)^\top\in\Xi$, $\Phi(\bx,\by,\bxi):=\sum_{i=0}^m y_i(F_i(\bx,\xi_i)-r_i),$ and $\phi(\bx,\by):=\E\left[\Phi(\bx,\by,\bxi)\right]$, where to ease notation we suppress the dependence of $\phi$ and $\Phi$ on the level parameter $r$ since it is always equal to a fixed value when these functions are invoked. Therefore, \eqref{eq:fssaddleold} can be reformulated as 
\begin{eqnarray}
\label{eq:fssaddleold1}
H(r)=\min_{\bx\in\mathcal{X}}\max_{\by\in\mathcal{Y}} \phi(\bx,\by).
\end{eqnarray}	 
Let $\hat{\phi}(\bx,\by)$ be an SAA of $\phi(\bx,\by)$. Subgradients of $\hat{\phi}(\bx,\by)$ provide an unbiased estimate of subgradients of $\phi(\bx,\by)$ because there is no nonlinear operator (e.g., maximization) acting on the expectation defining $\phi$. The oracles that we discuss for SFLS in \S\S\ref{StochasticOracle}-\ref{sec:compFeasStochOracle} will thus solve \eqref{eq:fssaddleold1}. 

\section{Idealized Stochastic Oracle}\label{StochasticOracle}

In \S\ref{subsec:SPRAndSMD}, we present stochastic mirror descent (SMD) in the form a stochastic oracle. In \S\ref{subsec:analIdealStochOracle}, we establish that SMD is indeed a stochastic oracle that can be used in SFLS (i.e., Algorithm~\ref{alg:lsdecreasing}) and then highlight computational issues that prevent its use. The discussion here serves a dual role. First, it provides practical motivation and sets the stage for developing a tractable stochastic oracle in \S\ref{sec:compFeasStochOracle}. Second, it provides basic concepts on primal-dual methods needed throughout the paper, also making the paper more accessible to readers potentially unfamiliar with such methods.

\subsection{Stochastic Mirror Descent}\label{subsec:SPRAndSMD}

Stochastic mirror descent {(SMD)} (\citealp{nemirov09}) is a well-known primal-dual method for solving saddle-point problems such as \eqref{eq:fssaddleold1}. SMD updates primal and dual variables $\bx$ and $\by$ of \eqref{eq:fssaddleold1}, respectively, by employing stochastic subgradients of $\phi(\bx,\by)$ and a projection operator. Let $F'_i(\bx,\xi_i)\in\partial F_i(\bx,\xi_i)$ for $i=0,1,\dots,m$, where $\partial$ is the subgradient operator. We denote the stochastic subgradient vector of $\phi(\bx,\by)$ by
\begin{eqnarray*}
	G(\bx,\by,\bxi):=\left[
	\begin{array}{c}
		G_x(\bx,\by,\bxi)\\
		-G_y(\bx,\by,\bxi)
	\end{array}
	\right]
	:=\left[
	\begin{array}{c}
		\sum_{i=0}^my_iF'_i(\bx,\xi_i)\\
		-(F_0(\bx,\xi_0)-r_0,F_1(\bx,\xi_1)-r_1,\dots,F_m(\bx,\xi_m)-r_m)^\top
	\end{array}
	\right].
\end{eqnarray*}
The projection employed by SMD relies on a distance function, known as Bregman divergence, that has as its argument $\bz :=(\bx,\by)$ and operates over $\mathcal{Z}:=\mathcal{X}\times \mathcal{Y}$. The space $\mathcal{Z}$ is equipped with a convex and continuously differentiable distance generating function $\omega_z(\bz)$ modulus 1 and a set of nonzero subgradients $\mathcal{Z}^o:=\{\bz\in\mathcal{Z}|\partial \omega_z(\bz)\neq \emptyset\}$. The Bregman divergence $V(\bz',\bz):\mathcal{Z}^o\times\mathcal{Z}\rightarrow\mathbb{R}_+$ expressed using $\omega_z$ is\looseness=-1
\[V(\bz',\bz):=\omega_z(\bz)-[\omega_z(\bz')+\nabla\omega_z(\bz')^\top(\bz-\bz')].\]
The projection operator (or prox-mapping), for any $\bzeta\in\mathbb{R}^{d+m+1}$, and $\bz'\in \mathcal{Z}^o$, is defined as
$P_{\bz'}(\bzeta):=\argmin_{\bz\in \mathcal{Z}}\left\{\bzeta^\top(\bz-\bz')+V(\bz',\bz)\right\}$.  

Algorithm~\ref{alg:SMD} summarizes the steps of SMD presented in the form of a stochastic oracle. The inputs to this algorithm are a level parameter $r\in\mathbb{R}$, an optimality tolerance $\epsoracle  > 0$, a probability $\delta\in (0,1)$, an iteration limit $W(\delta,\epsoracle)$ (we specify this limit later in Proposition \ref{validIdealOracle}), and a step-length rule $\gamma_t$ for all $t \in \mathbb{Z}_+$. Line 2 sets the initial solution $\bz^{(0)} = (\bx^{(0)}, \by^{(0)})$. Algorithm~\ref{alg:SMD} executes lines 4 and 5 for $W(\delta,\epsoracle)$ iterations. At iteration $t$, line 4 constructs a stochastic subgradient $G(\bx^{(t)},\by^{(t)},\bxi^{(t)})$ using a sample $\bxi^{(t)}$ of the random variables underlying the expectations in the objective and constraints of \eqref{eq:gco}. Line 5 computes a step-length weighted average $\bar\bz^{(t)}$ of past solutions. It also uses the stochastic subgradient computed in line 4 and a projection operator to find an updated solution $\bz^{(t+1)}$. After exiting the for loop, line 7 uses the averaged primal solution $\bar\bx^{(t)}$ to compute an upper bound $\max_{\by\in\mathcal{Y}}\phi(\bar\bx^{(t)},\by)$ on $H(r)$. The pair $(U(\bar \bx^{(t)}),\bar\bx^{(t)})$ is returned in line 8. \looseness=-1

It is worth noting that the update in line 5 relies on subgradients of an SAA $\hat\phi(\bx,\by)$ (with a single sample), which provides unbiased subgradients of $\phi(\bx,\by)$, unlike the biased subgradients that arise when working with SAAs of $\Pc(r,\bx)$ in the primal problem \eqref{eq:gcols}. In other words, a key benefit of the primal-dual reformulation \eqref{eq:fssaddleold} is that its objective $\phi(\bx,\by)$ allows the computation of unbiased subgradients after using SAAs to replace exact expectations.

\begin{algorithm}[h!]
	\caption{Stochastic Mirror Descent (SMD)} 
	\label{alg:SMD}
	\begin{algorithmic}[1]
		\STATE {\bfseries Inputs:} Level parameter $r\in\mathbb{R}$, optimality tolerance $\epsoracle > 0$, probability $\delta\in(0,1)$, an iteration limit $W$, and a step length rule $\gamma_t$ for all $t \in \mathbb{Z}_+$. 
		\STATE Set $\bz^{(0)}:=(\bx^{(0)},\by^{(0)})\in\argmin_{\bz\in\mathcal{Z}} \omega_z(\bz)$.
		\FOR{$t = 0,1,\ldots,W$}
		\STATE Sample $\bxi^{(t)}=(\xi_0^{(t)},\xi_1^{(t)},\dots,\xi_m^{(t)})^\top$ and compute $G(\bx^{(t)},\by^{(t)},\bxi^{(t)})$ .
		\STATE Execute 
		\begin{align*}
		&\bar\bz^{(t)} := 	(\bar\bx^{(t)},\bar\by^{(t)}):=\dfrac{\sum_{s=0}^{t}\gamma_s\bz^{(s)}}{\sum_{s=0}^{t}\gamma_s},\\ 
		&\bz^{(t+1)} := (\bx^{(t+1)},\by^{(t+1)}) :=P_{\bz^{(t)}}(\gamma_tG(\bx^{(t)},\by^{(t)},\bxi^{(t)})).
		\end{align*}
		\ENDFOR
		\STATE Compute $U(\bar \bx^{(t)}) = \max_{\by\in\mathcal{Y}}\phi(\bar\bx^{(t)},\by)$.
		\RETURN {$(U(\bar \bx^{(t)}),\bar\bx^{(t)})$}
	\end{algorithmic}
\end{algorithm}

\subsection{Validity of Stochastic Oracle and Computational Issues}\label{subsec:analIdealStochOracle}

We analyze below the validity of SMD as a stochastic oracle and also discuss its computational tractability. Our analysis, based on \cite{nemirov09}, requires specifying the distance generating function $\omega_z$ introduced in \S\ref{subsec:SPRAndSMD} and stating a standard assumption. 

To define $\omega_z$, we equip $\mathcal{X}$ and $\mathcal{Y}$
with their own distance-generating functions $\omega_x:\mathcal{X}\rightarrow \mathbb{R}$ modulus $\alpha_x$ with respect to norm $\|\cdot\|_x$ and $\omega_y:\mathcal{Y}\rightarrow \mathbb{R}$ modulus $\alpha_y$ with respect to norm $\|\cdot\|_y$. This means that $\omega_x$ is $\alpha_x$-strongly convex, continuous on $\mathcal{X}$, and continuously differentiable on the set of non-zero subgradients $\mathcal{X}^o:=\{\bx\in\mathcal{X}|\partial \omega_x(\bx)\neq \emptyset\}$. Similarly, $\omega_y$ is $\alpha_y$- strongly convex, continuous on $\mathcal{Y}$, and continuously differentiable on $\mathcal{Y}^o:=\{\by\in\mathcal{Y}|\partial \omega_y(\by)\neq \emptyset\}$. Typical choices for $\|\cdot\|_x$ and $\|\cdot\|_y$ are $\|\cdot\|_2$ and $\|\cdot\|_1$, respectively. In addition, it is common to set $w_x(\bx) = \frac{1}{2}\|\bx\|^2_2$ and $\omega_y(\by)=\sum_{i=0}^{m}y_i\ln y_i$. Defining the diameters of the sets $\mathcal{X}$ and $\mathcal{Y}$ as $D_x:=\sqrt{\max_{\bx\in\mathcal{X}}\omega_x(\bx)-\min_{\bx\in\mathcal{X}}\omega_x(\bx)}$ and $D_y  := \sqrt{\max_{\by\in\mathcal{Y}}\omega_y(\by)-\min_{\by\in\mathcal{Y}}\omega_y(\by)}$,  the distance-generating function associated with $\mathcal{Z}$ is 
$$
\omega_z(\bz):=\frac{\omega_x(\bx)}{2D_x^2}+\frac{\omega_y(\by)}{2D_y^2}.
$$ 

Next, the following standard assumption is needed to analyze SMD as well as other methods in the rest of the paper. Denote by $g(\bx,\by)$ expectation of the $(d+m+1)$-dimensional vector $G(\bx,\by,\bxi)$, that is, a deterministic subgradient. Moreover, let $\|\cdot\|_x$ and $\|\cdot\|_y$ represent the dual norms of $\|\cdot\|_{*,x}$ and $\|\cdot\|_{*,y}$, respectively.
\begin{assumption}\label{assum:sgd}
	For any $(\bx,\by,\bxi)\in\mathcal{X}\times \mathcal{Y}\times \Xi$, there exist $F'_i(\bx,\xi_i)\in\partial F_i(\bx,\xi_i)$ for $i=0,1,\dots,m$ such that 
	is well defined and satisfies
	\begin{eqnarray*}
		g(\bx,\by)
		\in
		\left[
		\begin{array}{c}
			\partial_x\phi(\bx,\by)\\
			\partial_y[-\phi(\bx,\by)]
		\end{array}
		\right],
	\end{eqnarray*}
	where $\partial_x$ and $\partial_y$ represent the sub-differentials with respect to $\bx$ and $\by$, respectively.
	Moreover, there exist positive constants $M_x$, $M_y$ and $Q$ such that 
	\begin{align}
	\label{eq:lighttailgradx}
	\E\left[\exp(\|G_x(\bx,\by,\bxi)\|_{*,x}^2/M_x^2)\right]&\leq\exp(1),\\
	\label{eq:lighttailgrady}
	\E\left[\exp(\|G_y(\bx,\by,\bxi)\|_{*,y}^2/M_y^2)\right]&\leq\exp(1),\\
	\label{eq:lighttailobj}
	\E\left[\exp(\left|\Phi(\bx,\by,\bxi)-\phi(\bx,\by)\right|^2/Q^2)\right]&\leq\exp(1),
	\end{align}
	for any $\bx\in\mathcal{X}$ and $\by\in\mathcal{Y}$, which indicate that $G_x$ and $G_y$  have a light-tailed distribution and their moments are bounded.
\end{assumption}

Proposition \ref{validIdealOracle} presents the iteration complexity of SMD, which follows from results in \cite{nemirov09}, and in addition, establishes that SMD is a valid stochastic oracle, that is, it satisfies Definition \ref{def:feasOracle}. The proof of this proposition relies on establishing that the primal-dual gap $U(\bar \bx^{(t)}) - L(\bar \by^{(t)})$ is guaranteed to be less than a given $\epsoracle > 0$ with a probability of at least $1-\delta$ for a given $\delta \in (0,1)$, where $L(\bar \by^{(t)}):= \min_{\bx\in\mathcal{X}}\phi(\bx,\bar\by^{(t)})$ and $U(\bar \bx^{(t)})$ is computed in Algorithm \ref{alg:SMD}. We also require the following constants:
\begin{align}
\label{eq:defM}
M &:=\sqrt{\frac{2D_x^2}{\alpha_x}M_x^2+\frac{2D_y^2}{\alpha_y}M_y^2};\\
\Omega(\delta) &:= \max\left\{\sqrt{12\ln\left(\dfrac{24}{\delta}\right)},\dfrac{4}{3}\ln\left(\dfrac{24}{\delta}\right)\right\}.\label{eq:defOmega}
\end{align}

\begin{proposition}\label{validIdealOracle}
	Given an input tuple $(r, \epsoracle,\delta, \gamma_t)$, the SMD solution $(\bar\bx^{(t)},\bar\by^{(t)})$ satisfies $U(\bar \bx^{(t)}) - L(\bar \by^{(t)}) \leq \epsoracle$ with probability at least $1- \delta$ in at most 
	\[W(\delta,\epsoracle) := \max\left\{6,\left(\frac{8\left(10 M \Omega(\delta)+4.5 M\right)}{\epsoracle}\ln\left(\dfrac{4\left(10 M \Omega(\delta)+4.5 M\right)}{\epsoracle}\right)\right)^2-2\right\}\]
	gradient iterations. 
	As a consequence, SMD is a valid stochastic oracle with $W\geq W(\delta,\epsoracle)$.
\end{proposition}
When solving \eqref{eq:fssaddleold1}, the dependence of the iteration complexity on $\epsoracle$ in Proposition \ref{validIdealOracle} has an additional $\ln (1/\epsoracle)$ term compared to the known SMD complexity dependence of $1/\epsoracle^2$ for solving an unconstrained version of this problem. Moreover, the analogous complexity dependence on $\delta$ inside logarithmic terms (see definition of $\Omega(\delta)$) in this proposition is comparable to the unconstrained case.\looseness=-1

We note that SMD is a valid stochastic oracle, exhibits a favorable iteration complexity, and is based on unbiased subgradients of $\phi(\bx,\by)$. Nevertheless, SMD is not directly implementable because the upper bound $U(\bar\bx^{(t)})$ is challenging to compute exactly as the definition of $\phi(\bx,\by)$ embeds expectations. Replacing these expectations by an SAA leads to a biased estimate of the upper bound $U(\bar\bx^{(t)})$. This bias can be reduced by using a large number of samples but doing this would lead to an approach with high data complexity, which we would like to avoid. In other words, although our saddle-point formulation facilitates the computation of unbiased subgradients needed by SMD to obtain a near optimal and high probability feasible solution, its upper bound $U(\bar\bx^{(t)})$, which serves as the constant $U(r)$ returned by the oracle (see Definition \ref{def:feasOracle}), cannot be computed. 

The aforementioned bound computation challenge is further exacerbated if one wishes to change the stopping criterion of Algorithm \ref{alg:SMD} (i.e., line 3) from a maximum iteration limit to a bound on the primal-dual gap $U(\bar \bx^{(t)}) - L(\bar \by^{(t)})$. In the latter case, implementing SMD would also entail the computation of the lower bound $L(\bar \by^{(t)})$, which suffers from analogous bias and data complexity issues when expectations in its definition are replaced by SAAs. In addition, the optimization problem over $\bx$ in the definition of $L(\bar \by^{(t)})$ is in general a high-dimensional non-smooth convex optimization problem and solving such a problem multiple times is computationally burdensome. Therefore, it is apriori unclear how one should go about designing a computationally tractable oracle to overcome these issues and what the iteration complexity of such an oracle would be.

\section{Tractable Stochastic Oracle}\label{sec:compFeasStochOracle}
In this section, we design a computationally viable stochastic oracle by combining SMD and an online validation technique \citep{lan2012validation}, and in particular, extending the latter technique originally proposed for minimization problems to handle min-max saddle point problems. This oracle overcomes the issues highlighted at the end of \S\ref{subsec:analIdealStochOracle} by defining bounds that are (i) tractable to compute with low data complexity and (ii) do not suffer from the bias issue when replacing expectations in their definitions by SAAs, as was the case with the bounds $U(\bar \bx^{(t)})$ and $L(\bar \by^{(t)})$. We present our algorithm in \S\ref{subsec:OnValidSMD} and prove that it is a stochastic oracle in \S\ref{subsec:ValidOnValidSMD}, where we also analyze its complexity. \looseness=-1
\subsection{Online Validation Based Stochastic Mirror Descent}\label{subsec:OnValidSMD}
Algorithm \ref{alg:SMDwithOnlineValidation} contains the steps of our proposed online validation based stochastic mirror descent (OVSMD) scheme, which differs from Algorithm \ref{alg:SMD} only in line 7, where the upper bound $U(\bar \bx^{(t)})$ on $H(r)$ is replaced by $\hat u_*^{(t)}$. The quantity $\hat u_*^{(t)}$ is an approximation of the following upper bound obtained using the online validation technique:
\begin{equation*}
u^{(t)}_* := \max_{\by\in\mathcal{Y}}\left\{\dfrac{1}{\sum_{s=0}^{t}\gamma_s}\sum_{s=0}^{t}\gamma_s\left[\phi(\bx^{(s)},\by^{(s)})+g_y(\bx^{(s)},\by^{(s)})^\top(\by-\by^{(s)})\right]\right\}.
\end{equation*}
This upper bound holds because 
\begin{align*}
\frac{1}{\sum_{s=0}^{t}\gamma_s}\sum_{s=0}^t\gamma_s\left[\phi(\bx^{(s)},\by^{(s)})+g_y(\bx^{(s)},\by^{(s)})^\top(\by-\by^{(s)})\right]
\geq\frac{\sum_{s=0}^{t}\gamma_s\phi(\bx^{(s)},\by)}{\sum_{s=0}^{t}\gamma_s}\geq\phi(\bar\bx^{(t)},\by),\end{align*} 
where the first inequality is true because $g_y$ is a subgradient with respect to $\by$ of the function $\phi(\bx,\by)$, which is concave in $\by$, and the second inequality follows directly from the convexity of $\phi(\bx,\by)$ in $\bx$. Therefore, we have 
\begin{align}
u^{(t)}_*\geq U(\bar \bx^{(t)}) = \max_{\by\in\mathcal{Y}}\phi(\bar\bx^{(t)},\by)\geq H(r),\label{eq:upperboundonHr}
\end{align}
that is, $u^{(t)}_*$ is an upper bound on $H(r)$, albeit potentially weaker than $U(\bar \bx^{(t)})$. Computing $u^{(t)}_*$ requires the exact evaluations of $\phi$, $g_x$ and $g_y$, which are not in general available because they involve expectations. In contrast, the term $\hat u^{(t)}_*$ computed in line 7 of Algorithm \ref{alg:SMDwithOnlineValidation}, which is stochastic approximation of $u^{(t)}_*$, can be easily computed in an online manner by solving a simple linear optimization problem. 

\begin{algorithm}[h!]
	\caption{Online Validation based Stochastic Mirror Descent: OVSMD} 
	\label{alg:SMDwithOnlineValidation}
	\begin{algorithmic}[1]
		\STATE {\bfseries Inputs:} Level parameter $r\in\mathbb{R}$, probability $\delta\in(0,1)$, optimality tolerance $\epsoracle > 0$, an iteration limit $T$, and a step length rule $\gamma_t$ for all $t \in \mathbb{Z}_+$. 
		\STATE Set $\bz^{(0)}:=(\bx^{(0)},\by^{(0)})\in\argmin_{\bz\in\mathcal{Z}} \omega_z(\bz)$.
		\FOR{$t = 0,1,\ldots,T(\delta,\epsoracle)$}
		\STATE Sample $\bxi^{(t)}=(\xi_0^{(t)},\xi_1^{(t)},\dots,\xi_m^{(t)})^\top$ and compute $G(\bx^{(t)},\by^{(t)},\bxi^{(t)})$ .
		\STATE Execute 
		\begin{align*}
		&\bar\bz^{(t)} := 	(\bar\bx^{(t)},\bar\by^{(t)}):=\dfrac{\sum_{s=0}^{t}\gamma_s\bz^{(s)}}{\sum_{s=0}^{t}\gamma_s},\\ 
		&\bz^{(t+1)} := (\bx^{(t+1)},\by^{(t+1)}) :=P_{\bz^{(t)}}(\gamma_tG(\bx^{(t)},\by^{(t)},\bxi^{(t)})).
		\end{align*}
		\ENDFOR
		\STATE Compute 
		\begin{equation}\label{eqn:compUBExpression}
		\hat u_*^{(t)}:=\max_{\by\in\mathcal{Y}}\left\{\frac{1}{\sum_{s=0}^{t}\gamma_s}\sum_{s=0}^{t}\gamma_s\left[\Phi(\bx^{(s)},\by^{(s)},\bxi^{(s)})+G_y(\bx^{(s)},\by^{(s)},\bxi^{(s)})^\top(\by-\by^{(s)})\right]\right\}.
		\end{equation}
		\RETURN {$(\hat u_*^{(t)},\bar\bx^{(t)})$}
	\end{algorithmic}
\end{algorithm}

As discussed in \S\ref{subsec:analIdealStochOracle}, replacing the iteration limit based stopping criterion by one that approximates an optimality gap requires a lower bound on $H(r)$. Following a similar argument to the upper bounding case above, we define the lower bound
\begin{equation*}l^{(t)}_* := \min_{\bx\in\mathcal{X}}\left\{\dfrac{1}{\sum_{s=0}^{t}\gamma_s}\sum_{s=0}^{t}\gamma_s\left[\phi(\bx^{(s)},\by^{(s)})+g_x(\bx^{(s)},\by^{(s)})^\top(\bx-\bx^{(s)})\right]\right\}.
\end{equation*}
Since $\phi(\bx,\by)$ is convex in $\bx$, it follows that $l^{(t)}_* \leq L(\bar \by^{(t)}) = \min_{\bx\in\mathcal{X}}\phi(\bx,\bar\by^{(t)})\leq H(r)$. Although $l^{(t)}_*$ is in general a weaker lower bound than $L(\bar \by^{(t)})$, the former bound is computed by solving a linear optimization problem as opposed to the potentially challenging non-smooth convex optimization problem defining the latter bound. Finally, we employ an online validation based approximation of $l^{(t)}_*$ to avoid computing expectations and obtain 
\begin{equation}\label{eqn:compLBExpression}\hat l_*^{(t)}:=\min_{\bx\in\mathcal{X}}\left\{\frac{1}{\sum_{s=0}^{t}\gamma_s}\sum_{s=0}^{t}\gamma_s\left[\Phi(\bx^{(s)},\by^{(s)},\bxi^{(s)})+G_y(\bx^{(s)},\by^{(s)},\bxi^{(s)})^\top(\bx-\bx^{(s)})\right]\right\}.\end{equation}

Despite the computational tractability of $\hat u^{(t)}_*$ and $\hat l_*^{(t)}$, these are stochastic quantities and subject to noise. Hence they do not always provide valid bounds on $H(r)$. In \S\ref{subsec:ValidOnValidSMD}, we show that $\hat l^{(t)}_*$ and $\hat u^{(t)}_*$ are nevertheless sufficiently close to $H(r)$ with high probability after a finite number of iterations (see Theorem \ref{thm:totalcomplexity}). \looseness = -1

\subsection{Validity of Stochastic Oracle and Iteration Complexity}\label{subsec:ValidOnValidSMD}
We establish here the validity of OVSMD (i.e., Algorithm \ref{alg:SMDwithOnlineValidation}) as a stochastic oracle and derive its iteration complexity. Proposition \ref{thm:luconvergecor} contains the two main ingredients underlying the analysis of OVSMD. Part (i) of this proposition shows that for a given $\epsoracle > 0$ the inequality $u^{(t)}_* - l^{(t)}_* \leq \epsoracle$ holds with high probability when $t$ is sufficiently large. In other words, the deterministic quantities $u^{(t)}_*$ and $l^{(t)}_* $ computed using the OVSMD solutions provide ``good'' deterministic estimates of the level set function $H(r)$. This is not directly useful since OVSMD can only compute stochastic approximations of these quantities, as already discussed in \S\ref{subsec:OnValidSMD}. Part (ii) of Proposition \ref{thm:luconvergecor} establishes that $\hat u^{(t)}_*$ and $\hat l^{(t)}_*$ are respectively close stochastic approximations of $u^{(t)}_*$ and $l^{(t)}_*$ at convergence with high probability. It then follows that the quantities $\hat u^{(t)}_*$ and $\hat l^{(t)}_*$ are ``good'' stochastic estimates of the level set function, and in particular, allows OVSMD to be used as a stochastic oracle.


\begin{proposition} 
	\label{thm:luconvergecor}	
	Given an input tuple $(r,\epsoracle, \delta, \gamma_t)$, OVSMD computes $(\bx^{(t)},\by^{(t)})$, $t = 1,2,3,\ldots,$ such that:
	\begin{itemize}
		\item[(i)] The inequality $\text{Prob}\{u^{(t)}_*-l^{(t)}_*>\epsoracle\}\leq\delta$ holds in at most
		$$
		\max\left\{6,\left(\frac{8\left(10 M \Omega(\delta)+4.5 M\right)}{\epsoracle}\ln\left(\dfrac{4\left(10 M \Omega(\delta)+4.5 M\right)}{\epsoracle}\right)\right)^2-2\right\}
		$$ 
		gradient iterations.
		\item[(ii)] The inequalities $\text{Prob}\{|\hat l^{(t)}_*-l^{(t)}_*|>\epsoracle\}\leq\delta$ and $\text{Prob}\{|\hat u^{(t)}_*-u^{(t)}_*|>\epsoracle\}\leq\delta$ hold in at most\looseness = -1
		$$
		\max\left\{6,\left(\frac{8\left(Q \Omega(\delta) +8 M \Omega(\delta)+2.5M\right)}{\epsoracle}\ln\left(\frac{4\left( Q \Omega(\delta)+8 \Omega(\delta) M+2.5M\right)}{\epsoracle}\right)\right)^2-2\right\}
		$$  
		gradient iterations. 
	\end{itemize}
\end{proposition}



Leveraging Proposition \ref{thm:luconvergecor}, Theorem \ref{thm:totalcomplexity} shows that OVSMD is a valid stochastic oracle and also presents its iteration complexity.

\begin{theorem} 	\label{thm:totalcomplexity}	
	Given an input tuple $(r,\epsoracle, \delta, \gamma_t)$, the OVSMD guarantees $\Pc(r,\bar\bx^{(t)}) - H(r) \leq \epsoracle$ and $\vert\hat u_*^{(t)} - H(r)\vert \leq \epsoracle$ with probability at least $1-\delta$ in at most 
	\small
	\begin{equation}
	\label{eq:TvalueH}
	T(\delta,\epsoracle) :=\max\left\{6,\left(\frac{16\left( Q\Omega(\delta)+10 M \Omega(\delta)+4.5M\right)}{\epsoracle}\ln\left(\frac{8\left( Q \Omega(\delta)+10 M \Omega(\delta)+4.5M\right)}{\epsoracle}\right)\right)^2-2\right\}
	\end{equation}
	\normalsize
	gradient iterations. As a consequence, OVSMD is a valid stochastic oracle with $T\geq T(\delta,\epsoracle)$. 
\end{theorem}
Despite OVSMD being a tractable oracle, the dependence of its iteration complexity on both $\epsoracle$ and $\delta$ is identical to the analogous dependence seen with the idealized SMD oracle analyzed in Proposition \ref{validIdealOracle}. Moreover, in terms of $\epsoracle$, OVSMD is only a $\ln(1/\epsoracle)$ worse than the known complexity of SMD in the unconstrained case, where feasibility is not a concern.


\section{SFLS with OVSMD as its Stochastic Oracle}\label{StopComplexityOracle}
In this section, we provide theoretical support for the use of OVSMD as SFLS's stochastic oracle in \S\ref{subsec:SFLSOVSMDAnal} and then discuss implementation guidelines in \S\ref{subsec:SFLSImplementation}.
\subsection{Theoretical Analysis}\label{subsec:SFLSOVSMDAnal}
Theorems \ref{generalcomplexity} and \ref{thm:totalcomplexity} can be used to derive the (gradient) iteration complexity of SFLS when using OVSMD as the stochastic oracle. We state this complexity in Corollary \ref{thm:FullComplexityWithOVSMD}. 
\begin{corollary}\label{thm:FullComplexityWithOVSMD}
	Given an input tuple $(r^{(0)}, \epsilon,\delta, \gamma_t,\theta)$, let $\epsilon_{\text{opt}}= -\frac{1}{\theta} H(r^{(0)})\epsilon$ and $\epsoracle= - \frac{\theta-1}{2\theta^2(\theta+1)}  H(r^{(0)}) \epsilon$. Moreover, suppose OVSMD with $T= T(\delta,\epsoracle)$ is chosen as the stochastic oracle $\mathcal{A}$. Then SFLS returns a relative $\epsilon$-optimal and feasible solution with probability of at least $1-\delta$ using at most 
	\[\dfrac{2\theta^2}{\beta}\ln\left(\dfrac{\theta^2}{\beta\epsilon}\right)\]
	OVSMD calls and \looseness = -1
	$$\mathcal{O}\left(\dfrac{\theta^2}{\beta\epsilon^2}\cdot \ln\left(\dfrac{\theta^2}{\beta\epsilon}\right)\cdot\ln^2\left(\dfrac{1}{\delta}\right)\cdot\ln^2\left(\dfrac{1}{\epsilon}\right)\right)$$
	gradient iterations. 
\end{corollary}
This complexity result is somewhat idealistic because the inputs to SFLS, namely $\epsilon_{\text{opt}}$ and $\epsoracle$, require knowledge of $H(r^{(0)})$, which is difficult to compute exactly. A possible resolution is to compute an upper bound on $H(r^{(0)})$, denoted by $\bar U$, such that $H(r^{(0)})\leq \bar U<0$. If $|\bar{U}|$ is much smaller than $|H(r^{(0)})|$, then the optimality tolerance $\epsoracle$ will be substantially more stringent and thus lead to a larger complexity than the iteration bound in Corollary \ref{thm:FullComplexityWithOVSMD}. Therefore, to obtain a complete theoretical assessment of the computational complexity of SLFS with OVSMD, it is important to incorporate the cost of finding a $\bar{U}$ that is comparable to $H(r^{(0)})$ (i.e., $|\bar{U}| = \Omega(|H(r^{(0)})|)$). \looseness=-1


\begin{algorithm}[h!]
	\caption{Estimating an upper bound on $H(r^{(0)})$ using OVSMD} 
	\label{alg:SMDwithOnlineValidationPDStop}
	\begin{algorithmic}[1]
		\STATE {\bfseries Inputs:} Level parameter $r^{(0)}>f^*$, initial approximation tolerance $\bar\alpha > 0$, probability $\delta\in(0,1)$, constant $\theta > 1$, and a step length rule $\gamma_h$ for all $h \in \mathbb{Z}_+$. 
		\STATE Set $h=0$ and $\alpha^{(0)} = \bar\alpha $.
		\REPEAT 
		\STATE Set $\delta^{(h)}=\dfrac{\delta}{2^{h+1}}$ and $\alpha^{(h)}=\dfrac{\alpha^{(0)}}{2^{h}}$.
		\STATE Compute $\hat u^{(h)}_*\leftarrow \text{OVSMD}(r^{(0)},\delta^{(h)},\alpha^{(h)},\gamma_h)$.
		\UNTIL{$\hat u_*^{(h)}+\alpha^{(h)}<0$ and $\dfrac{\hat u_*^{(h)}-\alpha^{(h)}}{\hat u_*^{(h)}+\alpha^{(h)}}\leq \theta$.}
		\RETURN {$\bar U=\hat u_*^{(h)}+\alpha^{(h)}$.}
	\end{algorithmic}
\end{algorithm}

Fortunately, OVSMD can itself be used to compute the desired $\bar{U}$. We discuss the intuition behind its use for this purpose and then formally state the result. Recall that $H(r^{(0)})<0$ since $r^{(0)} > f^*$. We consider obtain an upper bound $\bar U$ by solving \eqref{eq:gcols} with $r=r^{(0)}$ and a small enough optimality gap. By Theorem~\ref{thm:totalcomplexity}, OVSMD with $r=r^{(0)}$ can guarantee $H(r^{(0)})\leq \hat u_*^{(t)}+\epsoracle$ with high probability. This suggests setting $\bar U=\hat u_*^{(t)}+\epsoracle$. However, it is a priori unclear how small $\epsoracle$ should be in order to ensure $\bar U<0$ and $|\bar U|=\Omega(|H(r^{(0)})|)$. Therefore, we run OVSMD multiple times, starting from a tolerance $\alpha^{(0)} = \bar\alpha $, geometrically reducing this tolerance after each run, and stopping this procedure once $\bar U=\hat u_*^{(h)}+\alpha^{(h)}<0$ and ${(\hat u_*^{(h)}-\alpha^{(h)})}/{(\hat u_*^{(h)}+\alpha^{(h)})}\leq \theta$ hold. We can then use Theorem~\ref{thm:totalcomplexity} and the condition $\hat u_*^{(h)}+\alpha^{(h)}<0$ to show that  ${|H(r^{(0)})|}/{|\bar U|}\leq {(\hat u_*^{(h)}-\alpha^{(h)})}/{(\hat u_*^{(h)}+\alpha^{(h)})}\leq \theta$, which implies $|\bar U|=\Omega(|H(r^{(0)})|)$. We formalize the aforementioned approach in Algorithm~\ref{alg:SMDwithOnlineValidationPDStop}. \looseness=-1


Theorem \ref{cor:FullComplexityWithOVSMD_U} establishes the complexity of employing Algorithm \ref{alg:SMDwithOnlineValidationPDStop} to compute $\bar U$ and subsequently running SFLS leveraging this computation. 
\begin{theorem}\label{cor:FullComplexityWithOVSMD_U}
	Given an input tuple $(r^{(0)},\epsilon, \delta, \gamma_t,\theta)$, suppose we compute $\bar{U}$ using Algorithm \ref{alg:SMDwithOnlineValidationPDStop} and then execute SFLS to find a relative $\epsilon$-optimal and feasible solution with a probability of at least $1-\delta$ using $\epsilon_{\text{opt}}= -\frac{1}{\theta} \bar U\epsilon$, $\epsoracle= - \frac{\theta-1}{2\theta^2(\theta+1)}  \bar U \epsilon$, and OVSMD with $T= T(\delta,\epsoracle)$ as the stochastic oracle $\mathcal{A}$.  This procedure requires in total at most 
	$$\mathcal{O}\left(\dfrac{\theta^2}{\beta}\ln\left(\dfrac{\theta^2}{(1-\theta)\beta\epsilon}\right)\right)$$
	OVSMD calls and 
	$${\mathcal{O}}\left(\dfrac{1}{\beta^2}\ln^4\left(\frac{1}{\beta}\right)\ln^2\left(\dfrac{1}{\delta}\right)\right) + \mathcal{O}\left(\dfrac{\theta^2}{\beta\epsilon^2}\cdot \ln\left(\dfrac{\theta^2}{\beta\epsilon}\right)\cdot\ln^2\left(\dfrac{1}{\delta}\right)\cdot\ln^2\left(\dfrac{1}{\epsilon}\right)\right) $$
	gradient iterations. 
\end{theorem}
Theorem \ref{cor:FullComplexityWithOVSMD_U} provides a realistic theoretical assessment of the computational burden of solving SOECs using SFLS. Interestingly, it shows that running Algorithm \ref{alg:SMDwithOnlineValidationPDStop} to compute $\bar{U}$ before executing SFLS and replacing the unknown term $H(r^{(0)})$ in the definitions of $\epsoracle$ and $\epsilon_{\text{opt}}$ with the computed $\bar{U}$ value  does not 
change the overall big-$\mathcal{O}$ oracle and gradient iteration complexities in Corollary \ref{thm:FullComplexityWithOVSMD}, except for logarithmic terms. \looseness=-1

The complexity of SFLS (combined with OVSMD) in Theorem \ref{cor:FullComplexityWithOVSMD_U} is comparable in terms of its dependence on $\epsilon$ and $\delta$ to the complexity of the algorithm in \citet{yu2017online}, which does not ensure feasibility. This suggests that our procedure is efficient at ensuring feasibility. The cost of ensuring feasibility, however, appears in the dependence of the SFLS iteration complexity on the condition measure $\beta$. Such dependence is absent in approaches that do not ensure feasibility. 

Another relevant comparison is with the deterministic feasible level set approach (DFLS) of \cite{lin2018level} and its variant in \cite{pmlr-v80-lin18c}, which are both applicable to solve deterministic constrained convex optimization problems. The complexity of DFLS based methods depend on the number of data points that define expectations and thus lead to large data complexity, and in particular, have infinite complexity when expectations are defined over continuous random variables. In contrast, the complexity of SFLS in Theorem \ref{cor:FullComplexityWithOVSMD_U} does not depend on the number of data points. In addition, compared to DFLS, the iteration complexity of SFLS has only additional logarithmic factors involving $\epsilon$ and $\delta$, which is encouraging, as the stochastic level set algorithm (i.e., Algorithm \ref{alg:lsdecreasing}) and OVSMD oracle need to contend with several challenges that arise due to the presence of expectations in SOECs. 

In summary, our theoretical analysis of SFLS and comparison with known complexities of state-of-the-art approaches suggests that SFLS is effective in terms of iteration complexity at computing a high probability feasible solution path for SEOCs, a much broader and challenging class of problems than deterministic constrained convex programs. Moreover, a fully stochastic approach such as SFLS is theoretically necessary to achieve low data complexity in this context. 

\subsection{Implementation Guidelines}\label{subsec:SFLSImplementation}

As is common with first-order methods, the implementation of SFLS requires parameter tuning. A direct implementation of SFLS in a manner consistent with Theorem~\ref{cor:FullComplexityWithOVSMD_U} requires selecting $r^{(0)}$, $\epsilon$, $\delta$, $\theta$ and $\gamma_t$; estimating constants $M$ and $Q$ (needed to define $T=T(\delta,\epsoracle)$ in OVSMD); and then computing $\bar{U}$. While these parameters can be estimated or approximated, we suggest a simpler implementation strategy that largely side-steps such tuning. Firstly, we avoid stopping SFLS by pre-specifying optimality tolerance $\epsilon_{\text{opt}}$ and instead stop it based on an outer iteration limit. This is possible because the SFLS outer iterations only affect the suboptimality of the incumbent feasible solution, that is, being a feasible level set method, SFLS can return feasible and implementable solutions when terminated after any number of outer iterations. Secondly, instead of choosing the number of inner iteration in OVSMD as $T=T(\delta,\epsoracle)$ based on pre-specified $\delta$ and $\epsoracle$, we directly specify $T$. According to \eqref{eq:TvalueH}, $T(\delta,\epsoracle)$ is strictly monotonically decreasing in $\epsoracle$ and thus in $\epsilon$ so that a relative $\epsilon$-optimal and feasible solution with $\epsilon=\tilde{\mathcal{O}}(\frac{1}{\sqrt{T}})$ can be guaranteed. Corollary \ref{cor:implStop} establishes that the convergence of SFLS in this implementation. 
\begin{corollary}\label{cor:implStop}
Suppose we have an input tuple $(r^{(0)},\gamma_t, \theta)$ and the iteration limit in OVSMD is $T$. Given $\delta\in(0,1)$, SFLS finds a relative $\epsilon$-optimal and feasible solution with $\epsilon\leq \mathcal{O}\left(\frac{\theta^4\ln(1/\delta)\ln(T)\ln\left(T/\beta\right)}{\beta\sqrt{T}}\right)$ and with a probability of at least $1-\delta$ 	using at most $\mathcal{O}\left(\frac{\theta^2}{\beta}\ln\left(\frac{T}{\beta}\right)\right)$ OVSMD calls and $\mathcal{O}\left(\frac{\theta^2 T}{\beta }\ln\left(\frac{T}{\beta}\right)\right)$
	gradient iterations. 
\end{corollary}
Overall, following the aforementioned strategy only requires the choice of $T,\theta$, $r^{(0)}$, and $\gamma_t$ -- a significant reduction in implementation burden.

For choosing $\theta$ and $T$, we consider a discrete set of values and tune the algorithm, that  is, we test the performance of SFLS for a few iterations or data passes for each value, and select the one that leads to the largest decrease in suboptimality. Selecting $r^{(0)}$ is easy when an initial feasible solution $\tilde{\bx}$ is available because we have $\E\left[F_0(\tilde{\bx},\xi_0)\right] > f^*$. In this case, we estimate $\E\left[F_0(\tilde{\bx},\xi_0)\right]$ using an SAA and then set $r^{(0)}$ to a larger value to account for approximation error and ensure we have $r^{(0)} > f^*$. If a feasible solution is not readily available, we can find one by applying a minor modification of Algorithm \ref{alg:SMDwithOnlineValidationPDStop} to solve 
\[\min_{\bx\in\mathcal{X}}\max_{\by\in\mathcal{Y}}\left\{\sum_{i=1}^m y_i(f_i(\bx)-r_i)\right\},\]
which does not include the term in \eqref{eq:fssaddleold} corresponding to $i = 0$, that is, $f_0 - r$. Finally, the step length can be specified as $\gamma_t = 1/(c\sqrt{t+1})$ for a given constraint $c > 0$, which is tuned. While $c$ is chosen as $M$ in our theoretical analysis to simplify proofs, analogous results hold for a generic constant $c > 0$. We omit these general results for the sake of brevity as they do not change the dependence of our iterations bounds on $\epsilon$, $\beta$, and $\delta$. \looseness=-1


\section{Numerical Experiments}
\label{sec:num}
In this section, we evaluate the numerical performance of SFLS on three diverse SOEC applications: (i) approximate linear programs for solving Markov decision processes, (ii) multi-class Neyman-Pearson classification, and (iii) learning with fairness constraints. SOECs in the first application contain expectations of continuous random variables while those in the second and third applications involve discrete random variables. Our first algorithmic benchmark is the stochastic subgradient method YNW of \citet{yu2017online} as it is the only first order approach (we are aware of) that can handle SOECs with multiple constraints. In addition, we also compare against the deterministic feasible level-set method (DFLS) of \cite{lin2018level} because it ensures a feasible solution path. Specifically, comparing SFLS and DFLS allows us to evaluate the benefits of the reduced data complexity in our stochastic approach. In \S\ref{subsec:AlgoImpl}, we describe our computational setup and then the performance of algorithms on applications in \S\S\ref{subsec:PersInvVontrol}-\ref{appl:LearnDataConst}.
\subsection{Computational Setup}\label{subsec:AlgoImpl}
We implemented SFLS, DLFS, and YNW in Matlab running on a 64-bit Microsoft Windows 10 machine with a 2.70 Ghz Intel Core i7-6820HQ CPU and 8GB of memory. We set $\omega_x(\bx)=\frac{1}{2}\|\bx\|_2^2$ and  $\omega_y(\by)=\sum_{i=0}^{m}y_i\ln y_i$ in all three algorithms. 
We followed the guidelines in \S\ref{subsec:SFLSImplementation} when implementing SFLS and thus had to choose only $r^{(0)}$, $\theta$, and $\gamma_t$. We based $r^{(0)}$ on the solution $\tilde\bx$. We tuned $\theta$ over the discrete set $\{1.1,2,5\}$ and $T$ over the discrete set $\{50,100,200,300\}$. We selected $\gamma_t = 1/(c\sqrt{t+1})$ and tuned $c$ over the set of possible values $\{0.05, 0.1, 1, 2, 5\}$. We employed a mini-batch technique to construct the stochastic gradients in SFLS and YNW. 


Similar to SFLS, DFLS solves the subproblem $\min_{\bx\in\mathcal{X}}\Pc(r^{(k)},\bx)$ approximately in the $k$th outer iteration and uses the returned solution $\bx^{(k)}$ to update $r^{(k)}$ as $r^{(k+1)}\leftarrow r^{(k)}+\Pc(r^{(k)},\bx^{(k)})/2$. 
Following \cite{lin2018level}, we use the standard subgradient descent method to solve this subproblem and the parameters $r^{(0)}$ and $\gamma_t$ and the inner iteration limit $T$ in DLFS are tuned in the same way as in SFLS as described above. To apply DLFS, we constructed a deterministic version of each SOEC using SAAs of expectations. We found, consistent with \cite{lin2018level}, that using SAAs in lieu of expectations over continuous random variables in the perishable control problem (first application) did not sufficiently represent the original problem even when using a large number of samples. We thus omitted DFLS as a benchmark for this application. This was not an issue for the remaining two applications because expectations are defined over discrete random variables. To avoid the quality of SAAs confounding our performance evaluation, we chose instances for these two applications such that expectations can be evaluated exactly, albeit requiring more time. 

We followed the guidance in \cite{yu2017online} to setup YNW. Specifically, we chose the control parameters $V$ and $\alpha$ as $V =\sqrt{T}$ and 
$\alpha=T$, respectively, as a function of the total number of iterations $T$, where $V$ is the weight of the gradient of the objective function and $\alpha$ is the weight of the proximal term in the updating equation of $\bx$ in YNW. Similar to SFLS, we used a mini-batch technique to construct the stochastic gradients and evaluate the objective values.
\subsection{Approximate Linear Programming for Markov Decision Processes}\label{subsec:PersInvVontrol}
Approximate linear programs (ALPs) address the well-known curse of dimensionality associated with directly solving large-scale Markov decision processes (MDPs; \citealt{puterman1994markov}) by computing a value function approximation. We illustrate how our SFLS method can be applied to tackle ALPs, and thus large-scale MDPs, by considering a challenging perishable inventory control problem with partial backlogging and lead time. We begin by presenting the MDP for this problem and refer the reader to \citet{lin2019revisiting} for its derivation and detailed application context.

Consider the management of orders for a single product with a finite life time of $I$ periods and an order lead time of $J$ periods, that is, the product takes $J$ periods to be delivered from when it is ordered and $I$ periods to perish from receipt. The state space of the MDP is represented by the vector \looseness=-1 \[\bs=(z_{0},z_1,\dots,z_{I-1},q_{1},q_{2},\dots,q_{J-1})\in\mathbb{R}^{I + J - 1},\]
where $q_j$, $1\leq j \leq J-1$, denotes the order quantities that will be received $j$ periods from now, and $z_i$, $0 \leq i \leq I-1$, the on-hand inventory with $i$ periods of lifetime remaining. The order quantity $a$ is at most $\bar{a}$ and belongs to the interval $[0,\bar{a}]$, which implies $z_i \in [0,\bar{a}]$ for $i = 1,\ldots,I-1$ and $q_j \in [0,\bar{a}]$ for $j = 1,\ldots,J-1$. The element $z_0$ of the state is bounded below by $l_s < 0$ to allow limited or partial backlogging, that is, any units backlogged beyond $|l_s|$ are lost sales. To ease exposition, we write $\bs \in \mathcal{S}$ and $a \in \mathcal{A}$ to capture the state and action domains, respectively, and use $\bs^0$ to represent the initial state. Assuming orders are served on a first-come-first-serve basis, the MDP state transitions as  
\[f(\bs,a) = (\max\{z_{1}-(G-z_{0})_+,l_s - \sum_{i=2}^{I-1}z_i\},z_{2},\ldots,z_{I-1},q_{1},q_{2},\ldots,q_{J-1},a),\]
\noindent where $G$ represents stochastic demand with distribution $P_G$.
Moreover, the cost associated with ordering $a$ at state $\bs$ is 

{\small\begin{eqnarray*}
		&&c(\bs,a)=\gamma^J c_p a  \\
		&&+\mathbb{E}\Bigg[c_h\left(\sum_{i=1}^{I-1}z_{i}-(G - z_{0})_+\right)_+ + c_b\left(G-\sum_{i=0}^{I-1}z_{i}\right)_+ + c_d\left(z_{0}-G\right)_+ + c_l\left(l_s + G - \sum_{i=0}^{I-1}z_{i}\right)_+\Bigg],
\end{eqnarray*}}
where the per unit lost sale, disposal, purchasing, holding, and backlogging costs are $c_l$, $c_d$, $c_p$, $c_h$, and $c_b$, respectively; $\mathbb{E}$ is taken over $G$; and $\gamma \in (0,1)$ is a discount factor. The infinite horizon (discounted cost) MDP formulated using the aforementioned components can be solved using the fixed point equations
\[V(\bs) = \max_{a \in \mathcal{A}} \;c(\bs,a) + \gamma \mathbb{E}[V(f(\bs,a))], \quad\forall \bs \in \mathcal{S}.\]

ALPs approximate the high-dimensional MDP value function $V(\bs)$  (\citealp{schweitzer_generalized_1985,de_farias_linear_2003}) using a linear combination of basis functions. We construct the ALP value function approximation using an intercept $\tau$ and $B$ basis functions $\phi_b : \Space\mapsto \R$, $b = 1,\ldots, B$, that is, $V(\bs) \approx \tau + \sum_{b = 1}^B \alpsoln_b\phi_b(\bs)$, where $\theta := (\alpsoln_1,\ldots, \alpsoln_B) \in \mathbb{R}^B$ is the basis function weight vector. It is common to require that the pair $(\tau, \theta)$ belongs to a compact set $\mathcal{X}$. The VFA weights are computed by solving
{\small
	\begin{eqnarray*}
		\nonumber
		&\displaystyle \max_{(\tau,\alpsoln) \in\mathcal{X}} & \tau + \sum_{b = 1}^B\alpsoln_b\left[\phi_b(\bs^0)\right]\\\nonumber
		&\text{s.t.} &(1-\gamma)\tau + \sum_{b = 1}^B \alpsoln_b \left(\phi_b(\bs) - \gamma \mathbb{E}\left[ \phi_b(f(\bs,a)) \right]\right)-c(\bs,a)\leq 0,\hspace{10pt} \forall (\bs, a)\in \Space\times\action.\nonumber
\end{eqnarray*}}
\normalsize
The feasibility of the ALP constraints is important because it ensures that the objective function of a feasible solution provides a lower bound on the optimal policy value, which can be used to assess the suboptimality of heuristic policies (see, e.g., Proposition 4 in \citealp{adelman2008relaxations}). Thus, in principle, methods to solve ALP would benefit from emphasizing feasibility as we do in SFLS. 

Since the linear program above is semi-infinite, constraint sampling is a popular strategy to approach its solution and obtain a high-probability feasible solution \citep{farias_constraint_2004}. Specifically, suppose we sample $m$ state-action pairs $(\bs_i,a_i), i = 1,\ldots,m$. The ALP with constraints corresponding to these samples takes the form of \eqref{eq:gco}: 
{\small
	\begin{eqnarray}
	\label{ALP}
	&\displaystyle \max_{\bx=(\tau,\alpsoln) \in\mathcal{X}} & f_0(\bx):=\tau + \sum_{b = 1}^B\alpsoln_b\left[\phi_b(\bs^0)\right]\\\nonumber
	&\text{s.t.} &f_i(\bx):=(1-\gamma)\tau + \sum_{b = 1}^B \alpsoln_b \left(\phi_b(\bs_i) - \gamma \mathbb{E}\left[ \phi_b(f(\bs_i,a_i)) \right]\right)-c(\bs_i,a_i)\leq 0,\quad i=1,2,\dots,m. \nonumber
	\end{eqnarray}}
\normalsize
We solve this linear program in our experiments. 

Following \citet{lin2019revisiting}, we constructed instances with $I = 2$ and $J=2$, chose $P_G$ to be a truncated normal in the interval $[0,10]$ with mean $5$ and the standard deviation 2, and fixed $c_p$, $c_l$, $\bar a$, $l_s$, $\gamma$, and $\bs^0$ equal to $20$, $100$, $10$, $-10$, $0.95$, and $(5,0,0)$, respectively. We experimented with three instances based on the triple $(c_h,c_d,c_b)$ being equal to $(2,10,10)$, $(5,10,8)$, and $(2,5,10)$. We employed eighteen basis functions ($B=18$): $z_0$, $z_1$, $q_1$, and  $\{(z_0 - \nu)_+, (z_0 + z_1 - 2\nu)_+, (z_0 + z_1 + q_1 - 3\nu)_+, (2\nu - z_0 - z_1 - q_1)_+, (\nu - z_1 - q_1)_+| \nu \in \{\mathbb{E}[G],G^{0.25},G^{0.5}\}\}$, where $G^{0.25}$ and $G^{0.5}$ are the 25-th and 50-th quartiles of the demand distribution. The domain for the basis function weights $\mathcal{X}$ was taken to be the box $[0,3000]\times [-5,5]^{B}$. We chose $m$ as 500.


In all methods, we use the initial solution $\tilde\bx=(\tilde\tau,\tilde\alpsoln)$ with $\tilde\tau=\min_{i=1,\dots,m}\frac{c(s_i,a_i)}{1-\gamma}$ and $\tilde\alpsoln=\mathbf{0}$ which is feasible for \eqref{ALP}.
Our SFLS implementation uses $r^{(0)}=f_0(\tilde\bx)>f^*$, $\theta=1.1$, and the step length rule $\gamma_t=5/\sqrt{t+1}$. We do not report results for DFLS because, as alluded to in \S\ref{subsec:AlgoImpl}, obtaining a good deterministic approximation using SAAs is non-trivial for the perishable inventory control problem. We use a mini-batch technique with a batch size of $100$ to construct stochastic estimates of the gradients and function values of $f_i$, $i = 0,\ldots,m$ in both SFLS and YNW. 
In SFLS, we choose the numbers of inner (i.e. $T$) and outer iterations to be $200$ and $100$, respectively, which leads to $20,000$ stochastic gradient steps (inner iterations) in total. Hence, we choose the total number of iterations in YNW as $T=20,000$ so that both methods evaluate the same number of stochastic gradients in total which lead to similar runtime (about 1100 seconds). 

\begin{figure}[t]
	\begin{tabular}[h]{@{}c|ccc@{}}
		& $(c_h,c_d,c_b)=(2,10,10)$ &$(c_h,c_d,c_b)=(5,10,8)$&$(c_h,c_d,c_b)=(2,5,10)$\\
		\hline \\
		\raisebox{10ex}{\small{\rotatebox[origin=c]{90}{$f_0(x)-f^*$}}}
		& \includegraphics[width=0.31\textwidth]{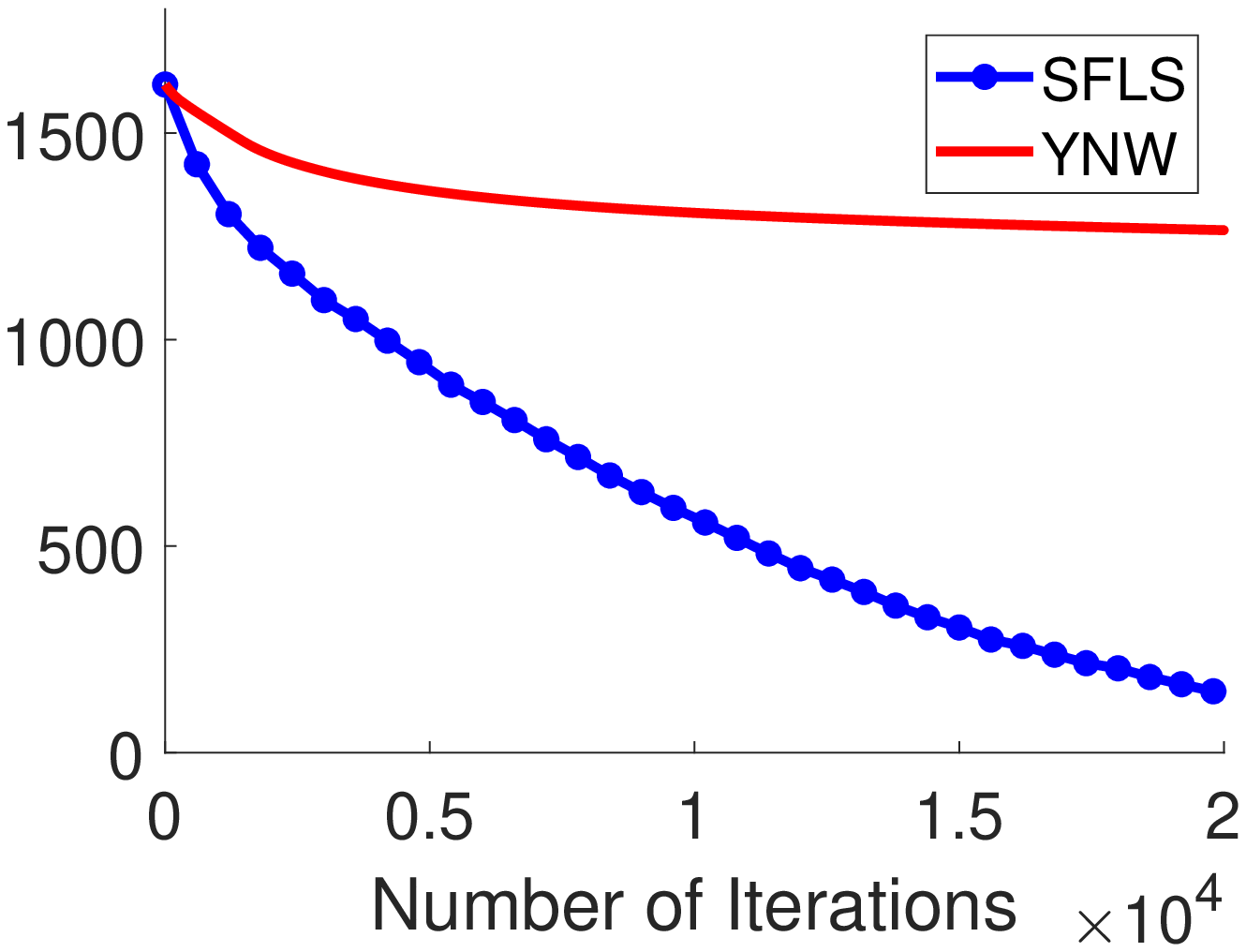}
		& \includegraphics[width=0.31\textwidth]{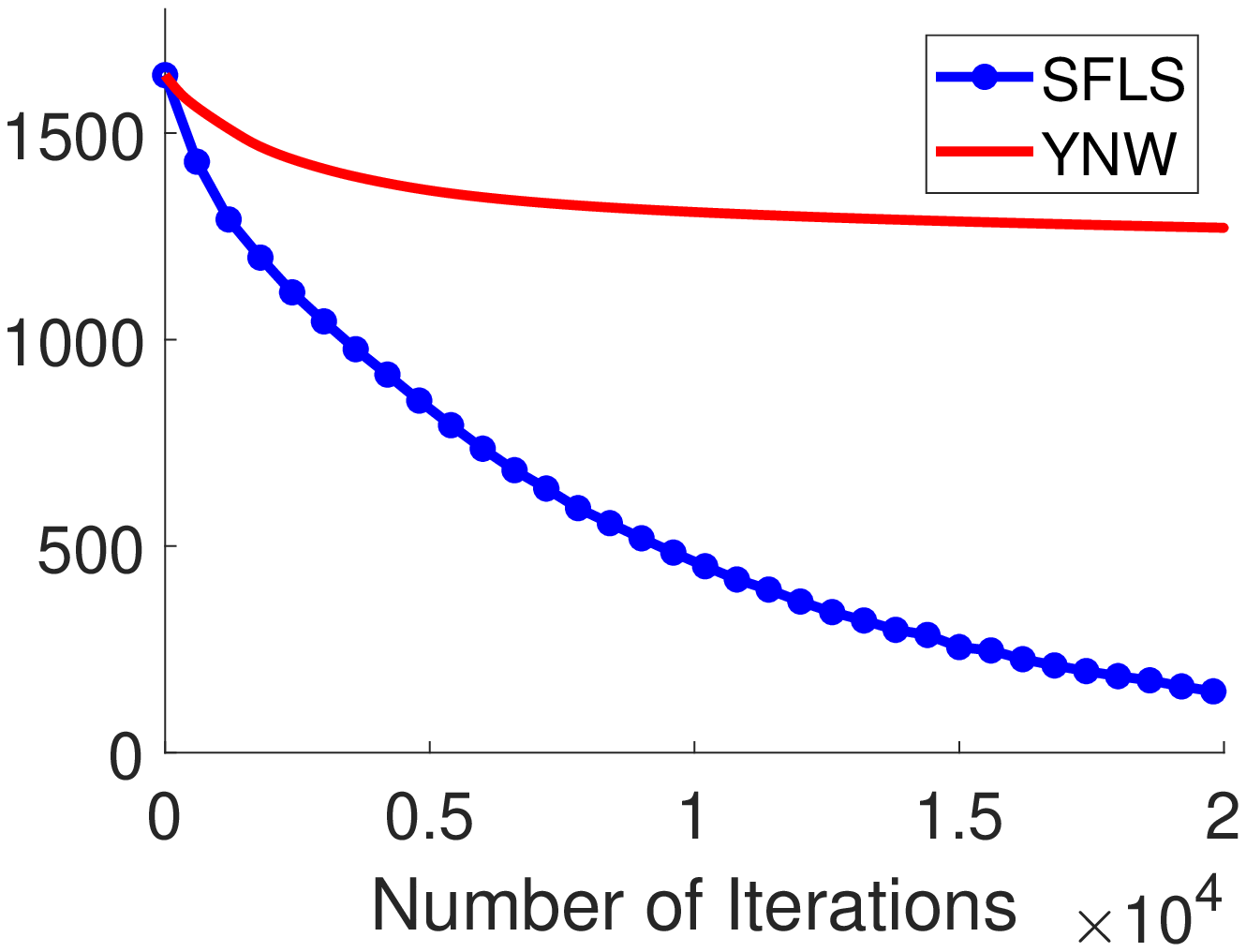}
		& \includegraphics[width=0.31\textwidth]{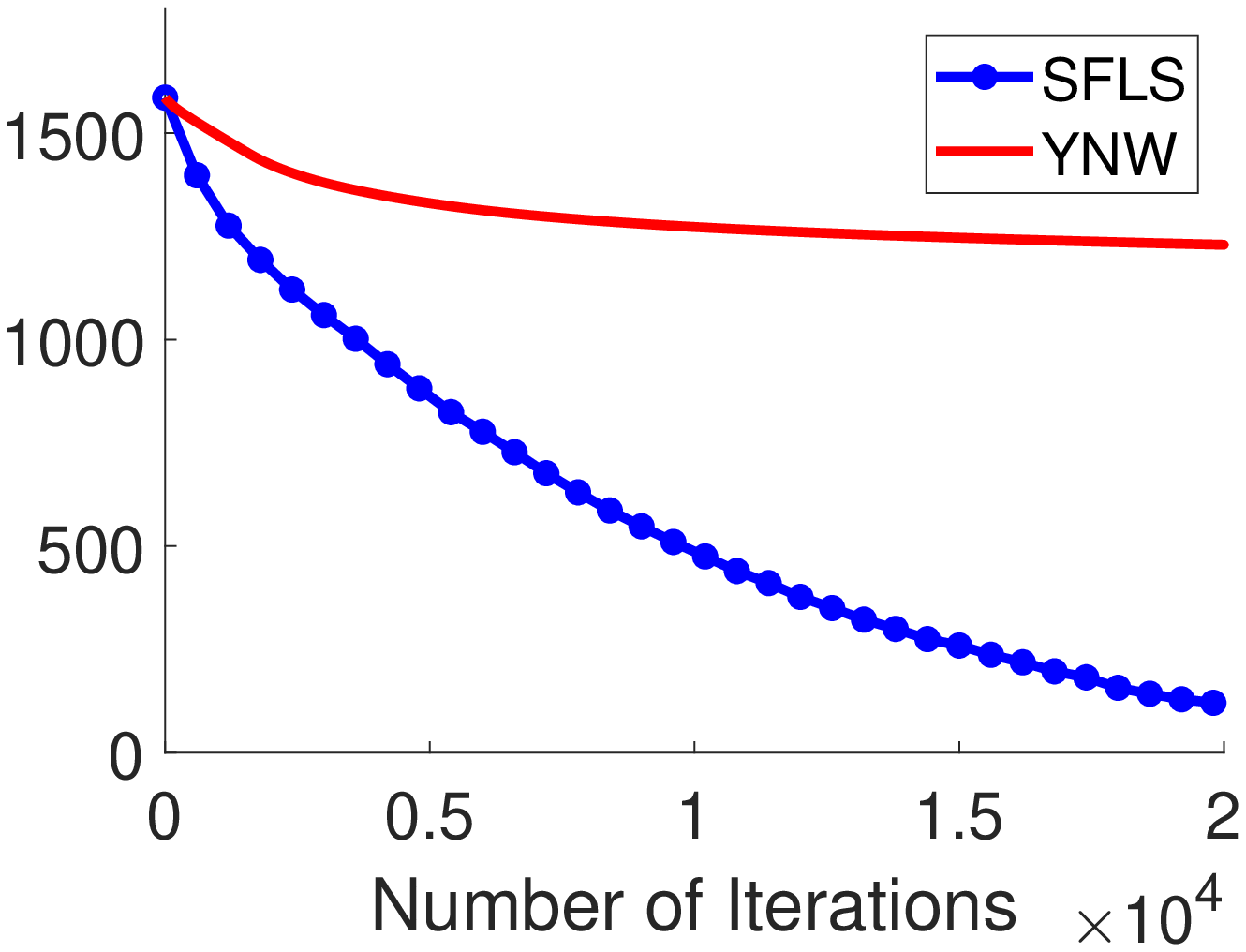}\\
		\raisebox{10ex}{\small{\rotatebox[origin=c]{90}{$\max\limits_{i=1,2,\dots,m}\{f_i(\bx)-r_i\}$}}}
		& \includegraphics[width=0.31\textwidth]{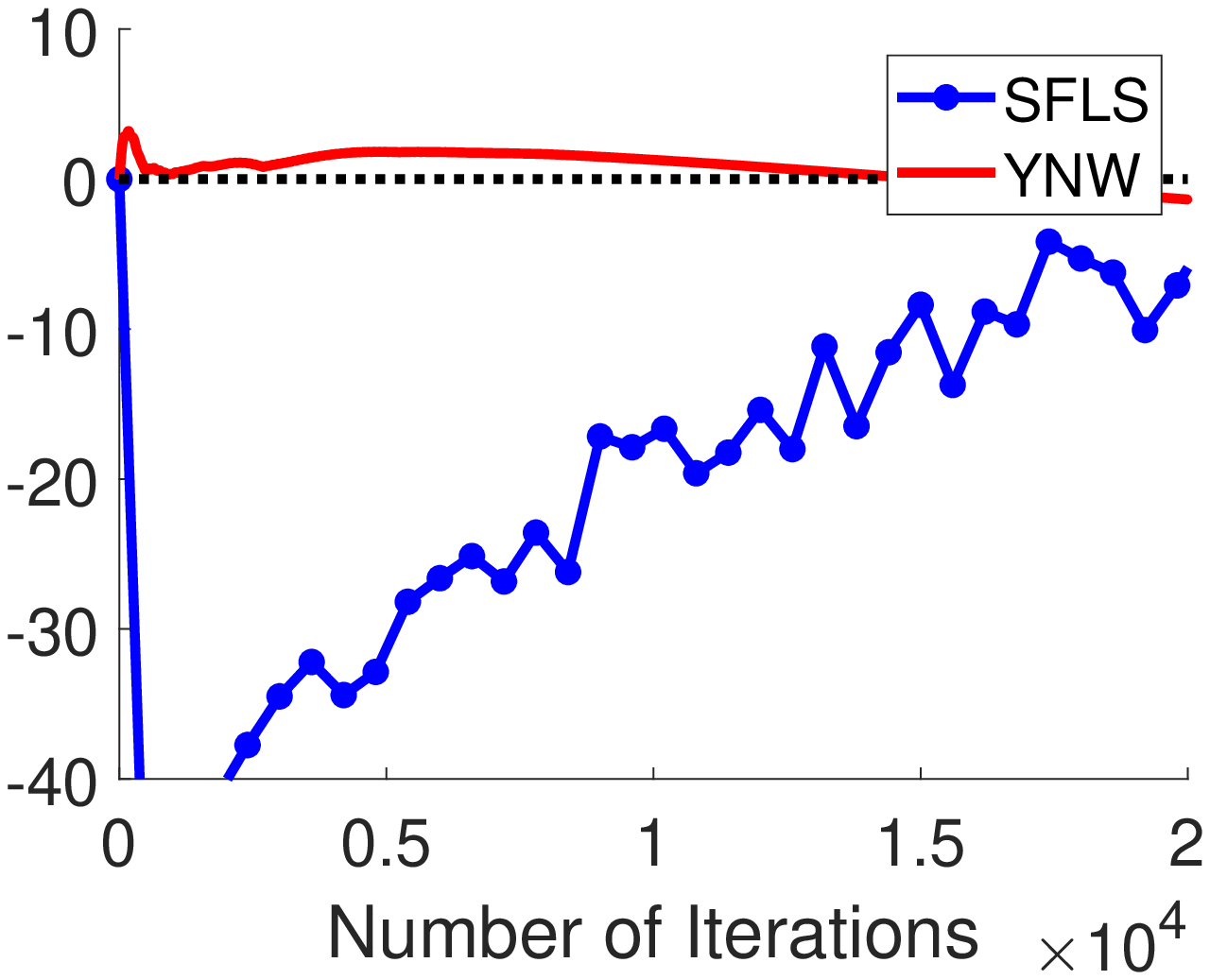}
		& \includegraphics[width=0.31\textwidth]{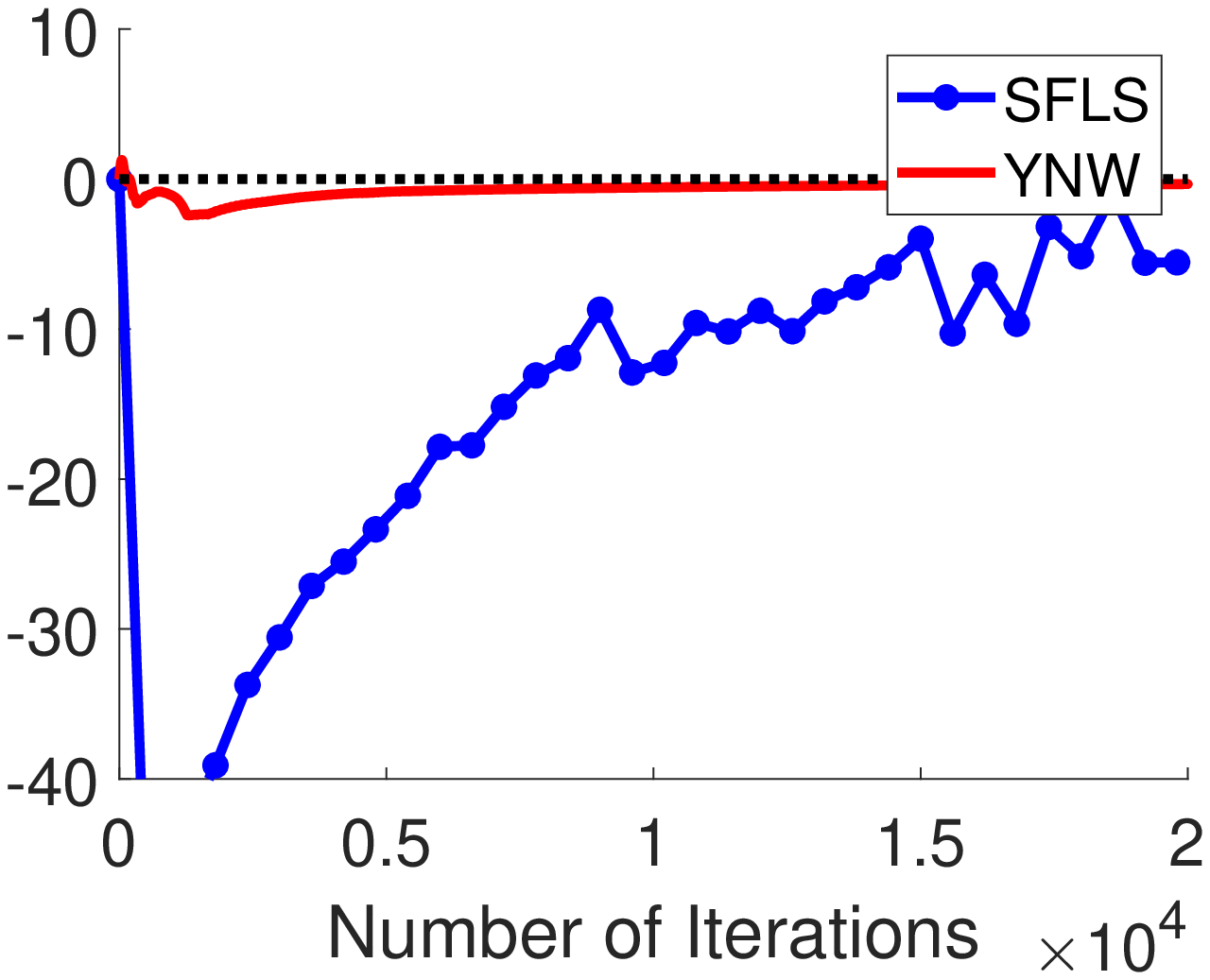}
		& \includegraphics[width=0.31\textwidth]{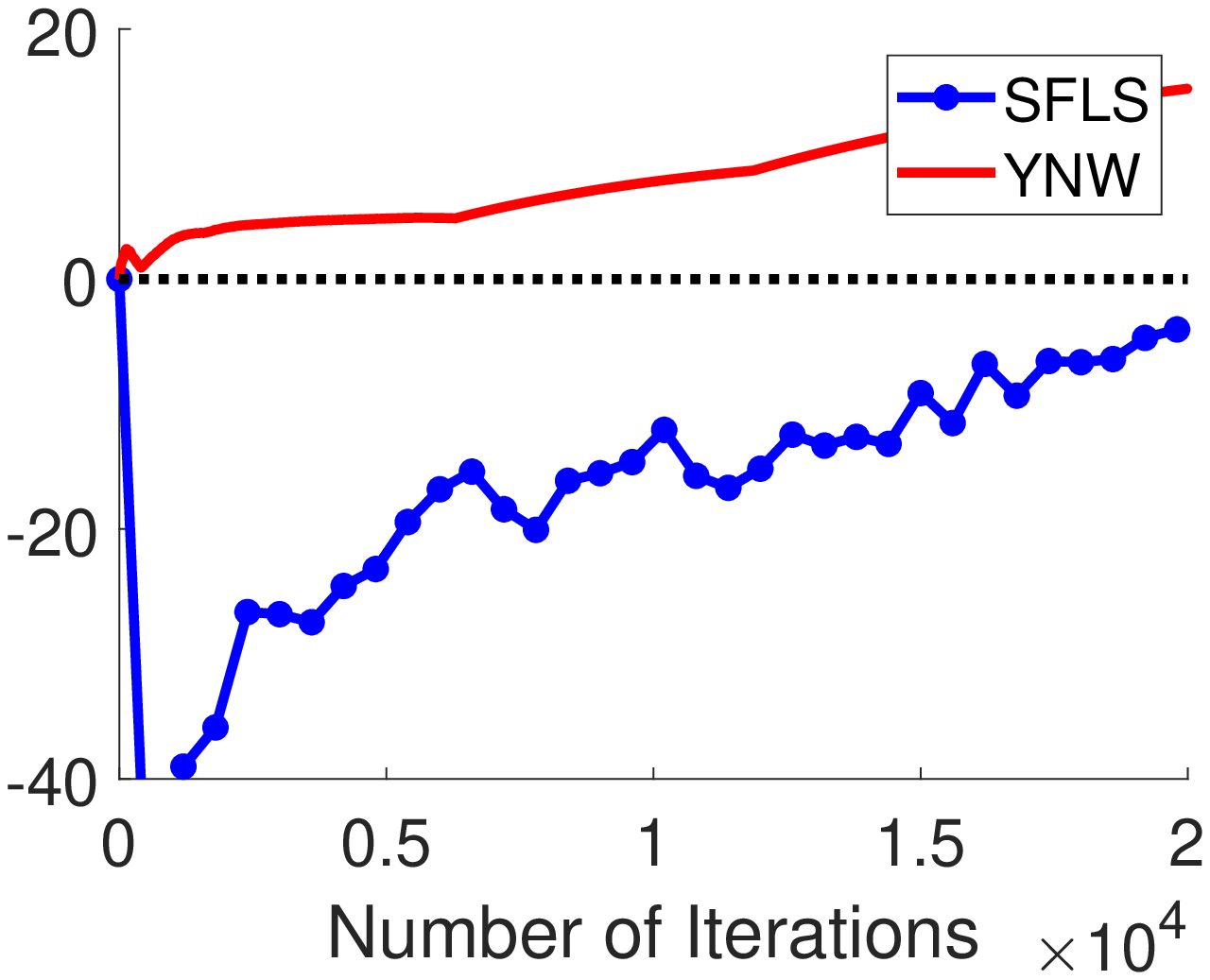}\\
	\end{tabular}
	\vspace{2ex}
	\caption{Performance of SFLS and YNW for solving approximate linear programs arising in perishable inventory control.
	}
	\label{fig:perishable}
\end{figure}

Figure \ref{fig:perishable} displays the performance of SFLS and YNW. The $y$-axes of the top subfigures report the optimality gap $f_0(\bx)-f^*$ while these axes in the bottom subfigures show the feasibility of solutions by plotting $\max_{i=1,2,\dots,m}\{f_i(\bx)-r_i\}$.  Here, $f_i(\bx)$ for $i=1,2,\dots,m$ are calculated by approximating the expectations in their definitions in \eqref{ALP} with $10,000$ samples of demand $G$. The optimal value $f^*$ is approximated by the objective value found by a separated run of SFLS with sufficient iterations ($400$ outer and $500$ inner iterations). We track these measures as a function of the number of iterations performed by each algorithm in the $x$-axis. To indicate the values of $f_0(x)-f^*$ and $\max_{i=1,2,\dots,m}\{f_i(\bx)-r_i\}$ corresponding to the high probability feasible solutions maintained at each SFLS (outer) iteration we use line markers in Figure \ref{fig:perishable}. The YNW curves have no line markers as there are no outer iterations ensuring feasibility. SFLS finds a feasible solution quickly and maintain a relatively large constraint slack but YNW does not always ensure feasibility. SFLS also reduces the suboptimality of solutions faster, suggesting that SFLS is able to balance optimality and feasibility well on these instances. 
\looseness=-1

\subsection{Multi-class Neyman-Pearson classification}\label{sec:NumClassification}

Another application that gives rise to \eqref{eq:gco} is Neyman-Pearson classification. In multi-class classification, there exist $m$ classes of data, where $\psi_i$, $i=1,2,\dots,m$, denotes a random variable defined using the distribution of data points associated with the $i$-th class. To classify a data point $\psi_i$ to one of the $m$ classes, we rely on the same number of linear models $\bx_i$, $i=1,2,\dots,m$. The predicted class for $\psi$ is 
$\argmax_{i=1,2,\dots,m}\bx_i^\top \psi$. High classification accuracy in this scheme requires $\bx_i^\top \psi_i-\bx_l^\top \psi_i$ with $i\neq l$ to be large and positive ~\citep{weston1998multi,crammer2002learnability}, that is, the classifiers have discriminatory power. Minimizing the expected loss $\E\left[\phi(\bx_i^\top \psi_i-\bx_l^\top \psi_i)\right]$ is one approach to promote this goal, where $\phi$ is a non-increasing convex loss function and $\E$ is expectation taken over $\psi_i$. 

Suppose misclassifying $\psi_i$ has a cost that depends on $i$ but not on the predicted class. We propose a model that prioritizes classes with relatively higher misclassification costs using constraints and simultaneously trains the set of $m$ linear models by solving\looseness=-1
\begin{eqnarray}
\label{eq:NPclassification_multi}
\min_{\substack{\|\bx_i\|_2\leq\lambda,\\ \forall i = 1,2,\ldots,m}} \sum_{l\neq 1}\E[\phi(\bx_1^\top \psi_1-\bx_l^\top \psi_1)],~\text{s.t.}~\sum_{l\neq i}\E[\phi(\bx_i^\top \psi_i-\bx_l^\top \psi_i)]\leq r_i,\quad i=2,3,\dots,m,
\end{eqnarray}
where it is assumed (without loss of generality) that class $1$ has the highest misclassification cost and the value of $r_i$ is chosen to capture the misclassification cost of class $i$. Here $\lambda$ is a regularization parameter. This formulation can be easily extended to handle the case where the mis-classification cost depends on both the true and predicted classes. Indeed, \eqref{eq:NPclassification_multi} is of the form \eqref{eq:gco}. Infeasible solutions may result in large misclassification costs for some classes, which is undesirable, and creates a need for methods that emphasize feasibility.


We created test instances using the multi-class classification LIBSVM datasets \emph{connect-4}, \emph{covtype}, and \emph{news20} from \citet{Chang2019}. We selected these instances as their size still allows us to run DFLS in the manner discussed in \S\ref{subsec:AlgoImpl}. We summarize in Table~\ref{tab:datasets} the number of classes, the number of data points in each class, and the number of features in these four datasets.
We chose the loss function~\eqref{eq:NPclassification_multi} to be the hinge loss $\phi(z)=(1-z)_+$. Let $\psi_i$ follow the empirical distribution over the dataset of class $i$ for $i=1,2,\dots,m$, which implies that all the expectations in~\eqref{eq:NPclassification_multi} become finite-sample averages over data classes. 
We set the parameters $\lambda=5$ and $r_i=m-1$ for $i=2,\dots,m$. 

In all methods, the solution $\tilde\bx=\mathbf{0}$ is used as the initial solution and it is feasible for \eqref{eq:NPclassification_multi}. To apply SFLS and DFLS, we chose $T=100$, $r^{(0)}=m$, and $\theta=1.1$ across all datasets. Note that 
$r^{(0)}=m>m-1=f_0(\mathbf{0})\geq f^*$ for \eqref{eq:NPclassification_multi}. In DFLS, we solve subproblems via standard subgradient descent method. In SFLS and DFLS, we choose step size $\gamma_t=0.05/\sqrt{t+1}$ for \emph{connect-4} and \emph{covtype} and choose $\gamma_t=1/\sqrt{t+1}$ for \emph{news20}.  Both SFLS and YNW employed a mini-batch size of $1000$ to construct the stochastic gradients and the objective values. We chose the number of iterations in YNW so that its total number of data passes is $200$ for \emph{connect-4} and \emph{news20} and $100$ for \emph{covtype}. Then, we also terminated SFLS and DFLS when the total data passes they performed exceed YNW.

\begin{table}[htp]
	\begin{center}
		\begin{tabular}{|c|c|c|c|}
			\hline
			Dataset	&Number of  classes & Number of instances & Number of  features \\
			\hline
			{connect-4}   &3 & 67557 &  126\\
			{covtype}& 7 & 581012 &  54   \\
			{news20} & 20 & 15935& 62061   \\
			\hline
		\end{tabular}
		\caption{Characteristics of multi-class classification datasets
			from LIBSVM library}
		\label{tab:datasets}
	\end{center}
\end{table}

\begin{figure}[t]
	\begin{tabular}[h]{@{}c|ccc@{}}
		& connect-4 & covtype  &news20\\
		\hline \\
		
		\raisebox{10ex}{\small{\rotatebox[origin=c]{90}{$f_0(x)-f^*$}}}
		& \includegraphics[width=0.31\textwidth]{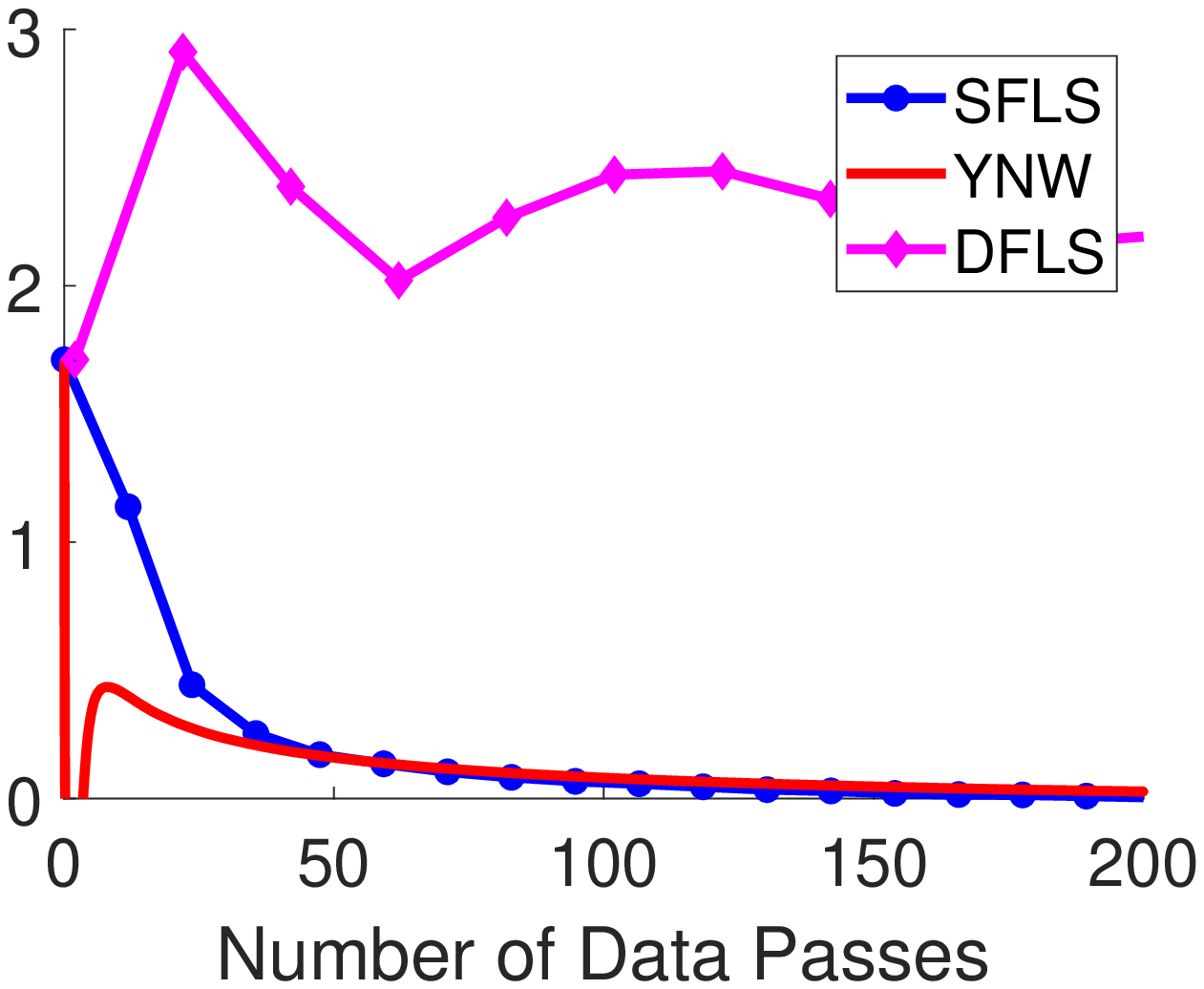}
		& \includegraphics[width=0.31\textwidth]{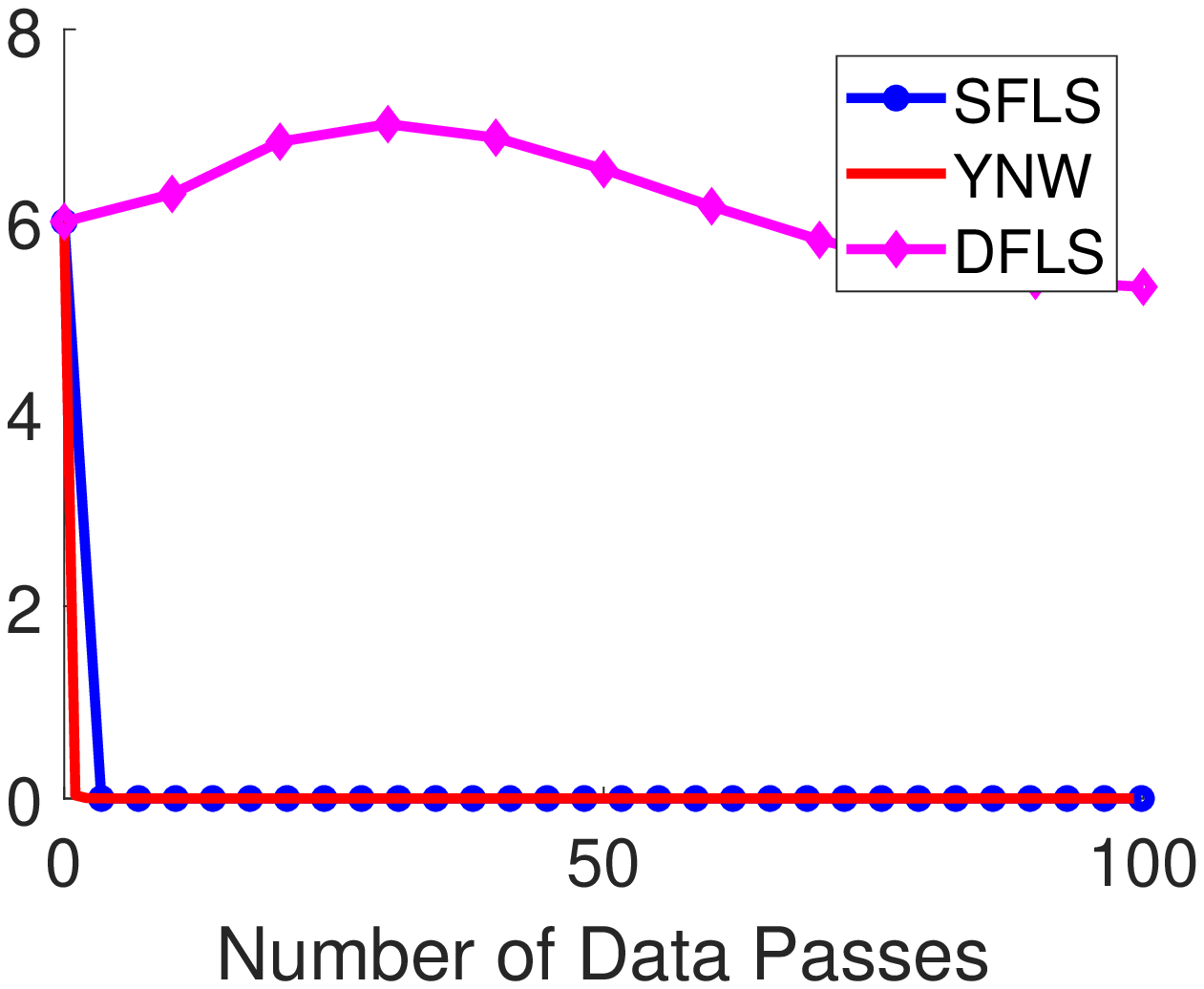}
		& \includegraphics[width=0.31\textwidth]{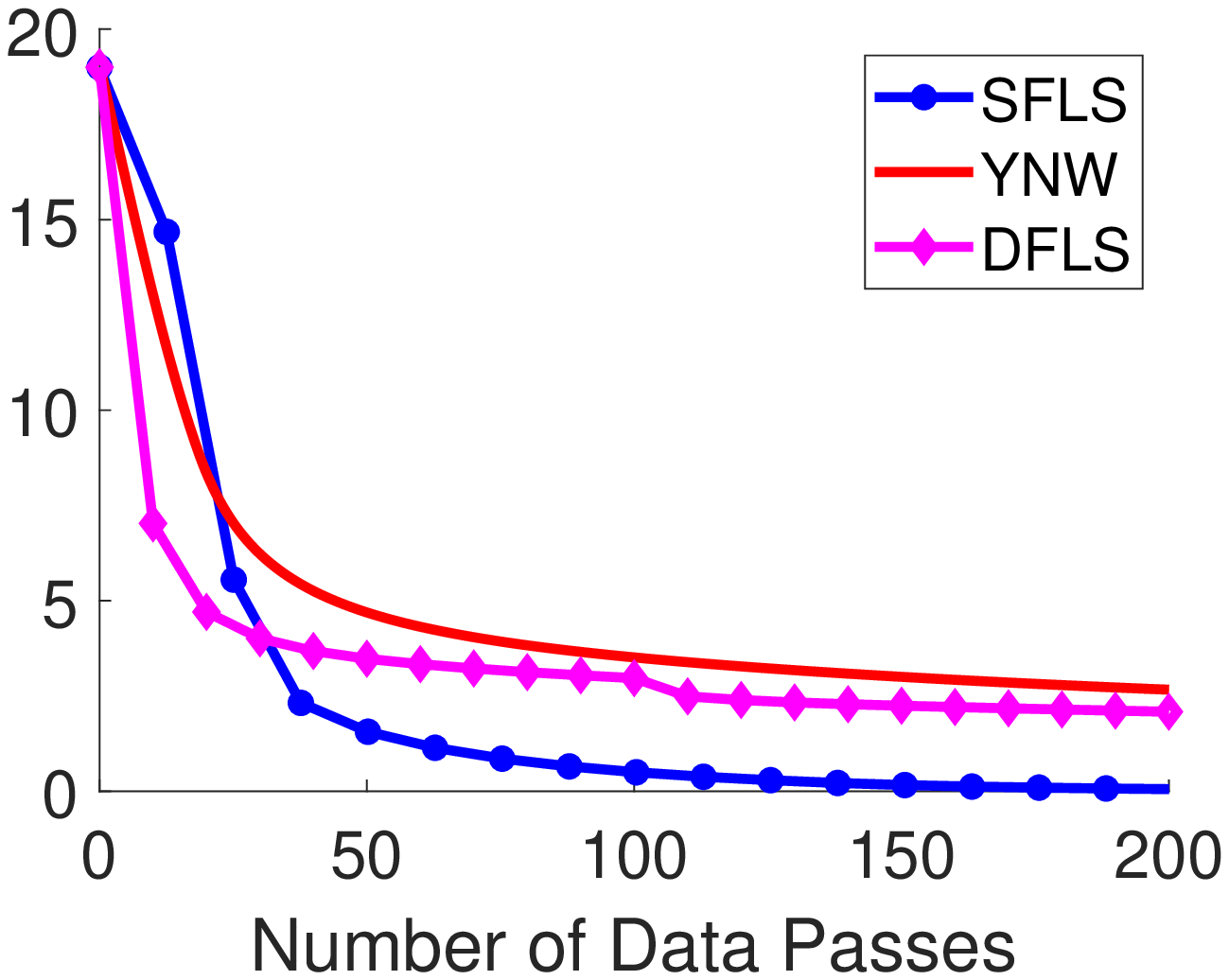}\\
		
		\raisebox{10ex}{\small{\rotatebox[origin=c]{90}{$\max\limits_{i=1,2,\dots,m}\{f_i(\bx)-r_i\}$}}}
		& \includegraphics[width=0.31\textwidth]{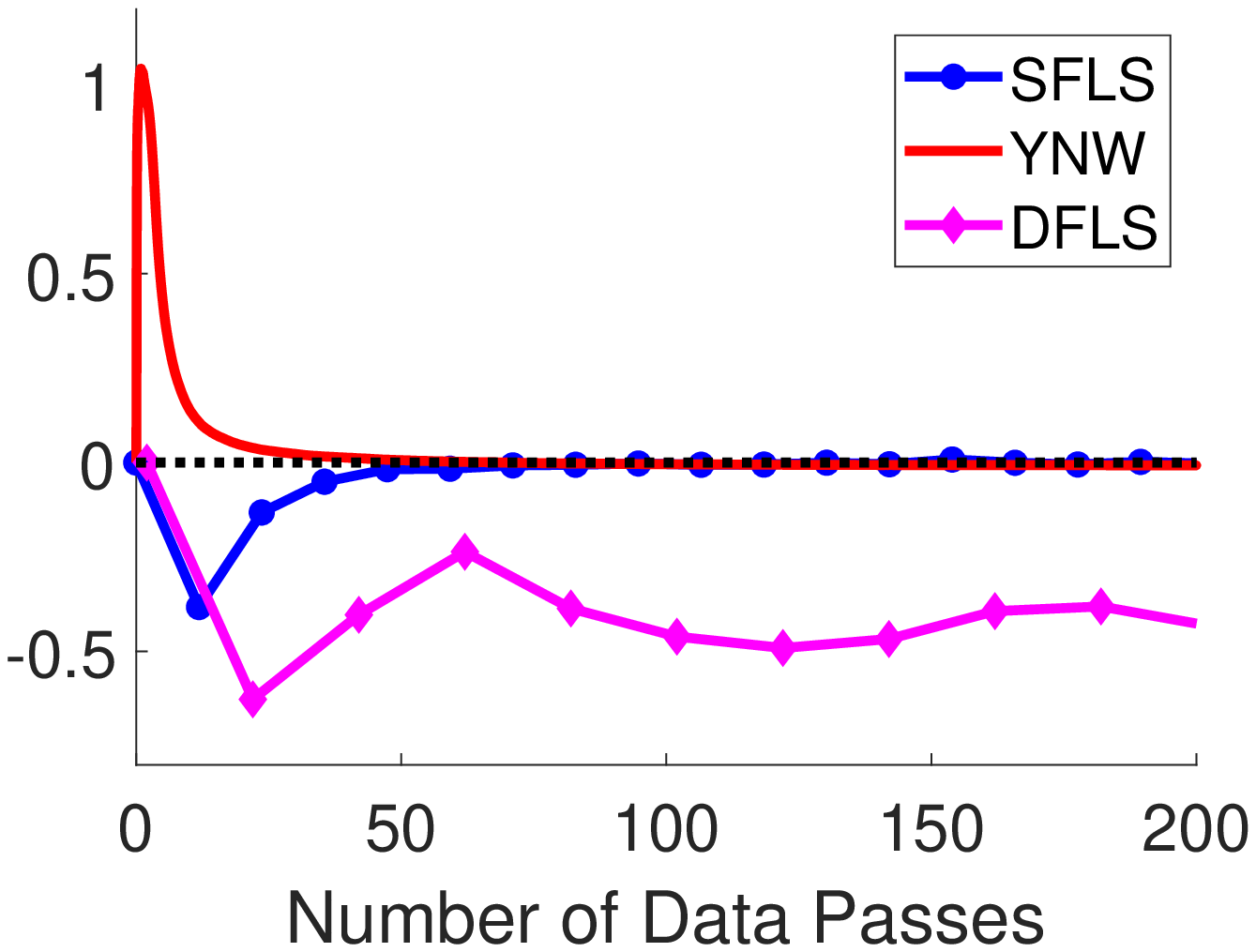}
		& \includegraphics[width=0.31\textwidth]{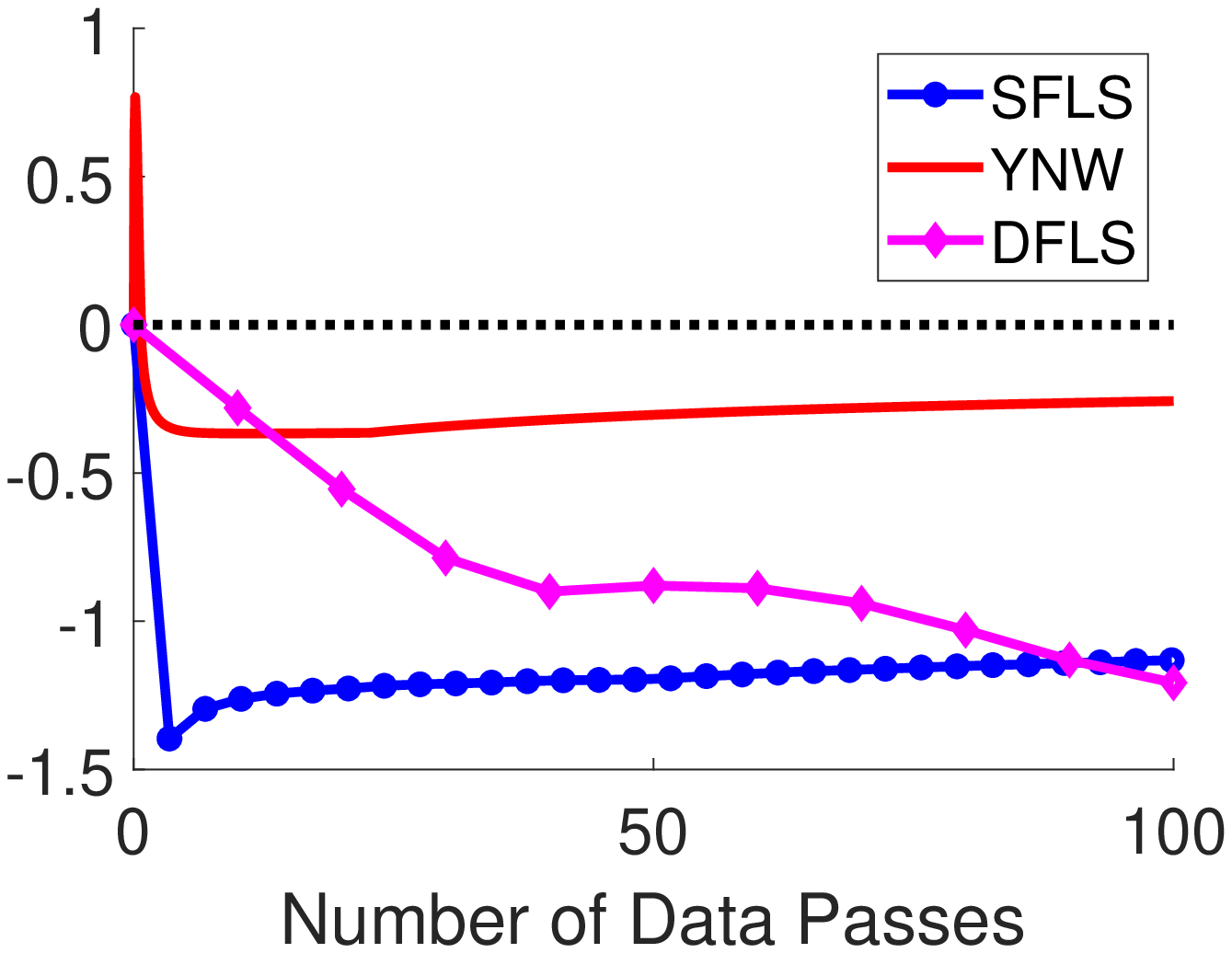}
		& \includegraphics[width=0.31\textwidth]{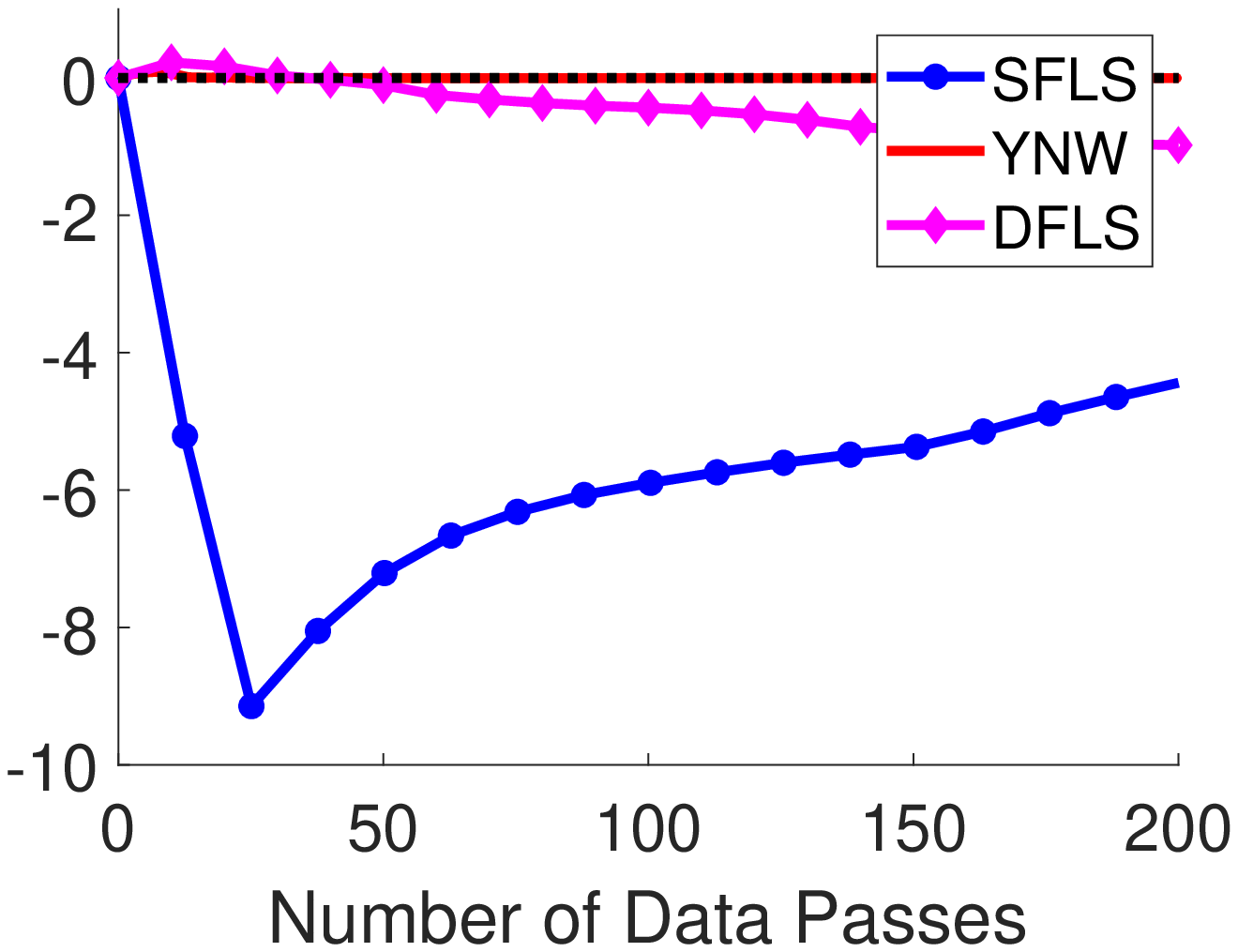}\\
	\end{tabular}
	\vspace{2ex}
	\caption{Performance of SFLS, YNW, and DFLS on the  multi-class Neyman-Pearson classification problem. 
	}
	\label{fig:connect4}
\end{figure}

Figure \ref{fig:connect4} displays the performance of each method. The $y$-axes of the first row reports the term $f_0(\bx)-f^*$, that is, it focuses on optimality, while this axis in the second row shows the feasibility of solutions by plotting $\max_{i=1,2,\dots,m}\{f_i(\bx)-r_i\}$. We track these measures as a function of the number of equivalent data passes performed by each algorithm in the $x$-axis, where a data pass involves going over the number of data points equal to the size of the training data. This is possible since the expectations in our instances are over discrete random variables. Tracking data passes allows us to assess algorithms in terms of data complexity. Similar to Figure \ref{fig:perishable}, we uses line markers to indicate the values of $f_0(x)-f^*$ and $\max\limits_{i=1,2,\dots,m}\{f_i(\bx)-r_i\}$ corresponding to the solutions maintained at each SFLS outer iteration, while YNW has no line marker since it does not maintain feasibility. Since DFLS needs two data passes in each inner iteration, it can only perform one or two outer iterations with the number of data passes in Figure \ref{fig:perishable}. Hence, for a better visualization, we use line markers to also indicate the inner iterations of DFLS instead of only outer iterations. In this figure, $f^*$ is approximately by the objective value returned by DFLS after a sufficient number of data passes (i.e. at least $5000$ data passes with $2T$ inner iterations.) 
\looseness=-1

On the connect-4 data set, SFLS maintains feasibility and reduces the optimaliy gap quite rapidly after a few data passes. Interestingly, despite providing an initial feasible solution, YNW decreases the optimality gap at the beginning by moving to a highly infeasible solution. The peformance of both methods on the covtype data are comparable. On the news20 data set, SFLS provides feasible solutions with smaller optimality gaps sooner than the benchmark method. The comparison of SFLS and YNW highlights the advantage of SFLS in terms of feasibility. Specifically, efficient methods that do not emphasize feasibility could lead to highly infeasible solutions if terminated prematurely (e.g., the connect-4 dataset). 

DFLS also maintains a feasible solution path on all the datasets, as expected. However, its optimality gap reduces at a much slower rate with the number of data passes compared to SFLS because it uses deterministic subgradients based on the entire data set. These results thus underscore the importance of developing methods, such as SFLS, with low data complexity to balance optimality and feasibility.

\looseness=-1

\subsection{Learning with Fairness Constraints}\label{appl:LearnDataConst}
We consider learning a classifier with fairness constraints. Other examples include training predictive models with constraints on coverage rates, churn rates, and stability. Please see \citet{goh2016satisfying} for further motivation and a non-convex formulation. Here we provide a convex formulation for these problems, which can be viewed as a tractable relaxation of the version in \citet{goh2016satisfying} that admits the SOEC structure \eqref{eq:gco}. 

Suppose $(\ba,b)$ is a data point from a distribution $\mathcal{D}$, where $\ba$ is a feature vector and $b\in\{1,-1\}$ is the class label. Let $\mathcal{D}_M$ and $\mathcal{D}_F$ denote two different distributions of features (that are not necessarily labeled), which may represent male and female individuals. The goal is to train a classifier $\ba^\top\bx$ that minimizes classification loss. The correct classification of data vector $\ba$ implies that $b \ba^\top \bx > 0$. One can train such a classifier subject to fairness constraints by solving
\begin{eqnarray}
\label{eq:DDCproblem}
\min_{\|\bx\|_2\leq\lambda} &&\E_{(\ba,b)\sim\mathcal{D}}[\phi(-b\ba^\top\bx)]\\\nonumber
\text{s.t.}&&\E_{\ba\sim\mathcal{D}_M}[\sigma(\ba^\top\bx)]\leq \E_{\ba\sim\mathcal{D}_F}[\sigma(\ba^\top\bx)]/\kappa,\\\nonumber
&& \E_{\ba\sim\mathcal{D}_F}[\sigma(\ba^\top\bx)]\leq \E_{\ba\sim\mathcal{D}_M}[\sigma(\ba^\top\bx)]/\kappa,
\end{eqnarray}
where  $\lambda$ is a regularization parameter, $\kappa\in(0,1]$ is a constant, $\phi$ is a non-increasing loss function, 
$$
\sigma(z)=\max\{0,\min\{1,\{0.5+z\}\},
$$ and $\sigma(\ba^\top\bx)\in[0,1]$ represents the probability of the (random) classifier $\bx$ predicting $\ba$ as positive. Therefore, $\E_{\ba\sim\mathcal{D}_M}[\sigma(\ba^\top\bx)]$ and $\E_{\ba\sim\mathcal{D}_M}[\sigma(\ba^\top\bx)]$ represent the percentages of instances in $\mathcal{D}_M$ and $\mathcal{D}_F$ predicted as positive, respectively. The first constraint guarantees that the percentage of the positively predicted instances in $\mathcal{D}_F$ is at least a $\kappa$ fraction of that in $\mathcal{D}_M$. The second constraint has similar interpretation. An analogous model was considered in \citet{goh2016satisfying} but it involves non-convex constraints. 

Observing that 
$\sigma(\ba^\top\bx)=1-\sigma(-\ba^\top\bx)$, we can reformulate the first constraint as $\E_{\ba\sim\mathcal{D}_M}[\sigma(\ba^\top\bx)]+\E_{\ba\sim\mathcal{D}_F}[\sigma(-\ba^\top\bx)]/\kappa\leq 1/\kappa$ and approximate $\sigma$ by $\max\{0,0.5+z\}=(0.5+z)_+$ so that we obtain a convex constraint  $\E_{\ba\sim\mathcal{D}_M}[(\ba^\top\bx+0.5)_+]+\E_{\ba\sim\mathcal{D}_F}[(-\ba^\top\bx+0.5)_+]/\kappa\leq 1/\kappa$. Applying an analogous convex approximation to the second constraint, we obtain the following convex formulation for training a classifier subject to fairness constraints:
\begin{eqnarray*}
	\min_{\|\bx\|_2\leq\lambda} &&\E_{(\ba,b)\sim\mathcal{D}}[\phi(-b\ba^\top\bx)]\\
	\text{s.t.}&&\E_{\ba\sim\mathcal{D}_M}[(\ba^\top\bx+0.5)_+]+\E_{\ba\sim\mathcal{D}_F}[(-\ba^\top\bx+0.5)_+]/\kappa\leq 1/\kappa,\\
	&&\E_{\ba\sim\mathcal{D}_F}[(\ba^\top\bx+0.5)_+]+\E_{\ba\sim\mathcal{D}_M}[(-\ba^\top\bx+0.5)_+]/\kappa\leq 1/\kappa.
\end{eqnarray*}
The left hand side of the first constraint will be large if the classifier $\bx$ is not ``fair'', that is, it makes $\ba^\top\bx$ very negative for most of $\ba$ from $\mathcal{D}_M$ but very positive for most of $\ba$ from $\mathcal{D}_F$. Similarly, the left hand side of the second constraint will be large if the model $\bx$ makes $\ba^\top\bx$ very positive for most of $\ba$ from $\mathcal{D}_F$ but very negative for most of $\ba$ from $\mathcal{D}_M$. Choosing an appropriate $\kappa$ ensures that the obtained model is fair to both  $\mathcal{D}_M$ and $\mathcal{D}_F$. Indeed, a solution that violates constraints in this formulation translates to a classifier that discriminates against one of the two classes.

\begin{figure}[h!]
	\centering
	\begin{tabular}[h]{@{}c|cc@{}}
		&a9a & LoanStats \\
		\hline \\
		\raisebox{10ex}{\small{\rotatebox[origin=c]{90}{$f(x)-f^*$}}}
		& \includegraphics[width=0.40\textwidth]{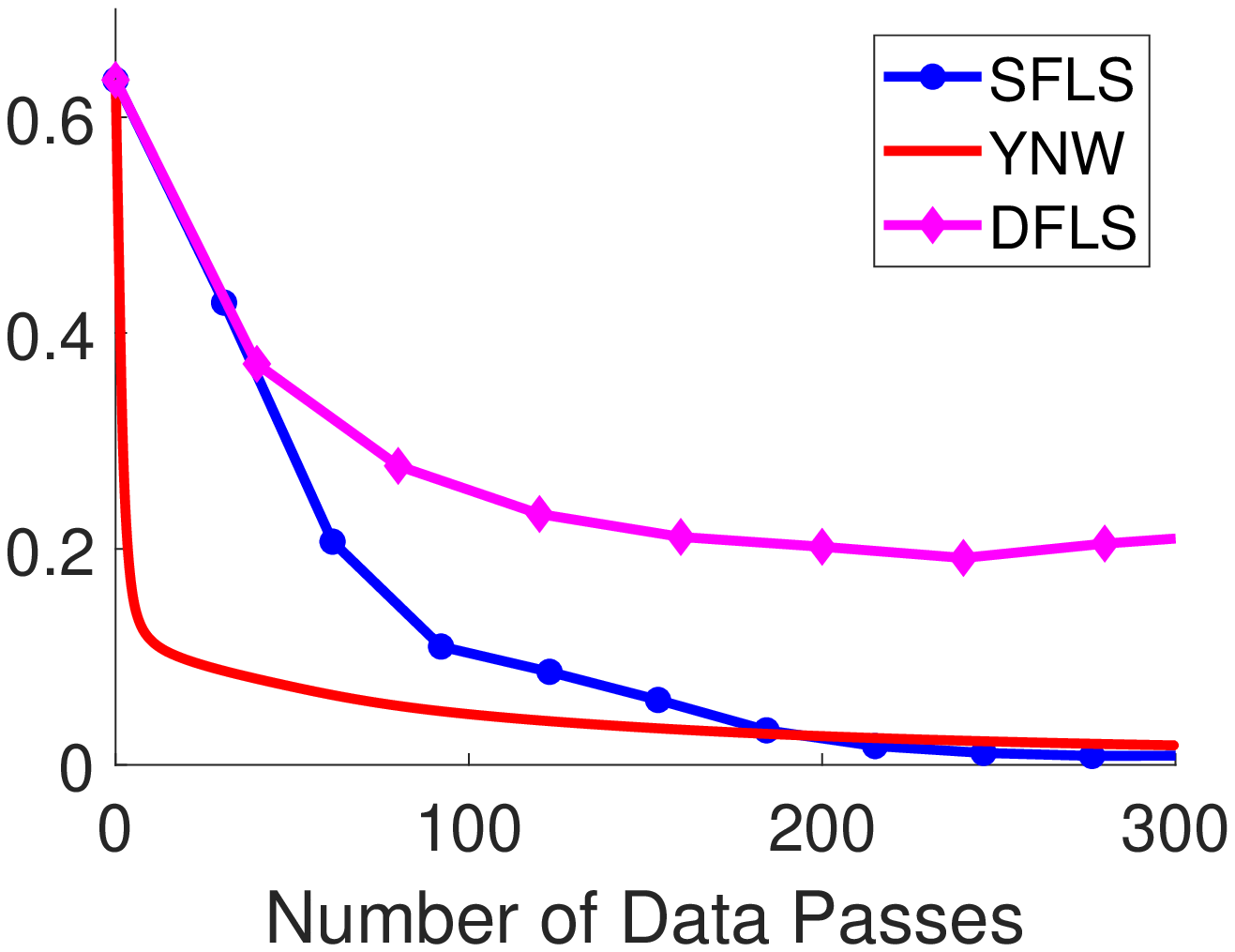}
		& \includegraphics[width=0.40\textwidth]{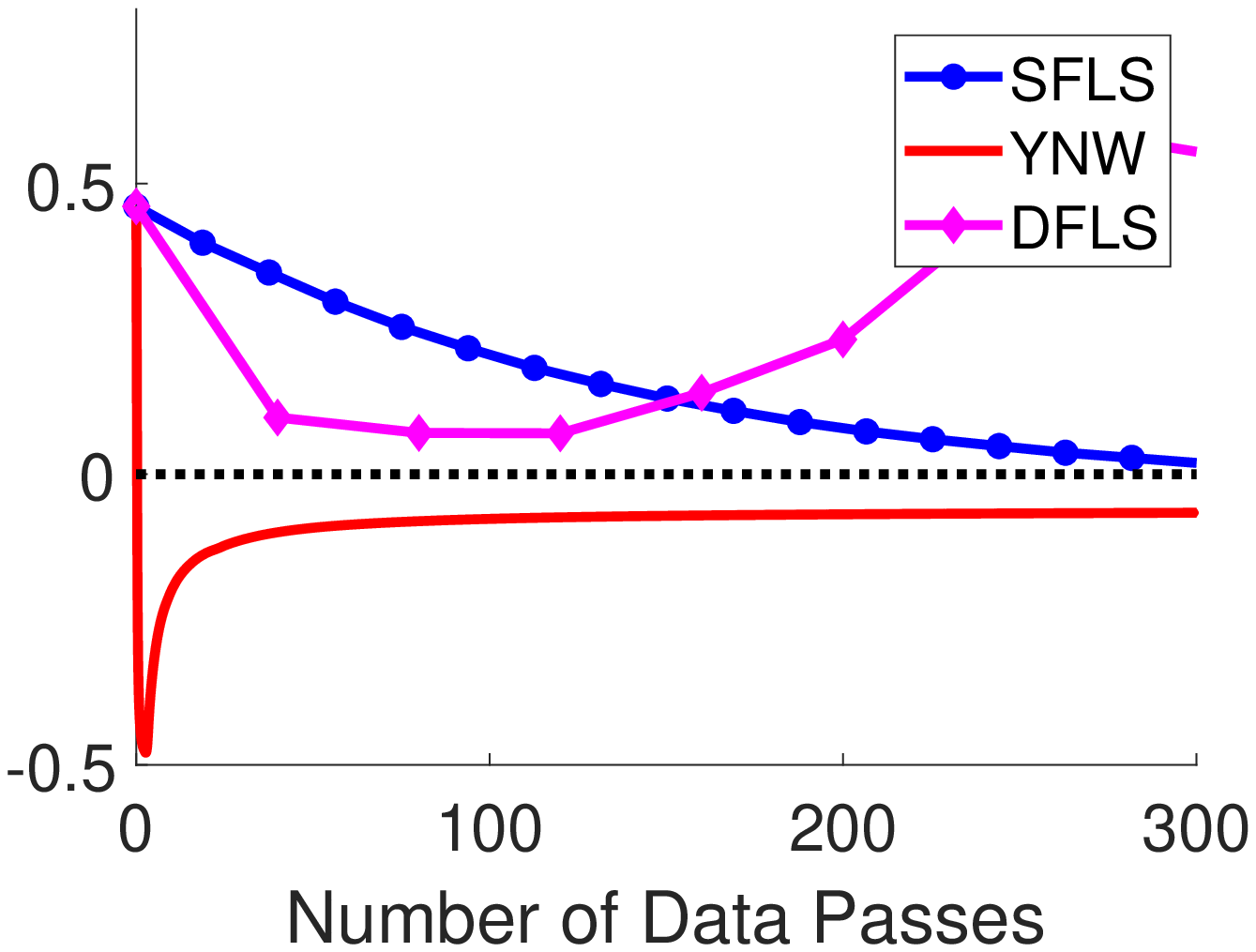}\\
		\raisebox{10ex}{\small{\rotatebox[origin=c]{90}{$\max\limits_{i=1,2,\dots,m}\{f_i(\bx)-r_i\}$}}}
		& \includegraphics[width=0.40\textwidth]{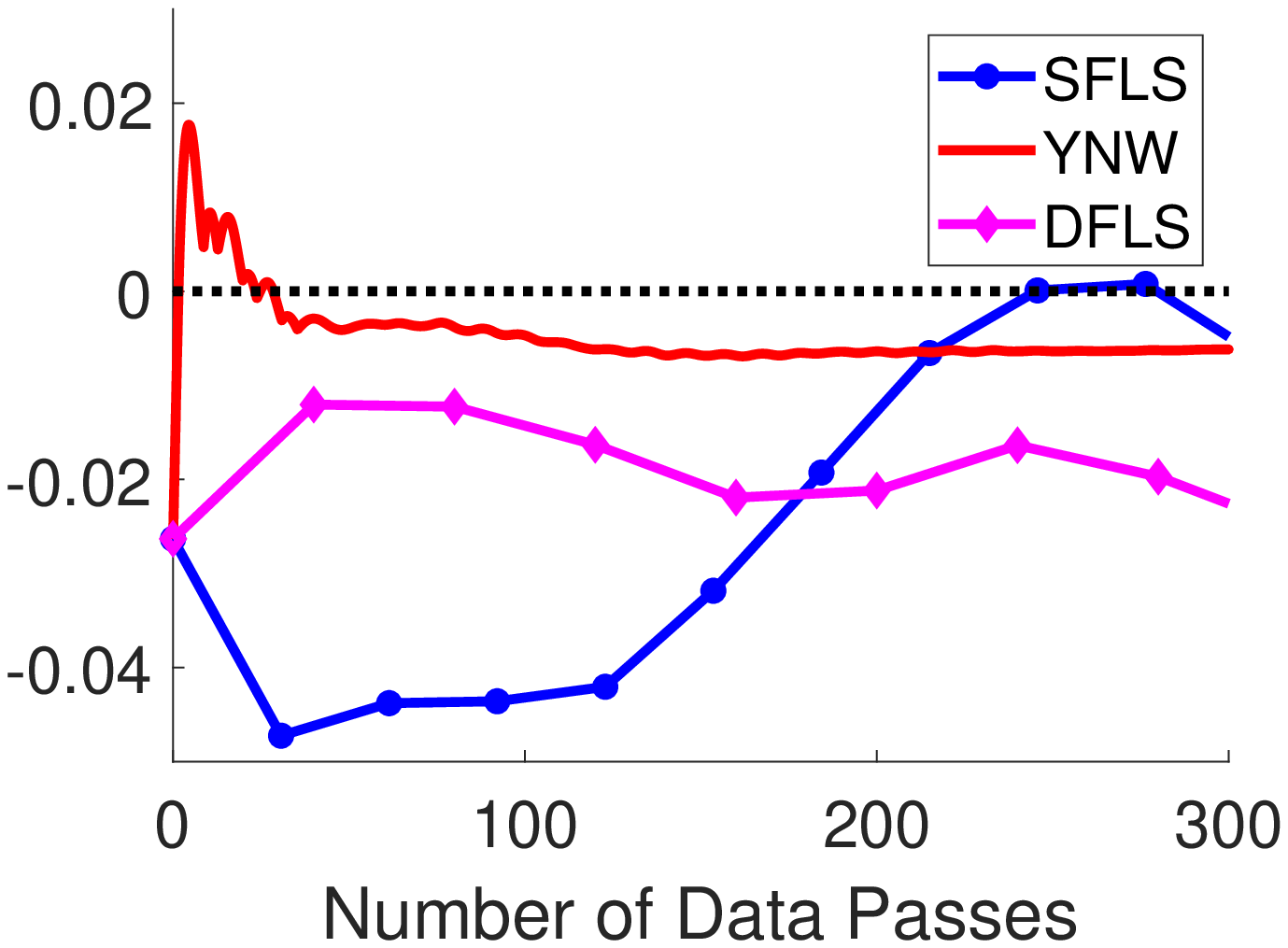}
		& \includegraphics[width=0.40\textwidth]{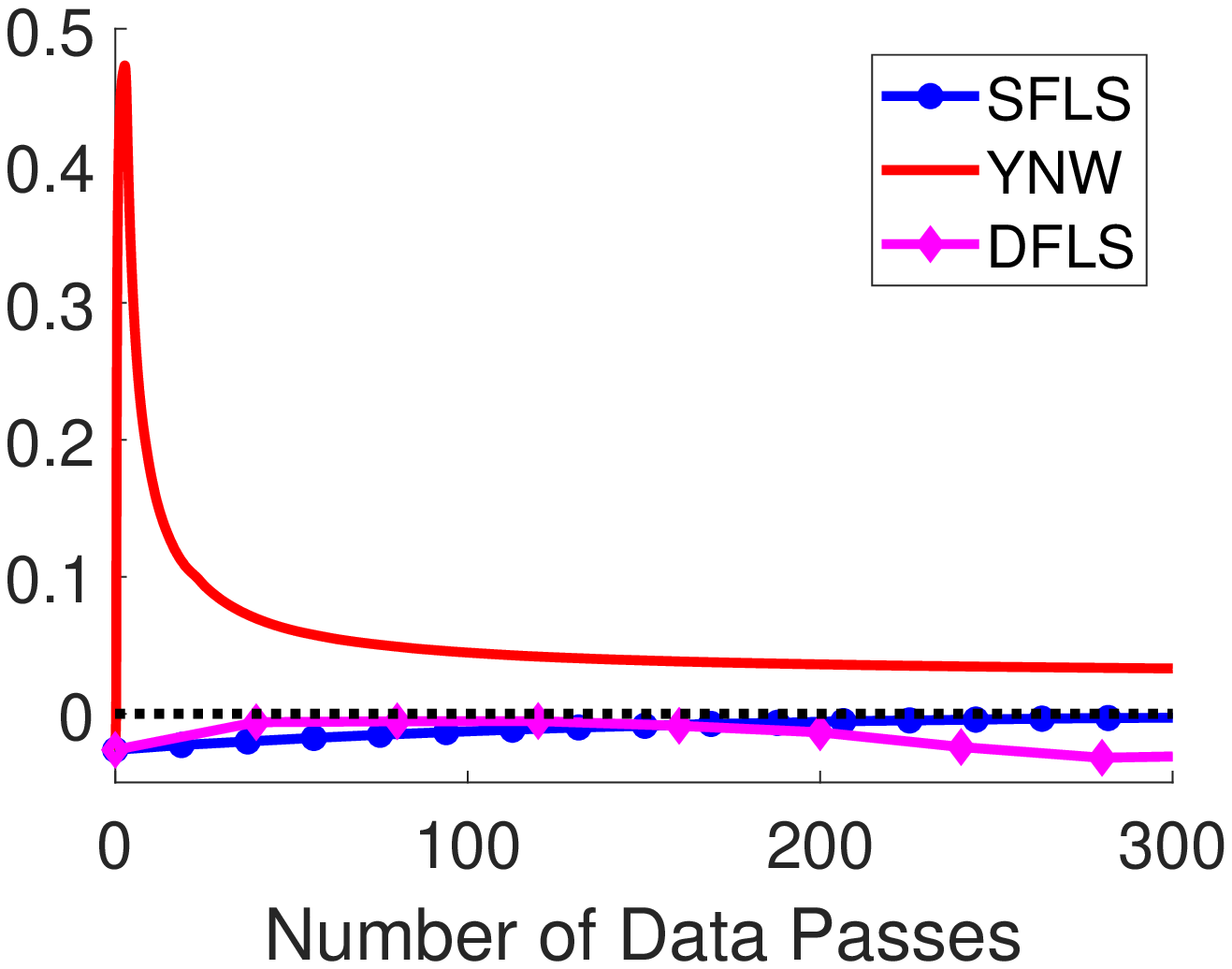}\\		
	\end{tabular}
	\vspace{2ex}
	\caption{Performance of SFLS, YNW, and DFLS for solving the classification problem with fairness constraints.
	}
	\label{fig:fairness}
\end{figure}

For testing, we considered the ``a9a'' dataset, also used by \citet{goh2016satisfying} and another dataset dubbed ``LoanStats'' from \cite{LendingClub2019}. We chose $\lambda=5$, $\kappa=0.95$, and $\phi(z)=(1-z)_+$ in each case. The distributions $\mathcal{D}$, $\mathcal{D}_M$, and $\mathcal{D}_F$ were defined as empirical distributions based on each dataset as described below. The goal in the a9a dataset is to predict people making more than 50,000 USD. Following  \citet{goh2016satisfying}, we used the 32,561 training instances ($\mathcal{D}$) and the 16,281 testing instances in the dataset to construct the objective function and constraints, respectively. Since we need male and female subsets to construct constraints, we further split the testing data into 14,720 male instances ($\mathcal{D}_M$) and 1,561 female instances ($\mathcal{D}_F$). 
The LoanStats dataset contains 
information of $128,375$ loans issued in the fourth quarter of 2018 and the goal is to predict if a loan will be approved or rejected. After creating dummy variables, each loan is represented by a feature vector of 250 dimensions. We randomly partitioned the dataset into a set of $63,890$ loans ($\mathcal{D}$) used to construct the objective function and a set of $64,485$ loans used to build the constraints. We further split the second set based on whether the feature ``homeOwnership'' equals ``Mortgage'' ($\mathcal{D}_M$) or some other value ($\mathcal{D}_F$) to obtain 31,966 and 32,519 loans in two subsets, respectively.

All methods are initialized at $\tilde\bx=\mathbf{0}$, which is feasible for \eqref{eq:NPclassification_multi}. In SFLS and DFLS, we chose $r^{(0)}=1$, $\gamma_t=0.1/\sqrt{t+1}$ and $\theta=1.1$ across all datasets. Note that $r^{(0)}=1=f_0(\mathbf{0})\geq f^*$. In SFLS, we chose $T=300$ and $T=200$ for a9a and LoanStats datasets, respectively. In DFLS, we chose $T=100$ and $T=50$ for a9a and LoanStats datasets, respectively. Both SFLS and YNW employed a mini-batch size of $500$ and $1000$ for a9a and LoanStats datasets, respectively.  Similar to \S\ref{sec:NumClassification}, we chose the number of iterations in YNW so that its total number of data passes is 300. Then, we also terminated SFLS and DFLS when the total data passes they performed exceed $300$.


Figure \ref{fig:fairness} displays the performance of SFLS, YNW, and DFLS as a function of data passes. The interpretation of the axes and line markers in this figure are analogous to the ones in Figure \ref{fig:connect4}. In this figure, $f^*$ is approximately by the objective value returned by DFLS after a sufficient number of data passes (i.e. at least $5000$ data passes with $2T$ inner iterations.) On the a9a dataset, SFLS maintains a feasible solution path, as expected, while the YNW solutions are initially infeasible and become feasible with more data passes. The YNW reduces optimality gap more rapidly at the beginning while SFLS catches up quickly. The objective function value of YNW cannot be interpreted as an optimality gap when its solutions are infeasible since the corresponding objective function value can be super optimal. This feature is clearly visible on the LoanStats data. Here most of the YNW solutions are infeasible and superoptimal, that is, $f(x) - f^*$ is non-positive. The SFLS solution path continues to be feasible and suboptimal on this dataset, with its suboptimality decreasing consistently after each outer iterations. DFLS also produces a feasible path but does not effectively reduce the optimality gap because its data complexity is high, that is, it requires a large number of data passes to achieve a small optimality gap. Similar to \S\ref{sec:NumClassification}, we once again find that the low data complexity of SFLS is critical to balance optimality and feasibility when solving an SOEC. 

\section{Conclusion}\label{sec:conclusion}
We consider constrained optimization models where both the objective function and multiple constraints contain expectations of random convex functions. These models, referred to as stochastic optimization problems with expectation constraints (SOECs), arise in several machine learning, engineering, and business applications. We develop a stochastic feasible level-set method (SFLS) to solve SOECs, propose a tractable oracle to be used with SLFS, and analyze related iteration complexities. SFLS's total iteration complexity is comparable to stochastic subgradient methods in terms of $\epsilon$ but depends on a condition number -- the cost of requiring feasibility. We evaluate the performance of SFLS across three applications involving approximate linear programming, multi-class classification, and learning classifiers with fairness constraints. We find that SFLS exhibits key advantages over existing methods. First, it ensures a feasible solution path with high probability while an existing state-of-the-art stochastic subgradient method can return highly infeasible solutions when terminated before conservative termination criteria are met. Infeasibilities may void the use of a solution in practice, especially if constraints model implementation requirements. Thus, the ability of SFLS to compute feasible solutions before convergence is practically relevant. Second, SFLS computes feasible solutions with small optimality gaps using only a few data passes owing to its low data-complexity, which is a desirable property when expectations are defined using large datasets that are expensive to scan. In contrast to SFLS, a recent deterministic feasible level set method exhibits high data complexity and large optimality gaps. Our theoretical and numerical findings bode well for the use of SFLS to solve SOECs and motivates further research into stochastic first order methods that emphasize feasibility.\looseness=-1


{
\bibliographystyle{informs2014}
	\bibliography{SIP}
}

\ECSwitch
\ECHead{Electronic Companion}
\section{Proofs of Theoretical Results}
In this section, we provide the proofs of all technical results in the paper. 

\proof{\bf Proof of Lemma \ref{lem:feasOracle}: }
Since $r>f^*$, it follows from Lemma~\ref{lem:knownPropsOfL}(c) that $H(r) \leq 0$. Therefore, since $\theta \geq 1$ we have $\epsilon  \leq -\frac{\theta-1}{\theta+1}H(r) \leq - H(r)$. Moreover, by Definition \ref{def:feasOracle}, we have $\Pc(r,\hat \bx)\leq H(r)+\epsilon$ with probability of at least $1-\delta$, which implies that $\hat{x}$ is a feasible solution to \eqref{eq:gcols} since $\Pc(r,\hat \bx)\leq H(r)+\epsilon \leq H(r) - H(r) \leq 0$.  \hfill\halmos\endproof

Proof of Theorem~\ref{generalcomplexity} depends on the following lemma. 
\begin{lemma}\label{lem:absOracle}
Given an input tuple $(r,\epsilon,\delta,\theta)$, a stochastic oracle $\mathcal{A}(r, \epsilon, \delta)$ with $0<\epsilon\leq -\frac{\theta-1}{\theta+1}H(r)$ returns $U(r)$ and $\hat\bx$ such that $\theta U(r)\leq H(r)\leq\Pc(r,\hat \bx)\leq U(r)/\theta$ with probability of at least $1-\delta$.\looseness = -1
\end{lemma}
\proof{Proof.}
The inequality $H(r)\leq \Pc(r,\hat \bx)$ holds by definition of $H(r)$. By definition of stochastic oracle (Definition \ref{def:feasOracle}) and the property of $\epsilon$, it follows that $\Pc(r,\hat \bx)\leq H(r) +\epsilon\leq \frac{2}{\theta+1}H(r)$, $H(r)\leq U(r)+\epsilon\leq U(r)-\frac{\theta-1}{\theta+1}H(r)$, and $ U(r) \leq H(r) +\epsilon\leq \frac{2}{\theta+1}H(r)$ hold with probability of at least $1-\delta$. Since $r>f^*$, Lemma~\ref{lem:knownPropsOfL}(c) implies that $H(r) \leq 0$. Therefore, using the inequality $ U(r) \leq  \frac{2}{\theta+1}H(r)$ we get $U(r)\leq0$ and $ \theta U(r) \leq \frac{\theta+1}{2}U(r) \leq H(r)$ since $\theta>1$. Finally, combining the inequalities $\Pc(r,\hat \bx)\leq \frac{2}{\theta+1}H(r)$ and $H(r)\leq U(r)-\frac{\theta-1}{\theta+1}H(r)$ (or equivalently $H(r)\leq \frac{\theta+1}{2\theta}U(r)$), we get $\Pc(r, \hat x)\leq \frac{2}{\theta+1}\cdot \frac{\theta+1}{2\theta}U(r)=U(r)/\theta$. \looseness = -1 
\hfill\halmos\endproof

In the proof of Theorem~\ref{generalcomplexity} we need the following property of the condition measure $\beta$. In particular, it can be easily verified from the convexity of $H(r)$ and $H(r) -\delta \leq H(r+\delta) \leq H(r)$ for any $\delta\geq0$ (Lemma 2.3.5 in \citealp{nesterov2004introductory}) that 
$\frac{H(r)}{r - f^*}$ is monotonically increasing in $r$ on $(f^*,r^{(0)}]$ and  
\begin{equation}\label{eqn:LstarConvexityRelationship}
-\beta = \frac{H(r^{(0)})}{r^{(0)} - f^*} \geq \frac{H(r)}{r - f^*}\geq -1,\quad \forall r\in(f^*,r^{(0)}]. 
\end{equation}

\proof{\bf Proof of Theorem \ref{generalcomplexity}:}
We first show that the Algorithm~\ref{alg:lsdecreasing} generates a feasible solution at each iteration with high probability. 
Let $K$ be the largest value of $k$ such that $r^{(k)}>f^*$ and  the following inequality holds:
\begin{eqnarray}
\label{eq:phasechange}
\epsoracle= -\frac{\theta-1}{2\theta^2(\theta+1)}  H(r^{(0)}) \epsilon\leq-\frac{\theta-1}{\theta+1}H(r^{(k)}).
\end{eqnarray} 
Notice that $K \geq 0$ since $0< \epsilon\leq 1 \leq2\theta^2$ and $H(r^{(0)}) \leq 0$. It follows from Lemma~\ref{lem:absOracle} that with a probability of at least $1-\delta^{(k)}$ we have,
\begin{eqnarray}
\label{eq:chainineq}
\theta U(r^{(k)})\leq H(r^{(k)})\leq\Pc(r^{(k)},\bx^{(k)})\leq U(r^{(k)})/\theta, \quad \text{for any $k\geq 0$}.
\end{eqnarray} 
Since $r^{(k+1)} = r^{(k)} + U(r^{(k)})/(2\theta)$, we have 
\begin{align}
\label{eq:LS1}
r^{(k+1)} - f^* &= r^{(k)} - f^* + U(r^{(k)})/(2\theta)\geq r^{(k)} - f^* + H(r^{(k)})/2\geq \frac{1}{2}(r^{(k)} - f^*),
\end{align}
and 
\begin{align}\label{eq:LS2}
r^{(k+1)} - f^*&= r^{(k)} - f^* + U(r^{(k)})/(2\theta)
\leq r^{(k)} - f^* + \dfrac{H(r^{(k)})}{2\theta^2}
\leq \left(1- \frac{\beta}{2\theta^2}\right)(r^{(k)} - f^*)
\end{align}
with a probability of at least $1-\delta^{(k)}$, where the last inequalities in both \eqref{eq:LS1} and \eqref{eq:LS2} follow from \eqref{eqn:LstarConvexityRelationship}. Inequality~\eqref{eq:LS1} and the condition $r^{(k)}>f^*$ imply that $r^{(k+1)} > f^*$.
Applying this argument recurrently  and using the fact that $\sum_{k=0}^\infty\delta^{(k)}=\delta$, we have \eqref{eq:LS1}, \eqref{eq:LS2} and $r^{(k+1)} > f^*$ holds for $k=0,1,\ldots, K$. Therefore, since $\epsoracle\leq -\frac{\theta-1}{\theta+1}H(r^{(k)})\leq -H(r^{(k)})$ for $k=0,1,\ldots,K$, Lemma~\ref{lem:feasOracle} implies the solution $\bx^{(k)}$ generated at iteration  $k=0,1,\ldots,K$ is feasible to~\eqref{eq:gco} with a probability of at least $1-\delta$.
We next show that \eqref{eq:phasechange} holds with a high probability until Algorithm~\ref{alg:lsdecreasing} terminates. By the definition of $K$, we know that \eqref{eq:phasechange} is violated when $k=K+1$, i.e. 
$-\frac{\theta-1}{2\theta^2(\theta+1)} H(r^{(0)}) \epsilon>-\frac{\theta-1}{\theta+1}H(r^{(K+1)})$. Since $r^{(k+1)} \leq r^{(k)}$ and $\frac{H(r)}{r - f^*}$ is monotonically increasing, we can show that 
\begin{equation}\label{eq:epsoptandrK}
-\frac{\theta-1}{2\theta^2(\theta+1)} H(r^{(0)}) \epsilon>-\frac{\theta-1}{\theta+1}H(r^{(K+1)})\geq -\frac{\theta-1}{\theta+1}H(r^{(K)})\frac{r^{(K+1)}-f^*}{r^{(K)}-f^*}\geq  -\frac{\theta-1}{2(\theta+1)}H(r^{(K)}),
\end{equation}
where the last inequality holds by \eqref{eq:LS1}. Using the definition of $\epsilon_{\text{opt}}$, \eqref{eq:epsoptandrK}, and \eqref{eq:chainineq} for $k=K$ (specifically, $H(r^{(K)}) \leq U(r^{(K)})/\theta$), we have 
$$
-\frac{\theta-1}{2\theta(\theta+1)}\epsilon_{\text{opt}} = \frac{\theta-1}{2\theta^2(\theta+1)}	H(r^{(0)}) \epsilon\leq \frac{\theta-1}{2(\theta+1)}H(r^{(K)})\leq \frac{\theta-1}{2\theta(\theta+1)}U(r^{(K)})
$$
which indicates that Algorithm~\ref{alg:lsdecreasing} must stop before $k= K+1$. Therefore, SFLS generates a feasible solution with a probability of at least $1-\delta$ at each iteration before termination.

We now proceed to establish that the terminal solution of SFLS is relative $\epsilon$-optimal solution. By definition of $\Pc(r^{(k)}, \bx^{(k)})$ and \eqref{eq:chainineq} it follows that $f_0(\bx^{(k)}) -  r^{(k)}  \leq \Pc(r^{(k)}, \bx^{(k)})\leq H(r^{(k)})/\theta^2 \leq 0$ for all $k$. Hence,
\begin{equation}\label{eqn:epsilonOptTemp}
f_0(\bx^{(k)}) - f^*  \leq r^{(k)} - f^*,\quad\text{for all }k=0,1,2,\ldots,K.
\end{equation}
Combining \eqref{eqn:epsilonOptTemp} and $r^{(k)} - f^* \leq (r^{(0)} - f^*)H(r^{(k)})/H(r^{(0)})$ derived from \eqref{eqn:LstarConvexityRelationship} stipulates that with a probability of at least $1-\delta$: 
\[\frac{f_0(\bx^{(k)}) - f^*}{r^{(0)} - f^*} \leq \frac{H(r^{(k)})}{H(r^{(0)})}\leq \frac{\theta U(r^{(k)})}{H(r^{(0)})},\] 
where we used \eqref{eq:chainineq} in the second inequality. Hence, at termination of Algorithm~\ref{alg:lsdecreasing} we get $\frac{f_0(\bx^{(k)}) - f^*}{r^{(0)} - f^*} \leq\epsilon$ since the algorithm stops when $\theta U(r^{(k)}) \geq H(r^{(0)})\epsilon$.

Finally we show that $K := \dfrac{2\theta^2}{\beta}\ln\left(\dfrac{\theta^2}{\beta\epsilon}\right)$. By recursively applying inequality~\eqref{eq:LS2} we get
\begin{equation}\label{eqn:alg1OuterIterCompTemp}
0\leq r^{(k)} - f^* \leq \left(1 - \frac{\beta}{2\theta^2}\right)^{k} (r^{(0)} - f^*),\quad\text{for all }k
\end{equation}
with probability of at least $1-\delta$, which implies $r^{(K)} - f^*\leq -\frac{H(r^{(0)})\epsilon}{\theta^2 }$ for the choice of $K$. Hence, we have 
$-U(r^{(K)}) \leq -\theta H(r^{(K)}) \leq \theta(r^{(K)} - f^*) \leq -\epsilon H(r^{(0)})/\theta$ where the first inequality follows by \eqref{eq:chainineq}, the second by~\eqref{eqn:LstarConvexityRelationship}, and the third by~\eqref{eqn:alg1OuterIterCompTemp}. This indicates that the stopping criterion of Algorithm~\ref{alg:lsdecreasing} holds with a probability of at least $1-\delta$ when $k=K$ and SFLS requires at most $K$ calls to oracle $\mathcal{A}$.\looseness = -1  
\hfill\halmos\endproof

\proof{\bf Proof of Proposition~\ref{validIdealOracle}:} The proof of the first part directly follows from Proposition 3.2 in \citealp{nemirov09}. We only show that SMD is a valid oracle. It is straightforward to see that the inequality $U(\bar x^{(t)}) - L(\bar y^{(t)}) \leq \epsoracle$ implies $\Pc(r, \bar x^{(t)})- H(r) \leq U(\bar x^{(t)}) - H(r) \leq U(\bar x^{(t)}) - L(\bar y^{(t)}) \leq \epsoracle$, where the first inequality holds since $U(\bar x^{(t)})$ is an upper bound on $\Pc(r, \bar x^{(t)})$ and the second since $L(\bar y^{(t)})$ is a lower bound on $H(r)$. This indicates that the conditions provided in Definition~\ref{def:feasOracle} are satisfied.  \looseness = -1
\hfill\halmos\endproof

To show part (i) of Proposition \ref{thm:luconvergecor}, we use known lemmas \ref{lem:sumterm1} and \ref{lem:sumterm} as well as prove lemmas  \ref{lem:simpResultsForDetAnal} and \ref{lem:luconverge}. To prove part (ii) of this proposition we need Lemma~\ref{lem:hatluconverge}. Before stating these lemmas, we present some required notation and representations, which we present next. We denote the diameter of $\mathcal{Z}$ with respect to $\omega_z$ by
$$
D_z:=\sqrt{\max_{\bz\in\mathcal{Z}}\omega_z(\bz)-\min_{\bz\in\mathcal{Z}}\omega_z(\bz)}=1.
$$
In addition,
for any $\bzeta_x\in\mathbb{R}^{d}$,  $\bzeta_y\in\mathbb{R}^{m+1}$,  $\bx'\in\mathcal{X}^o$, $\by'\in\mathcal{Y}^o$, and $\bz' = (\bx',\by')\in \mathcal{Z}^o$, it is easy to verify for $\bzeta=(\bzeta_x,\bzeta_y)$ that 
\begin{eqnarray}
\label{eq:relationship}
P_{\bz'}(\bzeta)=\left(P_{\bx'}^x(2D_x^2\bzeta_x),P_{\by'}^y(2D_y^2\bzeta_y)\right),
\end{eqnarray} 
where $P_{\bx'}^x(\bzeta_x):=\argmin_{\bx\in \mathcal{X}}\{\bzeta_x^\top(\bx-\bx')+V_x(\bx',\bx)\}$ and $P_{\by'}^y(\bzeta_y):=\argmin_{\by\in \mathcal{Y}}\{\bzeta_y^\top(\by-\by')+V_y(\by',\by)\}$. 

\begin{lemma}[Equation (2.37) and Lemma 6.1 in \citealp{nemirov09}]
	\label{lem:sumterm1}
	\begin{enumerate}
		\item	Let $\bzeta_x^{(t)}\in\mathbb{R}^{d}$, $t = 0,1,2,\ldots$ be a set of random variables, $\bv^{(0)}\in\mathcal{X}^o$ and $\bv^{(t+1)}=P^x_{\bv^{(t)}}(\bzeta_x^{(t)})$ for $t=0,1,2,\ldots$. For any $\bv\in\mathcal{X}$ and $t\geq 1$, we have
		\begin{eqnarray*}
			\sum_{s=0}^{t}(\bv^{(s)}-\bv)^\top\bzeta_x^{(s)}\leq V_x(\bv^{(0)},\bv)+\frac{1}{2\alpha_x}\sum_{s=0}^{t}\left\|\bzeta_x^{(s)}\right\|_{*,x}^2.
		\end{eqnarray*}
		\item Let $\bzeta_y^{(t)}\in\mathbb{R}^{m+1}$, $t = 0,1,2,\ldots$ be a set of random variables, $\bv^{(0)}\in\mathcal{Y}^o$ and $\bv^{(t+1)}=P^y_{\bv^{(t)}}(\bzeta_y^{(t)})$ for $t=0,1,2,\ldots$. For any $\bv\in\mathcal{Y}$ and $t\geq 1$, we have
		\begin{eqnarray*}
			\sum_{s=0}^{t}(\bv^{(s)}-\bv)^\top\bzeta_y^{(s)}\leq V_y(\bv^{(0)},\bv)+\frac{1}{2\alpha_y}\sum_{s=0}^{t}\left\|\bzeta_y^{(s)}\right\|_{*,y}^2.
		\end{eqnarray*}
		\item Let $\bzeta^{(t)}\in\mathbb{R}^{d+m+1}$, $t = 0,1,2,\ldots$ be a set of random variables, $\bv^{(0)}\in\mathcal{Z}^o$ and $\bv^{(t+1)}=P_{\bv^{(t)}}(\bzeta^{(t)})$ for $t=0,1,2,\ldots$. For any $\bv\in\mathcal{Z}$ and $t\geq 1$, we have
		\begin{eqnarray*}
			\sum_{s=0}^{t}(\bv^{(s)}-\bv)^\top\bzeta^{(s)}\leq V(\bv^{(0)},\bv)+\frac{1}{2}\sum_{s=0}^{t}\left\|\bzeta^{(s)}\right\|_{*,z}^2.
		\end{eqnarray*}
	\end{enumerate}
\end{lemma}

\begin{lemma}[Lemma 2 in \citealp{lan2012validation}]
	\label{lem:sumterm}
	Let $\bxi^{(t)}$ and $\sigma_t>0$ for $t=0,1,2,\ldots$ be respectively a sequence of i.i.d. random variables and deterministic numbers; $\bxi^{[t]}=(\bxi^{(0)},\bxi^{(1)},\dots,\bxi^{(t)})$; $\E_t$ the conditional expectation conditioning on $\bxi^{[t-1]}$ for $t\geq1$; and $\psi_t(\bxi^{[t]})$ be a measurable function of $\bxi^{[t]}$ such that either 
	\begin{itemize}
		\item[] {\bf Case A:} $\E_t\left[\psi_t\left(\bxi^{[t]}\right)\right]=0$ and $\E_t\left[\exp\left(\psi_t\left(\bxi^{[t]}\right)^2/\sigma_t^2\right)\right]\leq \exp(1)$, or
		\item[] {\bf Case B:} $\E_t\left[\exp\left(\left\vert\psi_t\left(\bxi^{[t]}\right)\right\vert/\sigma_t\right)\right]\leq \exp(1)$ ,
	\end{itemize}
	almost surely for all $t$. Then for any $\Omega>0$, we have the followings:\\
	In case A:
	$$
	\text{Prob}\left\{\sum_{s=0}^{t}\psi_s>\Omega\sqrt{\sum_{s=0}^{t}\sigma_s^2}\right\}\leq \exp(-\Omega^2/3).
	$$
	In case B: 
	\begin{eqnarray*}
		&&\text{Prob}\left\{\sum_{s=0}^{t}\psi_s>\|\sigma^{[t]}\|_1+\Omega\|\sigma^{[t]}\|_2\right\}
		\leq \exp(-\Omega^2/12)+\exp(-3\Omega/4),
	\end{eqnarray*}
	where $\sigma^{[t]}=(\sigma_0,\sigma_1,\dots,\sigma_{t})^\top$.
\end{lemma}

Lemma \ref{lem:simpResultsForDetAnal} shows that the stochastic subgradient $G(\cdot, \cdot, \cdot)$ has a light-tailed distribution and bounds the Bregmann distances. Define 
\begin{eqnarray*}
	\Delta_t&:=&G(\bx^{(t)},\by^{(t)},\bxi^{(t)})-g(\bx^{(t)},\by^{(t)})=
	\left[
	\begin{array}{ll}
		\Delta_t^x\\
		-\Delta_t^y
	\end{array}
	\right]:=
	\left[
	\begin{array}{ll}
		G_x(\bx^{(t)},\by^{(t)},\bxi^{(t)})-g_x(\bx^{(t)},\by^{(t)})\\
		g_y(\bx^{(t)},\by^{(t)})-G_y(\bx^{(t)},\by^{(t)},\bxi^{(t)})
	\end{array}
	\right].
\end{eqnarray*} 
\begin{lemma}\label{lem:simpResultsForDetAnal}The following inequalities hold: 
	\begin{align}
	\E_t\left[\exp\left(\|\Delta_t\|_{*,z}^2/(2M)^2\right)\right] \leq \exp(1),\label{eq:lighttailgraddeltaz}\\
	\E_t\left[\exp\left(\|\Delta_t^x\|_{*,x}^2/(2M_x)^2\right)\right]\leq\exp(1),\label{eq:lighttailgraddeltax}\\
	\E_t\left[\exp\left(\|\Delta_t^y\|_{*,y}^2/(2M_y)^2\right)\right]\leq\exp(1).\label{eq:lighttailgraddeltay}
	\end{align}
	Moreover, when $\bz'=(\bx',\by'):=\argmin_{\bz\in\mathcal{Z}}\omega_z(\bz)$, we have
	\begin{alignat}{3}
	\label{eq:Vxbound}
	\frac{\alpha_x}{2}\|\bx'-\bx\|_x^2&\leq V_x(\bx',\bx)&&\leq D_x^2,&&\quad\text{for all }\bx\in\mathcal{X},\\
	\label{eq:Vybound}
	\frac{\alpha_y}{2}\|\by'-\by\|_y^2&\leq V_y(\by',\by)&&\leq D_y^2,&&\quad\text{for all }\bx\in\mathcal{Y},\\
	\label{eq:Vzbound}
	\frac{1}{2}\|\bz'-\bz\|_z^2&\leq V(\bz',\bz)&&\leq D_z^2=1,&&\quad\text{for all }\bz\in\mathcal{Z}.
	\end{alignat}
\end{lemma}

\proof{Proof.} 
Applying Jensen's inequality and using the definitions of $\|\cdot\|_{*,z}$, $M$, and the inequalities \eqref{eq:lighttailgradx} and \eqref{eq:lighttailgrady}, we have 
\begin{eqnarray}
\nonumber
&&\E\left[\exp\left(\|G(\bx,\by,\bxi)\|_{*,z}^2/M^2\right)\right]\\\nonumber
&=&\E\left[\exp\left(\frac{\frac{2D_x^2}{\alpha_x}\|G_x(\bx,\by,\bxi)\|_{*,x}^2+\frac{2D_y^2}{\alpha_y}\|G_y(\bx,\by,\bxi)\|_{*,y}^2}{\frac{2D_x^2}{\alpha_x}M_x^2+\frac{2D_y^2}{\alpha_y}M_y^2}\right)\right]\\\nonumber
&\leq&\frac{\frac{2D_x^2}{\alpha_x}M_x^2\E\left[\exp\left(\|G_x(\bx,\by,\bxi)\|_{*,x}^2/M_x^2\right)\right]+\frac{2D_y^2}{\alpha_y}M_y^2\E\left[\exp\left(\|G_y(\bx,\by,\bxi)\|_{*,y}^2/M_y^2\right)\right]}{\frac{2D_x^2}{\alpha_x}M_x^2+\frac{2D_y^2}{\alpha_y}M_y^2}\\
&\leq&\exp(1).
\label{eq:lighttailgrad}
\end{eqnarray} 
Using \eqref{eq:lighttailgrad} and Jensen's inequality, it follows that
\begin{equation}\label{eq: gupperbound}
\left\|g(\bx^{(t)},\by^{(t)})\right\|_{*,z}^2\leq\E_t\left[\left \|G(\bx^{(t)},\by^{(t)},\bxi^{(t)})\right \|_{*,z}^2\right]\leq M^2.
\end{equation}
Hence, we have 
\begin{align}
\nonumber
&\E_t\left[\exp\left(\|\Delta_t\|_{*,z}^2/(2M)^2\right)\right]\\
\nonumber &\leq \E_t\left[\exp(2\|G(\bx^{(t)},\by^{(t)},\bxi^{(t)})\|_{*,z}^2/(2M)^2)\exp(2\|g(\bx^{(t)},\by^{(t)})\|_{*,z}^2/(2M)^2)\right],\\\nonumber
&\leq \E_t\left[\sqrt{\exp(\|G(\bx^{(t)},\by^{(t)},\bxi^{(t)})\|_{*,z}^2/M^2)}\exp(1/2)\right],\\\nonumber
&\leq \sqrt{\E_t\left[\exp(\|G(\bx^{(t)},\by^{(t)},\bxi^{(t)})\|_{*,z}^2/M^2)\right]}\exp(1/2),\\\label{eq:lighttailgraddelta}
&\leq \exp(1/2)\exp(1/2)=\exp(1),
\end{align}
where the first inequality follows from the definition of $\Delta_t$ and the inequality $\|a+b\|^2 \leq 2a^2 + 2b^2$ for any $a,b\in\mathbb{R}$, the second from \eqref{eq: gupperbound}, the third from Jensen's inequality for concave functions, and the fourth by inequalities \eqref{eq:lighttailgradx} and \eqref{eq:lighttailgrady}.
Following a similar argument, we can also show that $\E_t\left[\exp\left(\|\Delta_t^x\|_{*,x}^2/(2M_x)^2\right)\right]\leq\exp(1)$ and $\E_t\left[\exp\left(\|\Delta_t^y\|_{*,y}^2/(2M_y)^2\right)\right]\leq\exp(1)$.
Finally, inequalities \eqref{eq:Vxbound}, \eqref{eq:Vybound}, and \eqref{eq:Vzbound} follow because $\omega_x$, $\omega_y$ and $\omega_z$ are modulus $\alpha_x$, $\alpha_y$ and $1$, respectively.  
\hfill\halmos\endproof

\begin{lemma} 
	\label{lem:luconverge}	
	Let $\nu_{s,t}:=\dfrac{\gamma_s}{\sum_{s'=0}^{t}\gamma_{s'}}$. Given $\Omega>0$, Algorithm~\ref{alg:SMDwithOnlineValidation} computes $(\bx^{(t)},\by^{(t)})$, $t = 1,2,3,\ldots,$ such that 
	\begin{align}
	\nonumber
	&\text{Prob}\left\{u^{(t)}_*-l^{(t)}_*>4\sqrt{2}\Omega M\sqrt{\sum_{s=0}^t\nu_{s,t}^2}+\frac{2+2.5M^2\sum_{s=0}^t\gamma_s^2}{\sum_{s=0}^t\gamma_s} +2.5\Omega M^2\sqrt{\sum_{s=0}^t\gamma_s^2\nu_{s,t}^2}\right\}\\\label{eq:boundabsgaplu}
	&\leq  \exp(-\Omega^2/3)+\exp(-\Omega^2/12)+\exp(-3\Omega/4).
	\end{align}
	\normalsize
\end{lemma}
\proof{Proof.}
Since  $\bz^{(0)}\in\argmin_{\bz\in\mathcal{Z}} \omega_z(\bz)$ and $\bz^{(t+1)}=P_{\bz^{(t)}}(\gamma_tG(\bx^{(t)},\by^{(t)},\bxi^{(t)}))$ in Algorithm~\ref{alg:SMDwithOnlineValidation}, by Lemma~\ref{lem:sumterm1} we have, for any $\bz\in\mathcal{Z}$,
\begin{align}
\label{eq:lemma1withx}
\sum_{s=0}^{t}\gamma_s(\bz^{(s)}-\bz)^\top G(\bx^{(s)},\by^{(s)},\bxi^{(s)})&\leq V(\bz^{(0)},\bz)+\frac{1}{2}\sum_{s=0}^{t}\gamma_s^2\left\|G(\bx^{(s)},\by^{(s)},\bxi^{(s)})\right\|_{*,z}^2\nonumber\\
&\leq 1+\frac{1}{2}\sum_{s=0}^{t}\gamma_s^2\left\|G(\bx^{(s)},\by^{(s)},\bxi^{(s)})\right\|_{*,z}^2,
\end{align}
where the second inequality follows by \eqref{eq:Vzbound}. In addition, by definition of $\Delta_t$, for any $\bz\in\mathcal{Z}$ we have
\small
\begin{align}
\nonumber
&\frac{1}{\sum_{s=0}^{t}\gamma_s}\sum_{s=0}^{t}\gamma_s(\bz^{(s)}-\bz)^\top G(\bx^{(s)},\by^{(s)},\bxi^{(s)})\\
&=\frac{\sum_{s=0}^{t}\gamma_s(\bx^{(s)}-\bx)^\top g_x(\bx^{(s)},\by^{(s)})}{\sum_{s=0}^{t}\gamma_s}-
\frac{\sum_{s=0}^{t}\gamma_s(\by^{(s)}-\by)^\top g_y(\bx^{(s)},\by^{(s)})}{\sum_{s=0}^{t}\gamma_s}+
\frac{\sum_{s=0}^{t}\gamma_s(\bz^{(s)}-\bz)^\top \Delta_s}{\sum_{s=0}^{t}\gamma_s}\label{eq:getul}
\end{align}
\normalsize
Applying \eqref{eq:lemma1withx} to \eqref{eq:getul} and reorganizing terms lead to 
\begin{align*}
&\frac{\sum_{s=0}^{t}\gamma_s(\bx^{(s)}-\bx)^\top g_x(\bx^{(s)},\by^{(s)})}{\sum_{s=0}^{t}\gamma_s}-
\frac{\sum_{s=0}^{t}\gamma_s(\by^{(s)}-\by)^\top g_y(\bx^{(s)},\by^{(s)})}{\sum_{s=0}^{t}\gamma_s}\\
&\leq \frac{\sum_{s=0}^{t}\gamma_s(\bz-\bz^{(s)})^\top \Delta_s}{\sum_{s=0}^{t}\gamma_s}
+\frac{1+0.5\sum_{s=0}^{t}\gamma_s^2\left\|G(\bx^{(s)},\by^{(s)},\bxi^{(s)})\right\|_{*,z}^2}{\sum_{s=0}^{t}\gamma_s}.
\end{align*}
Maximizing both sides of the above inequality over $\bz\in\mathcal{Z}$ implies 
\begin{align}
\label{eq:boundulstar1}
u^{(t)}_*-l^{(t)}_*&\leq
\frac{\max\limits_{\bz\in\mathcal{Z}}\left[\sum_{s=0}^{t}\gamma_s(\bz-\bz^{(s)})^\top \Delta_s\right]}{\sum_{s=0}^{t}\gamma_s}
+\frac{1+0.5\sum_{s=0}^{t}\gamma_s^2\left\|G(\bx^{(s)},\by^{(s)},\bxi^{(s)})\right\|_{*,z}^2}{\sum_{s=0}^{t}\gamma_s}.
\end{align}
Let 
$\bv^{(0)}=\bz^{(0)}$ and $\bv^{(t+1)}=P_{\bv^{(t)}}(-\gamma_t\Delta_t)$  for $t=0,1,2,\ldots$. From Lemma~\ref{lem:sumterm1} it follows that for any $\bz\in\mathcal{Z}$, 
\begin{eqnarray}
\label{eq:lemma1withu}
-\sum_{s=0}^{t}\gamma_s(\bv^{(s)}-\bz)^\top \Delta_s\leq1+0.5\sum_{s=0}^{t}\gamma_s^2\|\Delta_s\|_{*,z}^2,
\end{eqnarray}
Rewriting $\bz - \bz^{(s)} = \bv^{(s)} - \bz^{(s)} + \bz - \bv^{(s)}$ and applying \eqref{eq:lemma1withu} to \eqref{eq:boundulstar1} yield
\small
\begin{align}
u^{(t)}_*-l^{(t)}_*\leq&
\frac{\sum_{s=0}^{t}\gamma_s(\bv^{(s)}-\bz^{(s)})^\top \Delta_s}{\sum_{s=0}^{t}\gamma_s}
+\frac{2+0.5\sum_{s=0}^{t}\gamma_s^2\left(\left\|G(\bx^{(s)},\by^{(s)},\bxi^{(s)})\right\|_{*,z}^2+\left\|\Delta_s\right\|_{*,z}^2\right)}{\sum_{s=0}^{t}\gamma_s}.\label{eq:boundulstar2}
\end{align}
\normalsize
We next find a probabilistic bound for the right hand side of the above inequality. 

\emph{Bound on $\frac{\sum_{s=0}^{t}\gamma_s(\bv^{(s)}-\bz^{(s)})^\top \Delta_s}{\sum_{s=0}^{t}\gamma_s}$:} By our choice of $\bz^{(0)}$, i.e. $\bz^{(0)}  = \argmin_{\bz\in\mathcal{Z}} \omega_z(\bz)$ and \eqref{eq:Vzbound}, for any $s = 0,1,\ldots,t$ we have 
\begin{equation}
\label{eq:bounduz}
\|\bv^{(s)}-\bz^{(s)}\|_z
\leq \|\bz^{(s)}-\bz^{(0)}\|_z+\|\bv^{(s)}-\bz^{(0)}\|_z\leq \sqrt{2V(\bz^{(0)},\bz^{(s)})}+\sqrt{2V(\bz^{(0)},\bv^{(s)})}\leq 2\sqrt{2}.
\end{equation} 
Define
$\psi_s :=\nu_{s,t}(\bv^{(s)}-\bz^{(s)})^\top \Delta_s$ and  $\sigma_s :=4\sqrt{2}M\nu_{s,t}$. Because $\bxi^{(s)}$ is independent of $\bv^{(s)}$ and $\bz^{(s)}$, we have
$\E_s[\psi_s]=0$. In addition, it can be verified that
$
\psi_s^2\leq\nu_{s,t}^2\|\bv^{(s)}-\bz^{(s)}\|_z^2\|\Delta_s\|_{*,z}^2\leq8\nu_{s,t}^2\|\Delta_s\|_{*,z}^2,
$
where the second inequality holds by \eqref{eq:bounduz}. Using this inequality and \eqref{eq:lighttailgraddelta}, we get
$\E_s\left[\exp\left(\psi_s^2/\sigma_s^2\right)\right]\leq \E_s\left[\exp\left(\|\Delta_s\|_{*,z}^2/(2M)^2\right)\right]\leq \exp(1)$.
Hence, it follows from Case A in Lemma~\ref{lem:sumterm} that
\begin{eqnarray}
\label{eq:boundcrosstermx1}
\text{Prob}\left\{\frac{\sum_{s=0}^{t}\gamma_s(\bv^{(s)}-\bz^{(s)})^\top \Delta_s}{\sum_{s=0}^{t}\gamma_s}>4\sqrt{2}\Omega M\sqrt{\sum_{s=0}^t\nu_{s,t}^2} \right\}\leq \exp(-\Omega^2/3).
\end{eqnarray}

\emph{Bound on $\frac{\sum_{s=0}^{t}\gamma_s^2\left(\left\|G(\bx^{(s)},\by^{(s)},\bxi^{(s)})\right\|_{*,z}^2+\left\|\Delta_s\right\|_{*,z}^2\right)}{\sum_{s=0}^{t}\gamma_s}$:} Define $\psi_s :=\gamma_s\nu_{s,t}\big(\left\|G(\bx^{(s)},\by^{(s)},\bxi^{(s)})\right\|_{*,z}^2+\left\|\Delta_s\right\|_{*,z}^2\big)$ and  $\sigma_s :=5M^2\gamma_s\nu_{s,t}$. We then have
\begin{align*}
\E\left[\exp\left(|\psi_s|/\sigma_s\right)\right]&=\E\left[\exp\left(\frac{\|G(\bx^{(s)},\by^{(s)},\bxi^{(s)})\|_{*,z}^2+\|\Delta_s\|_{*,z}^2}{5M^2}\right)\right]\\
&=\E\left[\exp\left(\frac{\left\|G(\bx^{(s)},\by^{(s)},\bxi^{(s)})\right\|_{*,z}^2/M^2+4\|\Delta_s\|_{*,z}^2/(4M^2)}{5}\right)\right]\\
&\leq\frac{1}{5}\E\left[\exp\left(\frac{\left\|G(\bx^{(s)},\by^{(s)},\bxi^{(s)})\right\|_{*,z}^2}{M^2}\right)\right]
+\frac{4}{5}\E\left[\exp\left(\frac{\|\Delta_s\|_{*,z}^2}{4M^2}\right)\right]\leq \exp(1),
\end{align*}
where the first inequality is from Jensen's inequality and the second inequality is from \eqref{eq:lighttailgrad} and \eqref{eq:lighttailgraddelta}.
Hence, from Case B in Lemma~\ref{lem:sumterm} it follows that
\begin{align}
\nonumber
&\text{Prob}\left\{\frac{\sum_{s=0}^{t}\gamma_s^2\left(\|G(\bx^{(s)},\by^{(s)},\bxi^{(s)})\|_{*,z}^2+\|\Delta_s\|_{*,z}^2\right)}{\sum_{s=0}^{t}\gamma_s}>
5M^2\sum_{s=0}^t\gamma_s\nu_{s,t}
+5\Omega M^2\sqrt{\sum_{s=0}^t\gamma_s^2\nu_{s,t}^2}\right\}\\\label{eq:boundcrosstermx2}
&\leq \exp(-\Omega^2/12)+\exp(-3\Omega/4).
\end{align}
\normalsize

The conclusion is hence obtained by upper bounding the right hand size of \eqref{eq:boundulstar2} using the union bound of \eqref{eq:boundcrosstermx1} and \eqref{eq:boundcrosstermx2}. 
\hfill\halmos\endproof

\begin{lemma} 
	\label{lem:hatluconverge}	
	Let $\nu_{s,t}:=\dfrac{\gamma_s}{\sum_{s'=0}^{t}\gamma_{s'}}$. Given $\Omega>0$, Algorithm~\ref{alg:SMDwithOnlineValidation} guarantees that
	
	\begin{align}
	\nonumber
	&\text{Prob}\left\{\left|\hat l^{(t)}_*-l^{(t)}_*\right|>\left(\Omega(\delta) Q+\frac{4\sqrt{2}\Omega(\delta) D_xM_x}{\sqrt{\alpha_x}}\right)\sqrt{\sum_{s=0}^t\nu_{s,t}^2}+\frac{0.5+\frac{4D_x^2M_x^2}{\alpha_x}\sum_{s=0}^t\gamma_s^2}{\sum_{s=0}^{t}\gamma_s}
	\right.\\
	&\quad\quad\quad\quad\quad\quad\left.+\frac{4\Omega(\delta) D_x^2 M_x^2}{\alpha_x}\sqrt{\sum_{s=0}^t\gamma_s^2\nu_{s,t}^2}
	\right\}\leq  6\exp(-\Omega(\delta)^2/3)+\exp(-\Omega(\delta)^2/12)+\exp\left(-3\Omega(\delta)/4\right).\label{eq:boundabsgapllhat}
	\end{align}
	and
	\begin{align}
	\nonumber
	&\text{Prob}\left\{\left|\hat u^{(t)}_*-u^{(t)}_*\right|>\left(\Omega(\delta) Q+\frac{4\sqrt{2}\Omega(\delta) D_yM_y}{\sqrt{\alpha_x}}\right)\sqrt{\sum_{s=0}^t\nu_{s,t}^2}+\frac{0.5+\frac{4D_y^2M_y^2}{\alpha_x}\sum_{s=0}^t\gamma_s^2}{\sum_{s=0}^{t}\gamma_s}
	\right.\\
	&\quad\quad\quad\quad\quad\quad\left.+\frac{4\Omega(\delta) D_y^2 M_y^2}{\alpha_x}\sqrt{\sum_{s=0}^t\gamma_s^2\nu_{s,t}^2}
	\right\}	\leq  6\exp(-\Omega(\delta)^2/3)+\exp(-\Omega(\delta)^2/12)+\exp\left(-3\Omega(\delta)/4\right).\label{eq:boundabsgapuuhat}
	\end{align}
	
\end{lemma}
\proof{Proof.}
Since the proofs of \eqref{eq:boundabsgapllhat} and \eqref{eq:boundabsgapuuhat} are very similar, we will only prove \eqref{eq:boundabsgapllhat}. Let 
$$l^{(t)}(\bx) : = \dfrac{1}{\sum_{s=0}^{t}\gamma_s}\sum_{s=0}^{t}\gamma_s\left[\phi(\bx^{(s)},\by^{(s)})+g_x(\bx^{(s)},\by^{(s)})^\top(\bx-\bx^{(s)})\right],$$ 
and 
$$\hat l^t(\bx):=\frac{1}{\sum_{s=0}^{t}\gamma_s}\sum_{s=0}^{t}\gamma_s\left[\Phi(\bx^{(s)},\by^{(s)},\bxi^{(s)})+G_x(\bx^{(s)},\by^{(s)},\bxi^{(s)})^\top(\bx-\bx^{(s)})\right].$$ 
Define $\delta_t := \Phi(\bx^{(t)},\by^{(t)},\bxi^{(t)})-\phi(\bx^{(t)},\by^{(t)})$. Using this definition and those of $l^{(t)}_* = \min_{\bx\in\mathcal{X}} l^{(t)}(\bx)$, $\hat l^{(t)}_* = \min_{\bx\in\mathcal{X}}\hat l^{(t)}(\bx)$, and $\Delta_t$ we have	
\begin{align}
\nonumber
\left|\hat l^{(t)}_*-l^{(t)}_*\right|
&= \left\vert\min_{\bx\in\mathcal{X}}\hat l^{(t)}(\bx)-\min_{\bx\in\mathcal{X}}l^{(t)}(\bx)\right\vert\\\nonumber
&\leq\max_{\bx\in\mathcal{X}}\left\vert\hat l^{(t)}(\bx)-\l^{(t)}(\bx)\right\vert\\
&\leq\max_{\bx\in\mathcal{X}}\left\vert\frac{\sum_{s=0}^{t}\gamma_s(\bx-\bx^{(s)})^\top\Delta_s^x}{\sum_{s=0}^{t}\gamma_s}\right\vert
+\left\vert\frac{\sum_{s=0}^{t}\gamma_s\delta_s}{\sum_{s=0}^{t}\gamma_s}\right\vert.
\label{eq:gapllhat}
\end{align}	
By \eqref{eq:relationship} and line 5 of Algorithm~\ref{alg:SMD}, we have $\bx^{(t+1)}=P^x_{\bx^{(t)}}(2D_x^2\gamma_tG_x(\bx^{(t)},\by^{(t)},\bxi^{(t)}))$. Let 
$\bw^{(0)}=\bv^{(0)}=\bx^{(0)}$, $\bw^{(t+1)}:=P^x_{\bw^{(t)}}(-2D_x^2\gamma_t\Delta_t^x)$ and $\bv^{(t+1)}:=P^x_{\bv^{(t)}}(2D_x^2\gamma_t\Delta_t^x)$ for $t=0,1,2,\ldots$. From Lemma~\ref{lem:sumterm1} and \eqref{eq:Vxbound} it follows that
\begin{eqnarray*}
	\label{eq:lemma1withvforx}
	-\sum_{s=0}^{t}\gamma_s(\bw^{(s)}-\bx)^\top \Delta_s^x\leq\frac{V_x(\bw^{(0)},\bx)}{2D_x^2}+\frac{D_x^2}{\alpha_x}\sum_{s=0}^{t}\gamma_s^2\|\Delta_s^x\|_{*,x}^2
	\leq\frac{1}{2}+\frac{D_x^2}{\alpha_x}\sum_{s=0}^{t}\gamma_s^2\|\Delta_s^x\|_{*,x}^2,\\
	\sum_{s=0}^{t}\gamma_s(\bv^{(s)}-\bx)^\top \Delta_s^x\leq\frac{V_x(\bv^{(0)},\bx)}{2D_x^2}+\frac{D_x^2}{\alpha_x}\sum_{s=0}^{t}\gamma_s^2\|\Delta_s^x\|_{*,x}^2 \leq\frac{1}{2}+\frac{D_x^2}{\alpha_x}\sum_{s=0}^{t}\gamma_s^2\|\Delta_s^x\|_{*,x}^2.
	\label{eq:lemma1withuforx}
\end{eqnarray*}
Writing $\bx - \bx^{(s)} = \bx- \bw^{(s)} + \bw^{(s)} - \bx^{(s)}$ and $\bx^{(s)} - \bx = \bv^{(s)} - \bx + \bx^{(s)} - \bv^{(s)}$, these two inequalities imply
\begin{eqnarray*}
	\label{eq:lemma1withvforxapplied}
	\sum_{s=0}^{t}\gamma_s(\bx-\bx^{(s)})^\top\Delta_s^x	\leq\sum_{s=0}^{t}\gamma_s(\bw^{(s)}-\bx^{(s)})^\top \Delta_s^x+\frac{1}{2}+\frac{D_x^2}{\alpha_x}\sum_{s=0}^{t}\gamma_s^2\|\Delta_s^x\|_{*,x}^2,\\
	\label{eq:lemma1withuforxapplied}
	\sum_{s=0}^{t}\gamma_s(\bx^{(s)}-\bx)^\top\Delta_s^x\leq\sum_{s=0}^{t}\gamma_s(\bx^{(s)}-\bv^{(s)})^\top \Delta_s^x+\frac{1}{2}+\frac{D_x^2}{\alpha_x}\sum_{s=0}^{t}\gamma_s^2\|\Delta_s^x\|_{*,x}^2.
\end{eqnarray*}
Hence,
\begin{align}
\left\vert\sum_{s=0}^{t}\gamma_s(\bx-\bx^{(s)})^\top\Delta_s^x\right\vert \leq& \max\left\{\left\vert\sum_{s=0}^{t}\gamma_s(\bw^{(s)}-\bx^{(s)})^\top \Delta_s^x\right\vert,~\left\vert\sum_{s=0}^{t}\gamma_s(\bx^{(s)}-\bv^{(s)})^\top \Delta_s^x\right\vert\right\} \nonumber \\&+\frac{1}{2}+\frac{D_x^2}{\alpha_x}\sum_{s=0}^{t}\gamma_s^2\|\Delta_s^x\|_{*,x}^2.\label{eq:boundonabs}
\end{align}
Applying~\eqref{eq:boundonabs} in \eqref{eq:gapllhat}, we get
\begin{eqnarray}
\nonumber
\left\vert\hat l^{(t)}_*-l^{(t)}_*\right\vert
&\leq&\max\left\{\left\vert\frac{\sum_{s=0}^{t}\gamma_s(\bw^{(s)}-\bx^{(s)})^\top \Delta_s^x}{\sum_{s=0}^{t}\gamma_s}\right\vert,~\left\vert\frac{\sum_{s=0}^{t}\gamma_s(\bx^{(s)}-\bv^{(s)})^\top \Delta_s^x}{\sum_{s=0}^{t}\gamma_s}\right\vert\right\}\\\label{eq:boundmaxx}
&&+\frac{0.5+(D_x^2/\alpha_x)\sum_{s=0}^{t}\gamma_s^2\|\Delta_s^x\|_{*,x}^2}{\sum_{s=0}^{t}\gamma_s}
+\left\vert\frac{\sum_{s=0}^{t}\gamma_s\delta_s}{\sum_{s=0}^{t}\gamma_s}\right\vert.
\end{eqnarray}
We next find a probabilistic bound for the right hand side of the above inequality.

\emph{Bounds on $\left\vert\dfrac{\sum_{s=0}^{t}\gamma_s(\bw^{(s)}-\bx^{(s)})^\top \Delta_s^x}{\sum_{s=0}^{t}\gamma_s}\right\vert$ and $\left\vert \dfrac{\sum_{s=0}^{t}\gamma_s(\bx^{(s)}-\bv^{(s)})^\top \Delta_s^x}{\sum_{s=0}^{t}\gamma_s}\right\vert$:} 
The inequality \eqref{eq:Vxbound} indicates that
\begin{align}
\nonumber
\|\bw^{(s)}-\bx^{(s)}\|_x&\leq \|\bx^{(s)}-\bx^{(0)}\|_x+\|\bw^{(s)}-\bx^{(0)}\|_x\leq \sqrt{\frac{2}{\alpha_x}V_x(\bx^{(0)},\bx^{(s)})}+\sqrt{\frac{2}{\alpha_x}V_x(\bx^{(0)},\bw^{(s)})}\\\label{eq:boundux}
&\leq\frac{2\sqrt{2}D_x}{\sqrt{\alpha_x}}.
\end{align} 
Define $\psi_s :=\nu_{s,t}(\bw^{(s)}-\bx^{(s)})^\top \Delta_s^x$ and  $\sigma_s :=\dfrac{4\sqrt{2}D_xM_x\nu_{s,t}}{\sqrt{\alpha_x}}$. Since $\bxi^{(s)}$ is independent of $\bw^{(s)}$ and $\bx^{(s)}$, we have
$\E_s[\psi_s]=0$. Furthermore,
\begin{equation}
\psi_s^2\leq\nu_{s,t}^2\left\|\bw^{(s)}-\bx^{(s)}\right\|_x^2\left\|\Delta_s^x\right\|_{*,x}^2\leq8\nu_{s,t}^2\left\|\Delta^x_s\right\|_{*,x}^2D_x^2/\alpha_x,\label{eq:boundonphi2}
\end{equation}
where the second inequality follows from \eqref{eq:boundux}. Using the definition of $\delta_s$, \eqref{eq:boundonphi2}, and \eqref{eq:lighttailgraddeltax}-\eqref{eq:lighttailgraddeltay}, it follows that
$\E_s[\exp(\psi_s^2/\sigma_s^2)]\leq \exp(1)$.
Hence Case A in Lemma~\ref{lem:sumterm} and union bound we get\looseness = -1
\begin{eqnarray}
\label{eq:boundcrosstermx3}
\text{Prob}\left\{\left|\sum_{s=0}^t\nu_{s,t}(\bw^{(s)}-\bx^{(s)})^\top \Delta_s^x\right|>\frac{4\sqrt{2}\Omega(\delta) D_xM_x}{\sqrt{\alpha_x}}\sqrt{\sum_{s=0}^t\nu_{s,t}^2} \right\}\leq 2\exp(-\Omega(\delta)^2/3).
\end{eqnarray}
With a similar argument, we can also show
\begin{eqnarray}
\label{eq:boundcrosstermx4}
\text{Prob}\left\{\left|\sum_{s=0}^t\nu_{s,t}(\bx^{(s)}-\bv^{(s)})^\top \Delta_s^x\right|>\frac{4\sqrt{2}\Omega(\delta) D_xM_x}{\sqrt{\alpha_x}}\sqrt{\sum_{s=0}^t\nu_{s,t}^2} \right\}\leq 2\exp(-\Omega(\delta)^2/3).
\end{eqnarray}

\emph{Bound on $\dfrac{\sum_{s=0}^{t}\gamma_s^2\|\Delta_s^x\|_{*,x}^2}{\sum_{s=0}^{t}\gamma_s}$:} Define $\psi_s :=\gamma_s\nu_{s,t}\|\Delta_s^x\|_{*,x}^2$ and  $\sigma_s :=4M_x^2\gamma_s\nu_{s,t}$. Using \eqref{eq:lighttailgraddeltax}-\eqref{eq:lighttailgraddeltay}, it is easy to verify that $\E_s\left[\exp\left(|\psi_s|/\sigma_s\right)\right]\leq \exp(1)$.
Hence, from Case B in Lemma~\ref{lem:sumterm} we have\looseness = -1

\small
\begin{align}
\text{Prob}\left\{\dfrac{\sum_{s=0}^{t}\gamma_s^2\|\Delta_s^x\|_{*,x}^2}{\sum_{s=0}^{t}\gamma_s} >
4M_x^2\sum_{s=0}^t\gamma_s\nu_{s,t}
+4\Omega(\delta) M_x^2\sqrt{\sum_{s=0}^t\gamma_s^2\nu_{s,t}^2}\right\}
\leq \exp(-\Omega(\delta)^2/12)+\exp(-3\Omega(\delta)/4).\label{eq:boundcrosstermx5}
\end{align}
\normalsize

\emph{Bound on $\left\vert\dfrac{\sum_{s=0}^{t}\gamma_s\delta_s}{\sum_{s=0}^{t}\gamma_s}\right\vert$:}
From definition of $\delta_s$ and \eqref{eq:lighttailobj}, it follows that
\begin{equation*}
\label{eq:lighttailobjtermx}
\E_s\left[\nu_{s,t}\delta_s\right] = 0\hspace{0.3cm}\text{and} \hspace{0.3cm}\E_s\left[\exp((\nu_{s,t}\delta_s)^2/(\nu_{s,t}Q)^2)\right]\leq\exp(1).
\end{equation*} 
Hence by Case A in Lemma~\ref{lem:sumterm} and union bound we get
\begin{eqnarray}
\label{eq:boundobjx}
\text{Prob}\left\{\left|\sum_{s=0}^t\nu_{s,t}\delta_s\right|>\Omega(\delta) Q\sqrt{\sum_{s=0}^t\nu_{s,t}^2} \right\}\leq 2\exp(-\Omega^2/3).
\end{eqnarray}	

The conclusion can be then obtained by upper bounding the right hand side of \eqref{eq:boundmaxx} using the union bound of
\eqref{eq:boundcrosstermx3}, \eqref{eq:boundcrosstermx4}, \eqref{eq:boundcrosstermx5}, and \eqref{eq:boundobjx}.
\hfill\halmos\endproof

\proof{\bf Proof of Proposition \ref{thm:luconvergecor}:}
(i) The definition of $\Omega(\delta)$ in \eqref{eq:defOmega} guarantees  $\exp(-\Omega(\delta)^2/3)+\exp(-\Omega(\delta)^2/12)+\exp(-3\Omega(\delta)/4)\leq \delta$.
Recall that $\nu_{s,t}=\dfrac{\gamma_s}{\sum_{s'=0}^{t}\gamma_{s'}}$.
With $\gamma_s=\dfrac{1}{M\sqrt{s+1}}$, it is straightforward to verify the following inequalities:

\begin{eqnarray}
\label{eq:nugamma1}
\sum_{s=0}^t\nu_{s,t}^2&=&\dfrac{\sum_{s=0}^t\frac{1}{s+1}}{\left(\sum_{s=0}^t\frac{1}{\sqrt{s+1}}\right)^2} \leq\dfrac{1+\ln(t+1)}{\left(2\sqrt{t+2}-2\right)^2},
\end{eqnarray}

\begin{eqnarray}
\label{eq:nugamma2}
\dfrac{\sum_{s=0}^t\gamma_s^2}{\sum_{s=0}^t\gamma_s}&=&\dfrac{1}{M}\cdot\dfrac{\sum_{s=0}^t\frac{1}{s+1}}{\sum_{s=0}^t\frac{1}{\sqrt{s+1}}} \leq \dfrac{1}{M}\cdot\dfrac{\left(1+\ln(t+1)\right)}{2\sqrt{t+2}-2},
\end{eqnarray}

\begin{eqnarray}
\label{eq:nugamma3}
\dfrac{1}{\sum_{s=0}^t\gamma_s}&=&\dfrac{M}{\sum_{s=0}^t\frac{1}{\sqrt{s+1}}} \leq\dfrac{M}{2\sqrt{t+2}-2},
\end{eqnarray}

\begin{eqnarray}
\label{eq:nugamma4}
\sum_{s=0}^t\gamma_s^2\nu_{s,t}^2&=& \left(\dfrac{1}{M}\right)^2\cdot\dfrac{\sum_{s=0}^t\frac{1}{(s+1)^2}}{\left(\sum_{s=0}^t\frac{1}{\sqrt{s+1}}\right)^2} \leq\dfrac{2\left(\frac{1}{M}\right)^2}{\left(2\sqrt{t+2}-2\right)^2}.
\end{eqnarray}
Applying these four inequalities to bound the terms in \eqref{eq:boundabsgaplu}, we get
\small
\begin{align}
\nonumber
& \text{Prob}\left\{u^{(t)}_*-l^{(t)}_*>4\sqrt{2}\Omega(\delta) M\sqrt{\sum_{s=0}^t\nu_{s,t}^2}+\frac{2+2.5M^2\sum_{s=0}^t\gamma_s^2}{\sum_{s=0}^t\gamma_s} +2.5\Omega(\delta) M^2\sqrt{\sum_{s=0}^t\gamma_s^2\nu_{s,t}^2}\right\}\\\nonumber
&\leq \text{Prob}\left\{u^{(t)}_*-l^{(t)}_*>\left(4\sqrt{2}\Omega(\delta) M+4.5 M+2.5\sqrt{2}\Omega(\delta) M\right)\frac{\left(1+\ln(t+1)\right)}{2\sqrt{t+2}-2}
\right\} \\\nonumber
&\leq \text{Prob}\left\{u^{(t)}_*-l^{(t)}_*>\left(10\Omega(\delta) M+4.5 M\right)\frac{\left(1+\ln(t+1)\right)}{2\sqrt{t+2}-2}
\right\} \\\label{eq:probsimplify}
&\leq  \delta.
\end{align}
Given $\epsoracle > 0$, let $\epsilon' := \epsoracle/\left(10\Omega(\delta) M+4.5 M\right)$. When $t\geq\max\left\{6,\left(\dfrac{8\ln(4/\epsilon')}{\epsilon'}\right)^2-2\right\}$, we have $	\dfrac{1+\ln(t+1)}{2\sqrt{t+2}-2}\leq\dfrac{2\ln(t+2)}{\sqrt{t+2}}$ and $\dfrac{2\ln(t+2)}{\sqrt{t+2}}$ is monotonically decreasing in $t$. Hence

\begin{equation}
\label{eq:tepsilon}
\dfrac{1+\ln(t+1)}{2\sqrt{t+2}-2}\leq\dfrac{2\ln(t+2)}{\sqrt{t+2}}
\leq\dfrac{\epsilon\ln((8/\epsilon')\ln(4/\epsilon'))}{4\ln(4/\epsilon')}
\leq\dfrac{\epsilon\ln(4/\epsilon')+\epsilon\ln(2\ln(4/\epsilon'))}{2\ln(4/\epsilon')}\leq \epsilon'.
\end{equation}
Using the above inequality in \eqref{eq:probsimplify} we get $\text{Prob}\left\{u^{(t)}_*-l^{(t)}_*> \left(10\Omega(\delta) M+4.5 M\right) \epsilon' = \epsoracle\right\}\leq \delta$ which completes the proof.\\[0.5EM]

(ii) We only prove this corollary for the lower bounds as the proof of upper bounds is similar.
The choice of $\Omega(\delta)$ guarantees  $6\exp(-\Omega(\delta)^2/3)+\exp(-\Omega(\delta)^2/12)+\exp(-3\Omega(\delta)/4)\leq \delta$.
Recall that $\nu_{s,t}=\dfrac{\gamma_s}{\sum_{s'=0}^{t}\gamma_{s'}}$.
Since $\gamma_s=\dfrac{1}{M\sqrt{s+1}}$, the inequalities \eqref{eq:nugamma1}, \eqref{eq:nugamma2}, \eqref{eq:nugamma3}, and \eqref{eq:nugamma4} hold. 
Applying these four inequalities to \eqref{eq:boundabsgapllhat} yields
\small
\begin{align}
\nonumber
&\text{Prob}\left\{\left|\hat l^{(t)}_*-l^{(t)}_*\right|>\left(\Omega(\delta) Q+8\Omega(\delta) M+2.5M\right)\frac{\left(1+\ln(t+1)\right)}{2\sqrt{t+2}-2}
\right\}\\\nonumber
&\leq\text{Prob}\left\{\left|\hat l^{(t)}_*-l^{(t)}_*\right|>\left(\Omega(\delta) Q+\frac{4\sqrt{2}\Omega(\delta) D_xM_x}{\sqrt{\alpha_x}}+0.5M+\frac{4D_x^2M_x^2}{M\alpha_x}+\frac{4\sqrt{2}\Omega(\delta) D_x^2 M_x^2}{M\alpha_x}\right)\frac{\left(1+\ln(t+1)\right)}{2\sqrt{t+2}-2}
\right\}\\\label{eq:boundabsgapllhatsimplify}
&\leq \delta,
\end{align}
\normalsize
where the first inequality follows from $M^2\geq\dfrac{2D_x^2M_x^2}{\alpha_x}$ by \eqref{eq:defM}.

Let $\epsilon' := \epsoracle/\left(\Omega(\delta) Q+8\Omega(\delta) M+2.5M\right)$. When $t\geq\max\left\{6,\left(\dfrac{8\ln(4/\epsilon')}{\epsilon'}\right)^2-2\right\}$, the inequality \eqref{eq:tepsilon} holds which can be applied to \eqref{eq:boundabsgapllhatsimplify} to show that $$\text{Prob}\left\{\left|\hat l^{(t)}_*-l^{(t)}_*\right|>\left(\Omega(\delta) Q+8\Omega(\delta) M+2.5M\right)\epsilon' = \epsoracle \right\}\leq \delta.$$
\hfill\halmos\endproof


\proof{\bf Proof of Theorem \ref{thm:totalcomplexity}.} We begin by establishing that the following inequalities hold with high probability in at most $T(\delta,\epsilon)$ number of iterations:
\begin{equation}\label{eq:conds}
u^{(t)}_*-l^{(t)}_*\leq \dfrac{1}{2}\epsoracle,\hspace{5pt} \left\vert\hat l^{(t)}_*-l^{(t)}_*\right\vert \leq \dfrac{1}{2}\epsoracle,\hspace{5pt} \text{and}\hspace{5pt} \left\vert\hat u^{(t)}_*-u^{(t)}_*\right\vert\leq \dfrac{1}{2}\epsoracle.
\end{equation} 
Given $\Omega(\delta)$, parts (i) and (ii) of Proposition \ref{thm:luconvergecor} imply that when
$$
t\geq\max\left\{6,\left(\frac{16\left(\Omega(\delta) Q+10\Omega(\delta) M+4.5M\right)}{\epsoracle}\ln\left(\frac{8\left(\Omega(\delta) Q+10\Omega(\delta) M+4.5M\right)}{\epsoracle}\right)\right)^2-2\right\},$$ we have 
$\text{Prob}\left\{u^{(t)}_*-l^{(t)}_*>\dfrac{\epsoracle}{2}\right\}\leq\delta/3$,  $\text{Prob}\left\{\left\vert\hat l^{(t)}_*-l^{(t)}_*\right\vert>\dfrac{\epsoracle}{2}\right\}\leq\delta/3$, and $\text{Prob}\Big\{\Big\vert\hat u^{(t)}_*-u^{(t)}_*\Big\vert>\dfrac{\epsoracle}{2}\Big\}\leq\delta/3$. 
Hence, using union bounds we get\looseness = -1
\begin{align*}
\text{Prob}\left\{u^{(t)}_*-l^{(t)}_*\leq\dfrac{\epsoracle}{2},\left\vert\hat l^{(t)}_*-l^{(t)}_*\right\vert\leq\dfrac{\epsoracle}{2},\left\vert\hat u^{(t)}_*-u^{(t)}_*\right\vert\leq\dfrac{\epsoracle}{2}\right\}\geq 1-\delta.
\end{align*}

To complete the proof we show that \eqref{eq:conds} implies $\Pc(r,\bar\bx^{(t)}) - H(r) \leq \epsoracle$ and $\left\vert\hat u_*^{(t)} - H(r)\right\vert \leq \epsoracle$. First note that we have
\begin{equation}\label{eq:relationshipbetweenpandH}
\Pc(r,\bar \bx^{(t)})\leq u^{(t)}_* \leq l^{(t)}_* + \dfrac{\epsoracle}{2}\leq  H(r) + \dfrac{\epsoracle}{2} \leq H(r) + \epsoracle,
\end{equation}
where the first inequality follows from~\eqref{eq:upperboundonHr}, the second from \eqref{eq:conds}, and the third holds since $l^{(t)}_*$ is a lower bound on $H(r)$. Using \eqref{eq:conds} and $u^{(t)}_* \leq H(r) + \dfrac{\epsoracle}{2}$, we get
\begin{equation}\label{upperboundonuhat}
\hat u^{(t)}_*\leq u^{(t)}_*+\dfrac{\epsoracle}{2}\leq H(r) + \dfrac{\epsoracle}{2} + \dfrac{\epsoracle}{2}= H(r) + \epsoracle. 
\end{equation}
In addition, 
\begin{equation}\label{lowerboundonuhat}
\hat u^{(t)}_*\geq u^{(t)}_*-\dfrac{\epsoracle}{2}\geq   H(r)- \dfrac{\epsoracle}{2}\geq H(r) - \epsoracle,
\end{equation}
where the first inequality holds by~\eqref{eq:conds} and the second since $u^{(t)}_*$ is an upper bound on $H(r)$. The inequalities \eqref{eq:relationshipbetweenpandH}-\eqref{lowerboundonuhat} complete the proof.\looseness = -1
\hfill\halmos\endproof

\proof{\bf Proof of Corollary \ref{thm:FullComplexityWithOVSMD}:} The proof directly follows from theorems \ref{generalcomplexity} and~\ref{thm:totalcomplexity} and definition of $\beta$.
 \hfill\halmos\endproof
 
Lemma~\ref{thm:totalcomplexity2} below shows the number of iterations required by Algorithm~\ref{alg:SMDwithOnlineValidationPDStop} to find the upper bound $\bar U$ on $H(r^{(0)})$.
 
 \begin{lemma} 
	\label{thm:totalcomplexity2}	
	Given an input tuple $(r^{(0)},\bar\alpha , \delta, \gamma_t,\theta)$, Algorithm~\ref{alg:SMDwithOnlineValidationPDStop} terminates with probability of at least $1-\delta$ after at most 
	\[\mathcal{O}\left(\ln\left(\frac{\theta}{(\theta-1)\beta}\right)\right)\]
	 OVSMD calls and  
	\[
	\mathcal{O}\left(\frac{1}{\beta^2}\left(\log_2\left(\frac{\bar\alpha }{\beta}\right)\right)^4\left(\log\left(\dfrac{1}{\delta}\right)\right)^2\right)
	\]
	gradient iterations. In addition, $H(r^{(0)})\leq \bar U<0 $ and ${|H(r^{(0)})|}/{|\bar U|}\leq \theta$ hold at termination.  
\end{lemma}
Note that we use the $\tilde{\mathcal{O}}$ complexity notation, which omits logarithmic terms, to simplify the expression for the gradient iteration complexity. 
\proof{Proof.}
We first prove that  Algorithm~\ref{alg:SMDwithOnlineValidationPDStop}  terminates with a probability of at least $1-\delta$.
Consider the $h$th iteration of this algorithm. Given $
\hat u_*^{(h)}$ returned by OVSMD, Theorem~\ref{thm:totalcomplexity} guarantees with a probability of at least $1-\delta^{(h)}$ that
\begin{eqnarray}
\label{eq:chainineq_U}
\hat  u_*^{(h)}-\alpha^{(h)}\leq  H(r^{(0)})\leq \hat u_*^{(h)}+\alpha^{(h)}.
\end{eqnarray} 
Since $\sum_{h=0}^\infty\delta^{(h)}=\delta$, using union bound it is clear that \eqref{eq:chainineq_U} holds for $h=0,1,2,\dots,$ with a probability of at least $1-\delta$. 
In addition, \eqref{eq:chainineq_U} implies that $\hat u_*^{(h)}+\alpha^{(h)}\leq H(r^{(0)})+2\alpha^{(h)}\leq 0$ when $\alpha^{(h)}\leq-H(r^{(0)})/2$.
Furthermore, when $\alpha^{(h)}\leq-\frac{\theta-1}{2\theta}  H(r^{(0)})$ (which also indicates that  $\alpha^{(h)}\leq-\frac{H(r^{(0)})}{2}$ since $\theta>1$ and $H(r^{(0)})\leq 0$),  we have 
\begin{eqnarray}
	\label{eq:chainineq_UU}
	\frac{\hat u_*^{(h)}-\alpha^{(h)}}{\hat u_*^{(h)}+\alpha^{(h)}}=  \frac{-\hat u_*^{(h)}+\alpha^{(h)}}{-\hat u_*^{(h)}-\alpha^{(h)}} =
\frac{-\hat u_*^{(h)}+\alpha^{(h)}}{-\hat u_*^{(h)}+\alpha^{(h)}-2\alpha^{(h)}}\leq\frac{-H(r^{(0)})}{-H(r^{(0)})-2\alpha^{(h)}} \leq\frac{1}{1-\frac{\theta-1}{\theta}}=\theta,
\end{eqnarray} 
where the first inequality follows from the inequality $-H(r^{(0)})\leq -\hat u_*^{(h)} + \alpha^{(h)}$ and the fact that the function $x/(x-2\alpha^{(h)})$ is a decreasing function in $x$.
\eqref{eq:chainineq_UU} indicates that as soon as $\alpha^{(h)}\leq-\frac{\theta-1}{2\theta}  H(r^{(0)})$, the stopping criteria of Algorithm~\ref{alg:SMDwithOnlineValidationPDStop} hold and the algorithm terminates with a probability of $1-\delta$.  Since $\alpha^{(h)} = \alpha^{(0)}/2^h=\bar\alpha/2^h$ and $\beta = \Omega(|H(r^{(0)})|)$, the inequality $\alpha^{(h)}\leq-\frac{\theta-1}{2\theta}  H(r^{(0)})$ can be guaranteed in at most $
J:=\log_2\left(\dfrac{2\theta\bar\alpha}{(\theta-1)|H(r^{(0)})|}\right)
= \mathcal{O}\left(\ln\left(\frac{\theta}{(\theta-1)\beta}\right)\right)$ iterations. Furthermore, the inequalities \eqref{eq:chainineq_U} and \eqref{eq:chainineq_UU} imply that at terminatopm $\bar U= \hat u_*^{(J)}+\alpha^{(J)}<0$ and $\dfrac{|H(r^{(0)})|}{|\bar U|}\leq\dfrac{\hat u_*^{(J)}-\alpha^{(J)}}{\hat u_*^{(J)}+\alpha^{(J)}}\leq\theta $.

We next compute the total number of gradient iterations taken by Algorithm~\ref{alg:SMDwithOnlineValidationPDStop}. Notice that by Theorem~\ref{thm:totalcomplexity}, the $h$-th call of OVSMD requires at most $T(\delta^{(h)},\alpha^{(h)})$ iterations. 
Therefore, the total number of iterations can be computed as 
\small
\begin{align*}
&\sum_{h=0}^J T\left(\delta^{(h)},\alpha^{(h)}\right)\\
=&\sum_{h=0}^J\max\left\{6,\left(\frac{16\left(\Omega(\delta^{(h)}) Q+10\Omega(\delta^{(h)}) M+4.5M\right)}{\alpha^{(h)}}\ln\left(\frac{8\left(\Omega(\delta^{(h)}) Q+10\Omega(\delta^{(h)}) M+4.5M\right)}{\alpha^{(h)}}\right)\right)^2-2\right\}\\
=&\sum_{h=0}^J\mathcal{O}\left(\frac{\Omega(\delta^{(h)})}{\alpha^{(h)}}\ln\left(\frac{\Omega(\delta^{(h)})}{\alpha^{(h)}}\right)\right)^2
= \sum_{h=0}^J\mathcal{O}\left(\frac{h2^h}{\bar\alpha}\log\left(\dfrac{1}{\delta}\right)\ln\left(\frac{h2^h}{\bar\alpha}\log\left(\dfrac{1}{\delta}\right)\right)\right)^2
=\sum_{h=0}^J\mathcal{O}\left(\frac{h^22^h}{\bar\alpha}\log\left(\dfrac{1}{\delta}\right)\right)^2\\
=  &\sum_{h=0}^J\mathcal{O}\left(\frac{J^42^{2h}}{\bar\alpha^2}\log^2\left(\dfrac{1}{\delta}\right)\right)
= \mathcal{O}\left(\frac{J^4\log^2(1/\delta)}{|H(r^{(0)})|^2}\right) = \mathcal{O}\left(\frac{\log^2(1/\delta)}{|H(r^{(0)})|^2}\left(\log_2\left(\frac{\bar\alpha}{|H(r^{(0)})|}\right)\right)^4\right)\\
=& \mathcal{O}\left(\frac{1}{\beta^2}\left(\log_2\left(\frac{\bar\alpha}{\beta}\right)\right)^4\left(\log\left(\dfrac{1}{\delta}\right)\right)^2\right),
\end{align*} 
where we used $\Omega(\delta^{(h)}) = \mathcal{O}\left(h\log\left(\dfrac{1}{\delta}\right)\right)$ and $\alpha^{(h)} = \frac{\bar\alpha}{2^h}$ in the second inequality, $J=\mathcal{O}\left(\log_2\left(\dfrac{\bar\alpha}{|H(r^{(0)})|}\right)\right)$ in the third and fourth equations, and $|H(r^{(0)})| = \Theta\left(\beta\right)$.  
\hfill\halmos\endproof
\smallskip
\proof{\bf Proof of Theorem \ref{cor:FullComplexityWithOVSMD_U}:} The proof of this theorem is a direct result of Corollary~\ref{thm:FullComplexityWithOVSMD} and Lemma~\ref{thm:totalcomplexity2}. In particular, it is straightforward to see that the total number of OVSMD calls is \[\mathcal{O}\left(\ln\left(\frac{\theta}{(\theta-1)\beta}\right)\right) +  \mathcal{O}\left(\dfrac{\theta^2}{\beta}\ln\left(\dfrac{\theta^2}{\beta\epsilon}\right)\right)= \mathcal{O}\left(\dfrac{\theta^2}{\beta}\ln\left(\dfrac{\theta^2}{(1-\theta)\beta\epsilon}\right)\right).\]
In addition, the total number of gradient iterations can be computed as 
\begin{align*}
&{\mathcal{O}}\left(\dfrac{1}{\beta^2}\ln^4\left(\frac{\bar\alpha}{\beta}\right)\ln^2\left(\dfrac{1}{\delta}\right)\right) + \mathcal{O}\left(\dfrac{\theta^2}{\beta\epsilon^2}\cdot \ln\left(\dfrac{\theta^2}{\beta\epsilon}\right)\cdot\ln^2\left(\dfrac{1}{\delta}\right)\cdot\ln^2\left(\dfrac{1}{\epsilon}\right)\right) 
\end{align*}
\hfill\halmos\endproof

\smallskip
\proof{\bf Proof of Corollary \ref{cor:implStop}:} 
Let $\delta^{(k)}=\dfrac{\delta}{2^k}$ for $k\geq0$ as defined in SFLS. With a little abuse of notation, we use $\Omega(n)$ to represent a quantity whose order of magnitude is at least $n$.
According to Theorem~\ref{thm:totalcomplexity}, for any $\delta\in(0,1)$ and $K\geq0$, there exists $\epsoracle$ satisfying $\Omega\left(\frac{\ln(1/\delta^K)}{\sqrt{T}}\right)\leq \epsoracle\leq \mathcal{O}\left(\frac{\ln(1/\delta^K)\ln(T)}{\sqrt{T}}\right)$ such that OVSMD is a valid stochastic oracle $\mathcal{A}\left(r^{(k)},\epsoracle,\delta^{(k)}\right)$ for iteration $k=0,1,\dots,K$ of SFLS. Let $\epsilon=-\frac{2\theta^2(\theta+1)}{(\theta-1)H(r^{(0)})}  \epsoracle$ such that $\Omega\left(\frac{K\theta^2\ln(1/\delta)}{\sqrt{T}}\right)\leq\epsilon\leq \mathcal{O}\left(\frac{K\theta^2\ln(1/\delta)\ln(T)}{\sqrt{T}}\right)$. Hence, there exists $K=\mathcal{O}\left(\frac{\theta^2}{\beta}\ln\left(\frac{T}{\beta}\right)\right)$
such that $K\geq\frac{2\theta^2}{\beta}\ln\left(\dfrac{\theta^2}{\beta\epsilon}\right)$. With such $K$ and $\epsilon$, according to Theorem~\ref{generalcomplexity}, SFLS generates a feasible solution at iteration $k=0,1,\dots,K$ and finds a relative $\epsilon$-optimal and feasible solution with $\epsilon\leq \mathcal{O}\left(\frac{\theta^4\ln(1/\delta)\ln(T)\ln\left(T/\beta\right)}{\beta\sqrt{T}}\right)$ with a probability of at least $1-\delta$ in at most $K$ outer iterations (calls of OVSMD), which corresponds to $KT=\mathcal{O}\left(\frac{\theta^2T}{\beta}\ln\left(\frac{T}{\beta}\right)\right)$ gradient iterations.
\hfill\halmos\endproof
\end{document}